\documentclass[a4paper,11pt]{amsart}\usepackage{amsfonts,amssymb,bbm}
\usepackage{enumerate,mathrsfs,hyperref}

\UseRawInputEncoding

\def\empty{\varnothing}

\def\A{\mathbb{A}} \def\C{\mathbb{C}}
\def\k{\Bbbk}
\def\N{\mathbb{N}}\def\Q{\mathbb{Q}}
\def\R{\mathbb{R}}\def\Z{\mathbb{Z}}
\def\bk{{\bar{\k}}}\def\bN{{\bar{\N}}}

\def\di{\partial}

\def\bl{\langle}\def\br{\rangle}

\def\liml{\lim\limits}
\def\suml{\sum\limits}

\def\prodl{\mathop\prod\limits}

\def\RpX{R^{\oplus p}_X}\def\RpY{R^{\oplus p}_Y}

\newcommand{\quot}[2]{{\footnotesize\left.\raisebox{1.2ex}{$#1$}\!\! \ensuremath\diagup \!\!\raisebox{-1.2ex}{$#2$}\right.}}
\newcommand{\quots}[2]{{\footnotesize\left.\raisebox{0.4ex}{$#1$}\! / \!\raisebox{-0.4ex}{$#2$}\right.}}

\renewcommand{\stackrel}[2]{\ \lower 0.4ex \hbox{$\mathrel{\mathop{#2}\limits^{#1}}$}\ }

\def\tC{\tilde{C}}\def\tc{\tilde{c}}\def\tf{{\tilde{f}}}
\def\tl{{\tilde{l}}} 
\def\tq{\tilde{q}}\def\tQ{\tilde{Q}}

\def\tU{{\tilde{U}}}
\def\ty{{\tilde{y}}}\def\tz{{\tilde{z}}}

\def\hR{{\widehat{R}}}
\def\hy{\hat{y}}\def\hz{\hat{z}}

\def\De{\Delta}

\def\la{\lambda}

\def\cA{\mathscr A}\def\ca{\mathfrak a}
\def\cb{\mathfrak b}\def\cC{\mathscr C}\def\cc{{\frak c}}

\def\cG{\mathscr G}\def\ch{\mathcal H}
\def\cK{{\mathscr K}\!}\def\cL{\mathscr L}\def\cO{\mathcal O}\def\cR{\mathscr{R}}
\def\cS{{\mathcal S}}
\def\cm{{\frak m}}

\def\um{{\underline{m}}}
\def\uy{{\underline{y}}}

\def\one{{1\hspace{-0.1cm}\rm I}}

\newcommand{\ber}{\begin{array}{l}}\newcommand{\eer}{\end{array}}
\newcommand{\bpm}{\begin{pmatrix}}\newcommand{\epm}{\end{pmatrix}}
\newcommand{\bbm}{\begin{bmatrix}}\newcommand{\ebm}{\end{bmatrix}}
\newcommand{\bM}{\begin{matrix}}\newcommand{\eM}{\end{matrix}}
\newcommand{\bee}{\begin{enumerate}}\newcommand{\eee}{\end{enumerate}}
\newcommand{\bei}{\begin{itemize}}\newcommand{\eei}{\end{itemize}}

\def\li{~\\ $\bullet$ }

\def\wrt{with respect to }\def\iff{if and only if }
\def\sset{\subset}\def\sseteq{\subseteq}\def\ssetneq{\subsetneq}\def\smin{\setminus}

\def\Maps{Maps(X,\!Y)}
\def\Mapsk{\rm{Maps}\big((\k^n,o),\!(\k^p,o)\big)}

\newtheorem{Lemma}{Lemma}[section]\newcommand{\bel}{\begin{Lemma}}\newcommand{\eel}{\end{Lemma}}
\newtheorem{Theorem}[Lemma]{Theorem}\newcommand{\bthe}{\begin{Theorem}}\newcommand{\ethe}{\end{Theorem}}
\newtheorem{Proposition}[Lemma]{Proposition}\newcommand{\bprop}{\begin{Proposition}}\newcommand{\eprop}{\end{Proposition}}
\newtheorem{Corollary}[Lemma]{Corollary}\newcommand{\bcor}{\begin{Corollary}}\newcommand{\ecor}{\end{Corollary}}
\newtheorem{Definition}[Lemma]{Definition}\newcommand{\bed}{\begin{Definition}}\newcommand{\eed}{\end{Definition}}
\newtheorem{Definition-Proposition}[Lemma]{Definition-Proposition}

\def\bpr{~\\{\em Proof.\ }}
\newcommand{\epr}{{\hfill\ensuremath\blacksquare}}
\newtheorem{Remark}[Lemma]{Remark}\newcommand{\beR}{\begin{Remark}\rm}\newcommand{\eeR}{\end{Remark}}
\newtheorem{Example}[Lemma]{Example}\newcommand{\bex}{\begin{Example}\rm}\newcommand{\eex}{\end{Example}}
\newtheorem{Problem}[Lemma]{Problem}\newcommand{\bprob}{\begin{Problem}\rm}\newcommand{\eprob}{\end{Problem}}
\newtheorem{Properties}[Lemma]{Properties}
\newtheorem{Assumptions}[Lemma]{Assumptions}

\newcommand{\bet}{\begin{tabular}{cccccccc}}\newcommand{\eet}{\end{tabular}}
\newcommand{\beq}{\begin{equation}}\newcommand{\eeq}{\end{equation}}

\def\into{{\hookrightarrow}}

\title[]{O\MakeLowercase{rbits of the left-right equivalence of maps in arbitrary characteristic}}
\author[]{D\MakeLowercase{mitry} K\MakeLowercase{erner}}
  \address{Department of Mathematics, Ben Gurion University of the Negev, P.O.B. 653, Be'er Sheva 84105, Israel.}
 \email{dmitry.kerner@gmail.com}
\date{\today.\ \  Filename: \jobname.tex}
\thanks{I was supported by the Israel Science Foundation,  grants No.  1910/18 and 1405/22}


\begin{document}
\begin{abstract}
The germs of maps $(\k^n,o)\stackrel{f}{\to}(\k^p,o)$ are traditionally studied up to the right ($\cR$), left-right ($\cA$) or contact ($\cK$) equivalence.
Various questions about the group-orbits $\cG f$ (for $\cG$ one of $\cR,\cK,\cA$) are reduced to their tangent spaces, $T_\cG f$. Classically the  passage $T_\cG f\rightsquigarrow \cG f$ was done by
 vector fields integration, hence it was bound to the $\R/\C$-analytic or $C^r$-category.

The purely-algebraic (characteristic-free) approach to the group-orbits of $\cR,\cK$ has been developed during the last decades.
But those methods could not address the (essentially more complicated) $\cA$-equivalence.
   Moreover, the characteristic-free results (for $\cR,\cK$) were weaker than those in characteristic zero,
   because of the (inevitable) pathologies of positive characteristic.

In this paper we close these omissions.
\bei
\item
 We establish the general (characteristic-free) passage $T_\cG f\rightsquigarrow \cG f$ for the groups $\cR$, $\cK,$ $\cA$.
  Submodules of $T_\cG f$ ensure (shifted) submodules of $\cG f$.
    For the $\cA$-equivalence this extends (and strengthens) various classical results of Mather, Gaffney, du Plessis,  and others.

\item Given a filtration on the space of maps one has the filtration on the group, $\cG^{(\bullet)}$, and on the tangent space,  $T_{\cG^{(\bullet)}}$.
 We establish the criteria of type ``$T_{\cG^{(j)}}f$ vs $\cG^{(j)} f$" in their strongest form,
 for arbitrary base field/ring, provided the characteristic is zero or high  for a given $f$.
  This brings the ``inevitably weaker" results of $char>0$ to the level of $char=0.$
\eei

\medskip

As an auxiliary step, important on its own, we develop the mixed-module structure of the tangent space $T_{\cA}f$ and establish various  properties of the annihilator ideal $\ca_\cA$, defining the instability locus of the map.
\end{abstract}
\maketitle
\setcounter{tocdepth}{1}
\tableofcontents

\section{Introduction}
\subsection{}
Consider $\k$-analytic map-germs, $f: (\k^n,o) {\to}(\k^p,o)$, here $\k$ is $\R$ or $\C.$
 Whitney's initial results on stable maps were greatly developed and extended by Thom, Mather and many others.
 (See e.g. \cite{AGLV}, \cite{Martinet.1982}, \cite{Mond-Nuno}.)
The maps are traditionally studied up to
 the  left-right ($\cA$) equivalence, with the auxiliary  right ($\cR$) and contact ($\cK$) equivalences. The first questions are about the group-orbits, $\cG f\sset \Mapsk$,
  in particular the finite determinacy.

  The classical methods relied heavily on vector fields integration. This chained the whole theory to the $\R,\C$-analytic (or $C^r$-differentiable) setting. Even the Nash function-germs, $\R\bl x\br ,$ could not be treated.

The study of $\cR,\cK$-equivalences by purely algebraic methods  (in zero and positive characteristic) has begun in 80's, e.g. \cite{Gre.Kron.90}.
 This was motivated  by the study of more complicated (e.g. non-isolated) singularities, by the applications in Algebraic Geometry/Commutative Algebra, and by the need to enable calculations in Computer Algebra.
 By now numerous results are available, \cite{Boubakri.Gre.Mark}, \cite{Boubakri.Gre.Mark2011},
  \cite{Greuel-Nguyen.2016}, \cite{Nguyen}, \cite{Greuel-Pham.2017}, \cite{B.K.motor}, \cite{BGK.20}.
 The results in positive characteristic are (non-surprisingly) weaker than those in characteristic zero.

The $\cA$-case (the most important for the study of maps) is essentially more complicated than the (auxiliary) $\cR,\cK$-cases.
 Even the case of characteristic zero was untouched.
 The $\cA$-orbits could not  be studied by the (algebraic) methods of the previous papers, e.g. because of the following difficulties.
 \bei
 \item The group-action $\cA\circlearrowright \Mapsk$ is neither additive, nor $\k$-multiplicative,
 unlike the actions $\cR,\cK^{lin}\circlearrowright \Mapsk$.

\item  The tangent space $T_\cA f$ is    not an $\cO_{(\k^n,o)}$-module, it is only a ``mixed  $(\cO_{(\k^n,o)},\cO_{(\k^p,o)})$-module".
\item The essential ingredient was the Thom-Levine lemma (``naturality of vector field integration") on the solution of a certain differential equation. Its algebraic version was absent.
 \item The classical Artin approximation (so helpful for $\cR,\cK$-equivalences) was absent (until recently) for the left-right equivalence.
 (See \S\ref{Sec.Preliminaries.Basic.Tools}.v for more detail.)
\eei

\subsection{The results} (The detailed description is in \S\ref{Sec.Intro.Contents})
 We study maps of germs, $Maps(X,(\k^p,o)),$  of arbitrary characteristic. Here $X=Spec(R_X)$ and $(\k^p,o)=Spec(R_Y)$ are formal/$\k$-analytic/$\k$-Nash germs of schemes. More precisely, $(R_X,R_Y)$ is one
  of the pairs     $(\quots{\k[[x]]}{J},\k[[y]])$, $(\quots{\k\{x\}}{J},\k\{y\})$, $(\quots{\k\bl x\br}{J},\k\bl y\br),$
   see \S\ref{Sec.Preliminaries.Rings.Germs}.ii. Here $\k$ is any field or an excellent Henselian ring. (The later case is important for deformations.)

The space of maps is an $R_X$-module, $Maps(X,(\k^p,o)):=Hom(R_{(\k^p,o)},R_X)\cong \cm\cdot \RpX$, for  the maximal ideal $\cm\sset R_X.$
 Fix a map $f\in \cm\cdot\RpX$ and let $\cG$ be one of the groups $\cR,\cK,\cA.$ Our results are of two types: conditions to ensure
  the large orbit $\cG f$ (in terms of the tangent space $T_\cG f$), and the structure of the tangent space $T_\cA f$.

 \subsubsection{The $\cG$-implicit function theorems} Fix   an ideal $\ca\sset R_X.$ Theorems \ref{Thm.IFT.R.case}, \ref{Thm.IFT.K.case}, \ref{Thm.IFT.A.case} read (roughly):
  \beq\label{Eq.1}
\text{If $\ca\cdot \RpX\sseteq T_\cG f$ \quad     then \quad $\cG f\supseteq \{f\}+\ca^2\cdot \RpX.$}
  \eeq
 This is characteristic-free, 
  and with no restrictions on the critical/singular/instability locus of the map $f.$

In many case $\ca$ can be chosen as an ideal that defines the non-reduced critical/singular/instability locus (with the natural scheme structure).
One can say roughly: $f$ is $\cG$-determined by its ``2-jet" on the critical/singular/instability locus. (The ``2-jet" means the image in $\quots{R_X}{\ca^2}$) Some results of this type are known for $\ca=\cm^d$ and $R_X=\C\{x\}$, e.g.
 \cite[Theorem 6.2, pg. 182]{Mond-Nuno}.
 The proofs were heavily analytic. Even the characteristic zero case is totally new.

 \subsubsection{The filtration criteria, ``$T_{\cG^{(j)}} f$ vs $ \cG^{(j)} f$"}
 Fix an ideal $I\sset R_X$, then the space of maps acquires the filtration $I^\bullet\cdot \RpX$. This induces filtrations on the group, $\cG^\bullet$, and on the tangent space, $T_{\cG^\bullet}$. Theorems \ref{Thm.Orbits.R.Filtration}, \ref{Thm.Orbits.K.Filtration}, \ref{Thm.Orbits.Aj.vs.TAj} read (roughly):
\beq\label{Eq.2}
\text{ $ \cG^{(j)} f \supseteq \{f\}+I^d\cdot  \RpX$   \quad    if and only if   \quad $ T_{\cG^{(j)}} f \supseteq I^d\cdot  \RpX.$}
 \eeq
   In this case we assume: either $char(\k)=0$ or $char(\k)>\frac{d-ord(f)}{j}.$

For $\cR,\cK$-equivalences this result is known in characteristic zero, \cite{BGK.20}, while for  $char(\k)>\frac{d-ord(f)}{j}$ it  strengthens the known results ``by a factor of 2". We remark, that the previous proofs of the ``only if"-part (in $char>0$) needed an additional assumption   ``$\k$ is infinite" and  non-trivial deformation/specialization arguments, \cite{Greuel-Pham.2017}.

For $\cA$-equivalence this was known only in the analytic/smooth-cases, $R_X=\C\{x\},\R\{x\},C^\infty(\R^n,o),$ with heavily analytic proofs, \cite{Bruce.du-Plessis.Wall}. Even the zero characteristic case is new. 

A remark, the ideals $\ca,I$ in \eqref{Eq.1} and \eqref{Eq.2} are not necessarily $\cm$-primary. I.e. there is no restriction on the dimension of the critical/singular/instability loci.

\medskip
 For many more results see \S\ref{Sec.Intro.Contents}.

  \subsubsection{The structure of $T_\cA f$} The tangent spaces $T_\cR f,T_\cK f\sseteq \RpX$ are $R_X$-submodules. Their properties can be addressed via the standard commutative/homological algebra. The tangent space $T_\cA f$ is not an $R_X$-module, it is only an $R_Y$-module
   (and  not finitely-generated when $dim(X)>p$ or when $f$ is not $\cA$-finite). The classical Nakayama/Artin-Rees lemmas do not hold for $T_\cA f.$
  We establish their substitutions  and numerous ``little tools" in \S\ref{Sec.TA}. (See \S\ref{Sec.Intro.Contents} for the list of results.)

 In particular we prove: if $f$ is a morphism of ``finite singularity type", then the image tangent space $T_\cA f$ is a {\em finitely generated} module over certain extension of $R_{(\k^p,o)},$ namely: the Noetherian, local  ring $f^\sharp R_{(\k^p,o)}+\ca_\cA.$

\subsection{The methods} The proofs of theorems \ref{Thm.IFT.R.case}, \ref{Thm.IFT.K.case} are heavily based on the implicit function arguments (hence their name ``$\cG$-implicit function theorem").
More precisely, we use a version of Newton lemma  of the type of \cite{Bourbaki.CA}-\cite{Tougeron1966}.
The proof of (the main) theorem  \ref{Thm.IFT.A.case} is much more complicated and demands additional tools. (Weierstra\ss\ division, arc-solution, Kostant-Rosenlicht theorem, and a special argument when the field
 $\quots{\k}{\cm_\k}$ is not algebraically closed.)

The proofs of theorems \ref{Thm.Orbits.R.Filtration},\ref{Thm.Orbits.K.Filtration},\ref{Thm.Orbits.Aj.vs.TAj} are based on the {\em algebraic}  Thom-Levine lemma (\S\ref{Sec.3.Related.Derivations.Autos}) and the {\em full}
 Baker-Campbell-Hausdorff expansion.

\subsection{}
Establishing these ``fundamentals" of $\cA$-equivalence opens the way to further results.
 In this paper I give only the basic/simplest examples, see \S\ref{Sec.Intro.Contents}.
  More serious applications contain (in arbitrary characteristic) the theory of unfoldings, the theory of stable maps,     Mather-Yau/Gaffney-Hauser theorems,  algebraization of (non-finitely-determined) maps, \cite{Kerner.Unfoldings}.
 Moreover, the new criteria on the orbit $\cA f$ and on the structure of the tangent space $T_\cA f$ seem useful   for the   classification/study of $\cA$-simple ($\C$-analytic) maps.

A remark: the proofs do not use the notion of ``geometric equivalence" of \cite{Damon84}, the methods seem applicable to the orbit-study of various other groups, e.g. volume preserving/symplectic equivalences, \cite{Domitrz-Rieger}.

\subsection{Acknowledgement} Thanks are to A. F. Boix,  G.-M. Greuel, D. Mond, M. A. S. Ruas for important advices.
 I was heavily influenced by their papers and also by those of T. Gaffney, A. du Plessis, L. Wilson.

 Special thanks are to M. Borovoi and Z. Rosegarten for their help with Kostant-Rosenlicht issues, \S\ref{Sec.Preliminaries.Orbits.Unipotent.Group.Closed}.


Finally I thank the  referee of \cite{BGK.20} who urged me  to treat the $\cA$-case (and not only $\cR,\cK$).

\subsection{The detailed structure, results, and contents of the paper}\label{Sec.Intro.Contents}

\bee[\!\S 1\!]\setcounter{enumi}{1}
\item is about the general preliminaries, to make the paper self-contained.
\bee[\!\S 2.1\!]
\item 
 fixes the rings we work with, i.e. $\quots{\k[\![x]\!]}{J},$  $\quots{\k\{x\}}{J},$  $\quots{\k\bl x\br}{J}.$
 \item  identifies  the space of maps, $Maps(X,(\k^p,o))= \cm\cdot \RpX$, the isomorphism of $R_X$-modules.

\item recalls ``the basic tools",
  Weierstra\ss\ division and finiteness,
 the full Baker-Campbell-Hausdorff formula, Artin approximation for $\k$-analytic/$\k$-Nash equations and for the left-right equivalence.

\item 
is about the ``local coordinate changes"  (i.e. the ring automorphisms $Aut_X:=Aut_\k(R_X)$) and
 the ``germs of vector fields" (i.e. the derivations  $Der_X:=Der_\k(R_X)$). We recall the corresponding filtrations, on the group, $\{Aut^\bullet_X\}$,
  and on the module, $\{Der^\bullet_X\}$.

\item 
 recalls the    exp/log maps, $Der_X\leftrightarrows Aut_X$ and their truncated versions, $jet_N Exp,$ $jet_N Ln$.
 In zero or high characteristic they substitute the  integration of vector fields.

\eee

\item sets-up the basic notions of $\cR,\cK,\cA$ equivalence over $\k.$ The cases $\cR,\cK$ are ``known in some sense",
   \cite{B.K.motor}, \cite{BGK.20}, but the $\cA$-case is new.
\bee[\!\S3.1\!]

\item The   group-actions on the space of maps, $\cG\circlearrowright Maps(X,(\k^p,o))= \cm\cdot \RpX$, for $\cG$ one of the groups $\cR,\cL,\cC,\cK,\cA:=\cL\times\cR$, are defined as in the classical case.

\item The   (extended) image tangent spaces $T_\cG f\sset \RpX$ are defined via the module(s) of derivations, e.g. 
 $T_\cR f:=Der_X(f),$ $T_\cK f:=T_\cR f+(f)\cdot \RpX.$ 
 As in the classical case, $T_\cR f,$ $T_\cK f$ are $R_X$-submodules of $\RpX$, while $T_\cA f$ is only a ``mixed" module.

\item Fix an ideal $I\sset R_X$ to get the filtration on the space of maps, $I^\bullet\cdot \RpX.$ It induces   (natural) filtrations of the group, $\cG\rhd \cG^{(\bullet)}$, and of the tangent space,
 $T_\cG\supseteq T_{\cG^{(\bullet)}}.$ 

A little twist: the filtrations $\cA^\bullet$, $T_{\cA^\bullet}$ and $\cK^\bullet$, $T_{\cK^\bullet}$ could be also defined in another way.
 Lemma \ref{Thm.TA.TK.distinct.filtrations} ensures: these distinct definitions are compatible
  when $\k$ is an infinite field. 
\item
This filtration $I^\bullet\cdot \RpX$ defines the (filtration) topology on $\RpX$. 
 The tangent space $T_\cG f$, and the  orbit  $\cG f$,  are closed in this topology for the groups $\cG= \cR,\cK$.
  The (standard) proof is based on the Artin-Rees lemma and the classical Artin approximation.

 The $\cA$-case is more complicated (in the absence of standard methods), see \S\ref{Sec.TA.filtration.closure}.

 \item 
gives  the algebraic version of the classical Thom-Levine  lemma ``Naturality of vector fields integration".
 That basic lemma was of key importance for $\cA$-equivalence in the $\R,\C$-analytic cases.
 And the algebraic version is used through the paper (when the characteristic is zero or high enough) in the same way.

  \item 
  introduces the annihilator ideal $\ca_\cG:=Ann[T^1_\cG f]\sset R_X$. Here $T^1_\cG f:=\quots{\RpX}{T_\cG f}$,
   this is the tangent space to $\cG$-miniversal unfoldings, when $f$ is $\cG$-finite.
 For functions ($p=1$) on smooth germs, i.e. $(\k^n,o)\to (\k^1,o)$, these are the classical Jacobian/Tjurina ideals, $\ca_\cR=Jac(f)$, $\ca_\cK=Jac(f)+(f)$.
 When the source $X$ is singular, one gets   Bruce-Roberts' versions of these ideals.

 More generally (for $p\ge1$), the ideal $\ca_\cR$ defines the critical locus of a map, $\ca_\cK$ defines the singular locus, $\ca_\cA$ defines
  the instability locus.
  A remark: $\ca_\cG$ is defined as the annihilator ideal, not as the Fitting ideal.

The annihilators $\ca_\cR,\ca_\cK$ are well known, but $\ca_\cA$ seems to be not studied previously.

    Lemma \ref{Thm.K.finite.maps} extends  the classical fact (from $\k=\C$ to local Henselian rings): the restriction $f|_{Crit(f)}:Crit(f)\to (\k^p,o)$ is a
   finite morphism iff $f$ is $\cK$-finite iff $T^1_\cR f$ is a f.g. module over $R_Y$.
\eee

\item   is the first  warmup, addressing the orbits of the group $\cR.$
\bee[\!\!\S4.1\!]
\item  contains the ``{\em$\cR$-implicit function theorem}". It puts a large part of the tangent space
 into the group-orbit: $\cR f\supseteq\{f\}+\ca\cdot T_\cR f$. This statement (and its proof) is characteristic-free.

 The proof goes in the style of implicit function theorem of \cite{Bourbaki.CA} and \cite{Tougeron1966}.

\item
 Numerous examples show that this bound is much stronger than the known bounds (even in the $\R,\C$-cases).
 As a trivial corollary we get the Morse lemma over an arbitrary field.

 A geometric corollary: the map $f$ is $\cR$-determined by its ``2-jet" on the critical locus,  $\cR f\supseteq\{f\}+\ca_\cR^2\cdot \RpX$.

\item 
 is the filtration criterion for $\cR$-orbits via tangent spaces: ``$\cR^{(j)}f\!\supseteq\! \{f\}\!+\!I^d\!\cdot\! \RpX$ \ iff  \  $T_{\cR^{(j)}}f\!\supseteq\! I^d\!\cdot\! \RpX$".

 (In positive characteristic the integers $j,d$ are restricted by the condition $ {d-ord(f)}<j\cdot char(\k)$, i.e. the
  characteristic must be large enough for a given $j,d$.)


A statement of this type has first appeared in \cite{Bruce.du-Plessis.Wall} for $\k=\R,\C$, and today it is  known in characteristic zero.
  In positive characteristic it strengthens the known criteria, \cite{Greuel-Pham.2017}, \cite{BGK.20},
 all being of type
 ``If $T_{\cR^{(j)}}f\supseteq I^d\cdot \RpX$ then
  $\cR^{(j+d)}f\supseteq \{f\}+I^{2d-ord(f)+1}\cdot \RpX$".

  To emphasize, the new statement is of ``iff" type.

  The proof uses the truncated exp/ln operators, $jet_N(Exp),$ $jet_N(Ln).$

\eee

\item   is the second warmup, addressing the $\cK$-equivalence.
\bee[\!\S5.1\!]
\item
 contains
 the ``{\em$\cK$-implicit function theorem}". It puts a large part of the tangent space
   into the orbit (roughly): $\cK f\supseteq\{f\}+(\ca^2 +\cm(f))\cdot\RpX$. As in the $\cR$-case, this is characteristic-free and
 much stronger than the known bounds (we give several examples).

 A geometric corollary: the map $f$ is $\cK$-determined by its ``2-jet" on the singular locus,  $\cK f\supseteq\{f\}+ \ca^2_\cK\cdot \RpX$.
\item
  is an immediate application:   the $\cK$-orbit of
  a reduced complete intersection curve germ $(C,o)\!\sset\! (\k^n,o)$, in terms of the semi-group of values.

\item
 is the filtration criterion for $\cK$-orbits: ``$\cK^{(j)}f\supseteq \{f\}+I^d\cdot \RpX$ iff $T_{\cK^{(j)}}f\supseteq I^d\cdot \RpX$".
  As in the $\cR$-case, for $char(\k)\!>\!0$  the characteristic is restricted by the condition $ d\!-\!ord(f)\! <\!j\!\cdot\! char(\k)$.

\quad  As in the $\cR$-case, this statement is known in characteristic zero, but  in positive (high enough) characteristic it strengthens the known criteria ``by a factor of 2".

   The proof uses (besides the truncated  maps $jet_N(Exp)$, $jet(Ln)$)  the  full
    Baker-Campbell-Hausdorff expansion, \S\ref{Sec.Preliminaries.Basic.Tools}.iv.
   (We bound the prime numbers appearing in the denominators.)
\eee

\rule{16cm}{0.5pt}

The $\cA$-equivalence is much more complicated. (Recall that the $\cK$-equivalence was initially introduced as an auxiliary,
 simpler version of $\cA$-equivalence.) A difficulty of $\cA$ (besides those mentioned in \S1.1) is that the classical
  arguments were heavily based on Thom-Levine's lemma, which holds only in zero or high characteristic, and on Mather-Malgrange preparation theorem, which assumes that $T_\cG f$ is an $R_X$-module.
   Without these results one needs special tools.

 The next two sections are the core of the paper.

\item is ``The basic theory of $T_\cA f$".
\bee[\!\S 6.1\!]
\item The annihilator ideal  $\ca_\cA\!=\!Ann[T^1_\cA f]\sset R_X$ (for $T^1_\cA f:=\quots{\RpX}{T_\cA f}$) is much more delicate than the corresponding   annihilators $\ca_\cR,\!\ca_\cK.$ We  establish     basic  properties  of  $\ca_\cA.$
\item The (extended) tangent space $T_\cA f$ is not an $R_X$-module. Considered as an $R_Y$-module it is often not finitely-generated (e.g. when $dim(X)>p$).
 We prove that $T_\cA f$ is a module over the extended ring $f^\#(R_Y)+\ca_\cA$. Moreover, $T_\cA f$ is finitely-generated over  $f^\#(R_Y)+\ca_\cA$
  in the most important case,
  ``$f$ is of finite singularity type" (i.e. $f$ is $\cK$-finite, i.e. the germ $V(f)\sset X$ has an isolated singularity).
\item As $T_\cA f$ is not an $R_X$-module, one does not have the usual Nakayama statement, ``$M\sseteq T_\cA f+\cm\cdot M$ implies  $M\sseteq T_\cA f$."
 We give weaker versions.

\item Similarly the assumption ``$T_{\cA^{(j)}}f\supseteq I^{d}\cdot \RpX$" does not imply $``T_{\cA^{(j+1)}}f\supseteq I^{d+1}\cdot \RpX".$
  Yet, we establish weaker statements of Artin-Rees type.
\item Unlike the tangent spaces $T_\cR f,T_\cK f$, the space $T_\cA f$ is not ``obviously filtration-closed".
  One cannot deduce $\overline{T_\cA f}=T_\cA f$ directly, just from the faithfulness/exactness of the completion functor $\liml_\leftarrow \quots{R_X}{I^\bullet}\otimes$.
 And one cannot apply the Artin approximation, as this is an ``inverse Artin problem".
  Yet, one has $\ca_\cA\cdot T_\cA f=\ca_\cA\cdot \overline{T_\cA f}.$ For maps of finite singularity type we prove:   $T_\cA f=\overline{T_\cA f}.$
\eee

\item is about large modules inside $\cA$-orbits.
\bee[\!\S 7.1\!]
\item recalls the auxiliary result, the Kostant-Rosenlicht theorem (``Orbits of unipotent algebraic groups are closed") over algebraically-closed fields.
  I add also an extension for arbitrary fields of characteristic zero (by M. Borovoi) and an example (by Z.Rosegarten) showing the pathology in positive characteristic, even when the field $\k$  is perfect.
\item contains the ``{\em $\cA$-implicit function theorem}", $\cA f\supseteq\{f\}+\ca^2\cdot \RpX+T_{\cL^{(1)}}f$.
 This is the main (and hardest) result of the paper.
 The proof is essentially more complicated than those in the $\cR,\cK$-cases.
 \bee
 \item The condition $\cA f\ni f+g$ is not an implicit function equation (unlike the conditions $\cR f\ni f+g$, $\cK f\ni f+g$).
   Yet, we convert it to an implicit function equation, modulo terms of high enough order.
    In this process we pass from  the source  $X$ to  the target  $(\k^p,o),$
  using the Weierstra\ss\ division.
 \item We obtain an arc solution of a system, i.e. a solution over the ring $R_Y[[t]]$.
  Then one passes to the finite jets, and finds  the ordinary solution for infinite number of $t$-values.
\item Use Zariski-closedness of the orbits of unipotent algebraic groups (by Kostant-Rosenlicht).
 If the field $\quots{\k}{\cm_\k}$ is not algebraically-closed, then one invokes an   argument of Galois-cohomology.
\eee
\item The first immediate corollary is for the $\cA$-orbits of stable maps: $\cA f\supseteq \{f\}+\ca^2_{\cK^{(o)}}\cdot\RpX$. This bound is completely new even in the
 classical case $R_X=\C\{x\},\R\{x\}$, and is much stronger than the known bounds.
 In fact this bound holds for a broader class of maps, a sufficient assumption is: $T_{\cA^{(o)}}f=T_{\cK^{(o)}}f,$ i.e. $\ca_{\cA^{(o)}}=\ca_{\cK^{(o)}}.$

 Then come other examples showing how to apply this $\cA$-$IFT$ theorem and how to extend numerous classical results to arbitrary characteristic.
  E.g. any $\cA$-finite map is $\cA$-finitely determined, various explicit determinacy bounds, and so on.

\item contains the filtration criterion,
``$\cA^{(j)}f\!\supseteq\! \{f\}\!+\!I^d \cdot  \RpX$ iff $T_{\cA^{(j)}}f\!\supseteq\! I^d\cdot \RpX$".
  As in the $\cR, \cK$-cases, for $char(\k)\!>\!0$ the integers $j,d$ are restricted by the characteristic,
   $ 2d\!-\!1\!-\!ord(f) \!<\!j\! \cdot\! char(\k).$

 This statement is known for $R_X=\k\{x\}$, $\k\in \R,\C$,  but is completely new in zero/positive characteristic.
As in the $\cR,\cK$-case, the proof uses the truncated maps $jet_N(Exp),$ $jet_N(Ln).$ But the additional (essential) ingredient is the algebraic Thom-Levine lemma of \S\ref{Sec.3.Related.Derivations.Autos}.

\item gives the geometric criterion of determinacy: $f$ is determined by its finite jet along the instability locus in the target, 
  $f(V(\ca_\cA))$.  (Or, along the locus $f^{-1}(f(V(\ca_\cA)))\cap Crit(f)$.)
 For $\cA$-finite maps (i.e.  $V(\ca_\cA)$ is a point) this extends the classical Mather-Gaffney criterion.
  Otherwise this seems to be new even in the $\R,\C$-cases.
\eee
\eee

\section{Notations, conventions, preliminaries}\label{Sec.Preliminaries}
\subsection{Local rings/germs of schemes}\label{Sec.Preliminaries.Rings.Germs}
 We study Maps($X,(\k^p,o)$), where $X=Spec(R_X),$ $(\k^p,o)=Spec(R_Y)$ are germs of schemes over the base ring $\k.$
 In the simplest case $\k$ is a field, classically $\k=\R,\C$. The case when $\k$ is not a field is important for deformations/unfoldings, e.g.
  $\k=\tilde\k[[t]]$, resp. $\tilde\k\{t\},\tilde\k\bl t\br,$
  see \cite{Kerner.Unfoldings}.
   Through the paper we use the multi-variables,
$x=(x_1,\dots,x_n),$ $y=(y_1,\dots,y_p),$ $t=(t_1,\dots,t_r).$

\bee[\!\bf i.\!]

\item   Denote by $\k^\times$ the group of invertible elements of $\k$. Sometimes we impose the condition ``$2,\dots,d\in \k^\times$".
 E.g. this holds if  $\k\supseteq\Q$ or $\k$ contains a field of characteristic larger than $d$.

\item   $(R_X,\cm),$    or just $(R,\cm),$ is one of the following local $\k$-algebras.
\bei
\item The formal power series $R=\k[[x]],$ or the quotient ring $R=\quots{\k[[x]]}{J}.$ Here  $\k$ is either a(ny) field or
    a complete local (Noetherian) ring, e.g.  $\k=\quots{\tilde\k[[t]]}{I}.$
 Then $R$ is complete \wrt the filtration $\cm^\bullet$.
\item The $\k$-analytic (locally convergent) power series $R= \k\{x\},$ or the quotient ring $R=\quots{\k\{x\}}{J}$.
   In this case  $\k$ is either a normed ring (e.g a normed field), complete \wrt its (non-discrete) norm, or  $\k=\quots{\tilde\k\{t\}}{I},$ where $\tilde\k$ is a normed ring (as before).
 The relevant cases of complete normed fields are $\R,\C,\Q_p$.

 \item The  algebraic power series $R= \k\bl x\br $  (i.e. power series in $x$ that satisfy monic polynomial equations over $\k[x]$)
 or the quotient ring $R=\quots{\k\bl x\br }{J}.$
   Here  $\k$ is either a(ny) field or a local excellent Henselian ring.
    Then $ \k\bl x\br $  is the Henselization of   the ring $\k[x]$ at the ideal $(x)$, constructed as the filtered colimit of \'etale $\k[x]$-algebras.
  See \cite{Moret-Bailly} for more detail.

The elements of $\R\bl x\br$ are exactly the Nash germs, see e.g. \cite{Ruiz} for the basic introduction.
 Occasionally we call the elements of $\quots{\k\bl x\br}{J}$  ``$\k$-Nash germs" rather than ``algebraic germ", to avoid any confusion with rings like $\k[x],$ $\k[x]_\cm.$
\eei

\

The assumptions on $\k$ (local complete,  resp. normed and complete for its norm, resp. local excellent Henselian) are needed to ensure the
 implicit function theorem,   Weierstra\ss\ finiteness, and Artin approximation, see \S\ref{Sec.Preliminaries.Basic.Tools}.

 Denote the maximal ideal of $\k$ by $\cm_\k.$
Occasionally we denote the image of $x$ in $R$ by the same letter. Thus the maximal ideal is $\cm=\cm_\k+(x)\sset R.$
  \quad  If $\k$ is a field then $\cm=(x).$

  The elements of $R$ are represented by  power series. Abusing notations we write $f(x),$ to emphasize the dependence on $x.$

 \item
 Geometrically the source of the map $X=Spec(R)\sseteq(\k^n,o)$ (or $X=Spec(R_X)$) is a formal/$\k$-analytic/$\k$-Nash germ, and $x$ denotes the local coordinates on the ambient space $(\k^n,o).$
\\
 The target of the map  is $Y\!:=\!(\k^p,o)\!:=\!Spec(R_Y),$ the formal/$\k$-analytic/$\k$-Nash germ (corresponding to the type of $X$).
  Namely, $R_Y$ is one of the rings
  $ {\k[[y]]} ,$ $ {\k\{y\}} ,$ $ {\k\bl y\br } ,$ with $\k$  a field or a ring, as before.

Below we use $Y$ and $(\k^p,o)$ interchangeably.
\item

Through the paper we   have to compose elements $f(z)\in R[[z]]$ with $g\in R.$ For $g\in \cm\sset R$ the element $f\circ g\in \hR$ is well defined.  For $g\!\in\! R\!\smin\! \cm$ the element $f\!\circ \! g $ can be not well defined. But in the proof of theorems \ref{Thm.IFT.R.case}, \ref{Thm.IFT.K.case}, \ref{Thm.IFT.A.case}
  the series $f$ are of special form, $f(z)\!\in\! R[[\cm\cdot z]][z]$, i.e. $f$ is a polynomial in $z$, whose coefficients are power series in $\cm\!\cdot\! z$.
  In this case the composition $f\!\circ\! g$ is well defined for any $g\!\in\! R. $

\item
Take a tuple of elements $f=(f_1,\dots, f_p)\in R^{\oplus p}$ and an ideal $I\sset R.$ Denote the $I$-order of the ideal $(f)\sset R$  by $ord(f)$. Namely, $(f)\sseteq I^{ord(f)}$ and $(f)\not\sseteq I^{ord(f)+1} $.
 Sometimes we use the ideal $I^{ord(f)-2}$. If $ord(f)\le 2$ or $I=R_X$ we put $I^{ord(f)-2}:=R_X$.

\item Recall the general fact: any $R_X$-submodule $M\sseteq \RpX$ is finitely generated. Proof: $\RpX$ itself is a f.g. module over the Noetherian ring $R_X$,
 therefore $\RpX$ is a Noetherian module. Hence $M$ is f.g.

\item We consider only local homomorphisms of local $\k$-algebras, i.e. $\phi: R_Y\to R_X$ with $\phi(\cm_Y)\sseteq \cm_X.$ Such homomorphisms are necessarily Krull-continuous, see e.g. Lemma 15.36.2 of \cite{Stacks}. Therefore they act by substitution, $\phi(f(y))=f(\phi(y)).$
     In particular, the homeomorphism $\phi$ is determined by its action on the generators $y=(y_1,\dots,y_p).$
\eee

\subsection{The identification  $\mathbf{\Maps=\cm \cdot \RpX}$ for $Y=(\k^p,o)$}\label{Sec.3.Maps.of.Germs.Identification}
 A map of germs $f:X\to Y$ is defined algebraically by the (local) homomorphism of $\k$-algebras, $f^\sharp\in Hom(R_Y,R_X).$
 Fix generators $y=(y_1,\dots,y_p)$ in $R_Y.$
 Any   homomorphism $R_Y\to R_X$ is determined by the image of $y $ in $\cm\sset R_X,$ see \S\ref{Sec.Preliminaries.Rings.Germs}.vii.
  Vice-versa, any such image extends to a homomorphism.
  Accordingly  we identify $\Maps\!=\! Hom_\k(R_Y,R_X)\!=\! \cm\! \cdot\! \RpX.$
   This identification equips the set $\Maps$ with the $R_X$-module structure.

If $\k$ is a field then $\cm=(x)\sset R_X$ and the maps fix the origin, $X\ni o\to o\in (\k^p,o)$.
 If $\k$ is a local ring then $\cm= \cm_\k+ (x)$, and $\Maps$ can be interpreted as families of maps, possibly displacing the origin.

\

Given a map $f\in \Maps=\cm \cdot \RpX\sset\RpX,$ the space $\RpX$ becomes a (not finitely-generated) $R_Y$-module via the composition,
\beq\label{Eq.3.R_Y.action.on.R_X}
q(y)\cdot \RpX:=f^\#(q(y))\cdot \RpX:=q(f_1,\dots,f_p)\cdot\RpX.
\eeq
We have the $R_Y$-submodule $f^\#(R^{\oplus p}_Y)\sseteq f^*(R^{\oplus p}_Y)= \RpX.$
Here $f^\#(y)\sseteq f^*(y)=(f)\sset R_X$, and the first inclusion is usually proper.

Recall the classical notations, $\omega f(\cm_p)$ for $f^\#(y)$, and $f^*(\cm_p)$ for $f^*(y)$, \cite{Mond-Nuno}.

\subsection{The basic tools}\label{Sec.Preliminaries.Basic.Tools} (They are used repeatedly in \S\ref{Sec.R.Orbits}-\ref{Sec.A.Orbits}.)
\bee[\!\bf i.\!]
 \item{\bf Implicit function theorem with unit linear part, $\mathbf{IFT_\one}$.}
We often have to resolve implicit function equations of type $z=f(z)$ for a vector of (formal/analytic/algebraic) power series:
\beq
f(z)\in  (\cm+(z)^2)\cdot R[[z]]^{\oplus p}, \quad
 \text{ resp. } f(z)\in(\cm+(z)^2)\cdot R\{z\}^{\oplus p}, \quad  \text{ resp. }  f(z)\in(\cm+(z)^2)\cdot R\bl z\br^{\oplus p}.
 \eeq
 Here $z=(z_1,\dots,z_p)$ are the unknowns and the  equations are non-polynomial in $z.$   We look for the solution $z\in \cm\sset \RpX.$
 There exists the unique formal solution $\hz(x)\in \cm\cdot \hat{R}^{\oplus p}.$
  It is    obtained, e.g. by the order-by-order procedure.

The implicit function theorem   holds for our rings and ensures: this formal solution  belongs to $R.$

 \noindent\bei
 \item For the ring  $\quots{\k[[x]]}{J}$, with $\k$ any complete local ring,  this statement is  trivial.
\item   $IFT_\one$ holds in the rings $\quots{\k[[x]]}{J}$, $\quots{\k\{x\}}{J}$, $\quots{\k\bl x\br}{J}$, for any (normed) field $\k$,   see example 2.2 of
  \cite{BGK.20}.
\item  $IFT_\one$ holds in the ring  \!$ \quots{\k\{x\}}{J}$\! for \!$\k$\! a normed ring,  complete \wrt its norm,\cite[pg.84]{Abhyankar}.

\item $IFT_\one$ holds in the ring  $\quots{\k\bl x\br}{J}$, for $\k$ a local, excellent, Henselian ring, \cite{Lafon}, see also page 4 of\footnote{\cite{Denef-Lipshitz} assume $\k$ a field or a DVR (and denote $\k$ by $R$). But their proof uses only the Weierstra\ss\ preparation, which holds for $\k$-local, excellent, Henselian.} \cite{Denef-Lipshitz}.

\eei


\noindent Sometimes the system to resolve  has the form $z=h(z)+c$, for a vector $c\in R^{\oplus p}$, where the entries of the vector $h(z)$
 belong to the subset $\cm\cdot R[[\cm\cdot z]][z]\sset R[[z]]$.
 This case is reduced to that of $z=f(z)$ by the substitution $\tz=z-c$.

\item{\bf Weierstra\ss\ division with remainder.}  Take a local subring $R\sseteq \k[[x]]$, here $\k$ is a local ring (e.g. a field).
 The Weierstra\ss\ division holds in $R $ if
 for any $x_n$-regular element  $f \in R $   of $x_n$-order $d_n$, and any $g\in R $,
 one can present $g=q\cdot f+r$, where $q\in R $ and $r\in  R \cap \k[[x_1,\dots,x_{n-1}]] [x_n]$, with $deg_{x_n} r<d_n$.
 This division with remainder holds in our rings of \S\ref{Sec.Preliminaries.Rings.Germs}.ii.
   For   $\k\{x\}$ see, e.g. \cite[pg.75]{Abhyankar}, for   $\k \bl x\br$   see \cite{Lafon.65}. For more general results and further references see \cite{Moret-Bailly}.

\item{\bf Weierstra\ss\ finiteness condition} holds for a pair of local $\k$-algebras $R_X,R_Y$ if for any homomorphism $\phi: R_Y {\to}R_X$
 and  any finitely generated $R_X$ module $M$, with the quotient $\quots{M}{\phi^*(\cm_Y)\cdot M}$ f.g. over $\k$, the module $M$ is f.g. over $R_Y$.

This condition is implied by the Weierstra\ss\ division condition. In particular it holds for the pairs
 $(\quots{\k[[x]]}{J},\k[[y]])$, $(\quots{\k\{x\}}{J},\k\{y\})$, $(\quots{\k\bl x\br}{J},\k\bl y\br)$.

A remark: the $C^\infty$-version of Weierstra\ss\ finiteness is a deep preparation theorem of Malgrange.

\item{\bf Baker-Campbell-Hausdorff formula.} Suppose $ \k\supseteq\Q$. Let $M$ be a $\k$-module (possibly not finitely generated), with a filtration $M_\bullet$. Suppose $M$ is $M_\bullet$-complete.
 Suppose some operators (not necessarily $\k$-linear) $\xi,\eta\circlearrowright M$  are filtration-nilpotent, i.e. $\xi,\eta(M_\bullet)\sseteq M_{\bullet+1}$.
  Then the operators $exp(\xi)$, $\exp(\eta)$, $\exp(\xi+\eta)$ are defined on   $M$. One has the classical relation
  $exp(\xi)exp(\eta)=exp(\sum^\infty_{l=1}p_l(\xi,\eta))$. Here $p_l(\xi,\eta)$ is a homogeneous polynomial of degree $l$ in
   non-commuting variables $\xi,\eta$.
The expansion begins as follows:
$exp(\xi)\cdot exp(\eta)=exp\big(\xi+\eta+\frac{[\xi,\eta]}{2}+\frac{[\xi,[\xi,\eta]]-[\eta,[\xi,\eta]]}{12} +\cdots\big)$.
 Below we use the full expansion, see, e.g. \cite[\S5 of chapter 5]{Jacobson} or \cite[pg.29]{Serre.Lie} or \cite[pg.160]{Bourbaki.Lie}:
   \beq\label{Eq.Baker.Campbell.Hausdorff}
\sum^\infty_{l=1}p_l(\xi,\eta)=\sum^\infty_{i=1}\frac{(-1)^i}{i}\sum_{\substack{r_1+s_1>0\\\dots\\r_i+s_i>0}}
 \frac{[\xi^{r_1}\eta^{s_1}\cdots \xi^{r_i}\eta^{s_i}]}{\Big(\suml^i_{\bullet=1}(r_\bullet+s_\bullet)\Big)\Big(\prodl^i_{\bullet=1} r_\bullet!\Big)
 \Big(\prodl^i_{\bullet=1} s_\bullet!\Big)}.
   \eeq
Here $[\xi^{r_1}\eta^{s_1}\cdots \xi^{r_i}\eta^{s_i}]$ is a certain combination of commutators with $\Z$ coefficients.
 Comparing the two sides we get $\sum^i_{\bullet=1}(r_\bullet+s_\bullet)=l$ and therefore:
 \beq
 p_l(\xi,\eta)=\frac{1}{l}\sum^l_{i=1}\frac{(-1)^i}{i}\sum_{\substack{\{r_\bullet+s_\bullet>0\}_\bullet\\\sum^i_{\bullet=1}(r_\bullet+s_\bullet)=l}}
 \frac{[\xi^{r_1}\eta^{s_1}\cdots \xi^{r_i}\eta^{s_i}]}{\Big(\prod^i_{\bullet=1} r_\bullet!\Big)
 \Big(\prod^i_{\bullet=1} s_\bullet!\Big)}.
 \eeq
 In particular, the rescaled polynomial $(l!)^4\cdot p_l(\xi,\eta)$ has only integer coefficients.

 We emphasize that the relation $exp(\xi)\cdot exp(\eta)=\dots$ is purely formal,
  and no linearity/algebraicity is assumed on the action $\xi,\eta\circlearrowright M$.

Recall that in the analytic case, i.e. $R=\quots{\k\{x\}}{J}$, $\k\supseteq \Q$, the series $\sum^\infty_{l=1}p_l(\xi,\eta)$
 can be nowhere convergent for some derivations $\xi,\eta\in Der_\k(R)$. However this series is locally convergent, i.e. the resulting
 derivation is analytic, $\sum^\infty_{l=1}p_l(\xi,\eta)\in Der_\k(R)$, if one assumes the filtration-nilpotence,
  $\xi,\eta(M_\bullet)\sseteq M_{\bullet+1}$.
   Indeed, in this case we get the (filtration-unipotent) analytic automorphisms, $exp(\xi),exp(\eta)\in Aut_\k(R)$.
    (See  e.g. the proof of part 2 of Lemma 3.17 in \cite{B.K.motor}.) Their product is a filtration-unipotent automorphism as well,
     $exp(\sum^\infty_{l=1}p_l(\xi,\eta))\in Aut_\k(R)$. And therefore $\sum^\infty_{l=1}p_l(\xi,\eta)\in Der_\k(R)$.

\item{\bf Artin approximation.} Let $R_X$ be one of $\quots{\k[[x]]}{J}$, $\quots{\k\{x\}}{J}$, $\quots{\k\bl x\br}{J}$, with $\k$
 a local ring,
 see \S\ref{Sec.Preliminaries.Rings.Germs}.ii. Take a (finite) system
 of implicit function equations, $F(y)=0$, here $F(y)\in R_X[[y]]^{\oplus N}$, resp. $R_X\{y\}^{\oplus N}$, resp. $R_X\bl y\br^{\oplus N}$.
 Suppose we have an order-by-order solution, i.e. a sequence $\ty_\bullet(x)\in R^{\oplus q}_X$ satisfying:
  $F(\ty_\bullet(x))\in \cm^\bullet\cdot R^{\oplus N}_X$.
  The Pfister-Popescu theorem, \cite{Pfister-Popescu} (see also \S2.4.1 of \cite{BGK.20}),  ensures a formal solution, i.e. $\hy(x)\in \hR_X^{\oplus q}$
   satisfying $F(\hy(x))=0$.

Suppose $R_X$ is one of   $\quots{\k\{x\}}{J}$, $\quots{\k\bl x\br}{J}$, see \S\ref{Sec.Preliminaries.Rings.Germs}.ii.
By the Artin approximation theorem, \cite{Artin.68}, \cite{Artin}, this formal solution is approximated by $R_X$-solutions. Namely, for each $d\ge 1$
 there exists $y_d(x)\in R_X^{\oplus q}$ satisfying: $F(y_d(x))=0$ and $\hy(x)-y_d(x)\in (x)^d\cdot \hR_X^{\oplus q}$.

A remark: the general Artin approximation holds for {\em polynomial} equations over any Henselian local ring. In our case the equations $F(y)=0$
 are non-polynomial, they are (algebraic/analytic) power series.
   Therefore we restrict to the particular rings  $\quots{\k\{x\}}{J}$, $\quots{\k\bl x\br}{J}$.
  See \cite[\S2.4.2]{BGK.20}, \cite{Rond}, \cite{Popescu} for more detail.

 The equation  of $\cA$-equivalence is $\Phi_Y\circ f\circ\Phi_X^{-1}=\tf$. This is not an implicit function equation.
  And the Artin approximation does not work in the analytic case, $\k\{x\}$,
 due to the example of \cite{Gabrielov}, see  \cite[Fact.1.4]{Shiota.1998}.
 The left-right version of the Artin approximation has been established:
 \bei
 \item     $R_X=\R\bl x\br,$ for any $f$, \cite[Fact.1.3]{Shiota.1998}, \cite[Theorem.4]{Shiota.2010}.
  \item $R_X=\R\{x\},$  for  $f$ of finite singularity type, \cite[Fact.1.7]{Shiota.1998}.
   \eei
   The proof for $\R\bl x\br$ is characteristic-free,  and   is valid over any field. \!But the proof for $\R\{x\}$ is  heavily based on the $C^\infty$-topology.
   \!These statements are extended in \!\cite{Kerner.LRAP}\! to the rings \!$\quots{\k\bl x\br}{J},\! \quots{\k\{ x\}}{J}.$
\eee

\subsection{\!\!Coordinate changes (ring automorphisms, $Aut_X$) and vector fields (ring derivations, $Der_X$)}\label{Sec.Preliminaries.Aut.Der}
\bee[\!\!\bf i.\!]
\item  Denote by \!$Aut_\k(R)$\! the group of (local) $\k$-linear automorphisms of $R$. The elements of $R$ are presentable by power series,
     the automorphisms act
 by substitution, \!$\phi(f(x))\!=\!f(\phi(x))$ for $\phi(x)\!\in\! \cm.$
   See \S\ref{Sec.Preliminaries.Rings.Germs}.

 Geometrically (for \!$J\!=\!0$) these are ``the   coordinate changes" of
 the germ $X\!=\!Spec(R).$\!   See \cite[\S2.3.2]{BGK.20} for more detail.
 Because of this identification we abbreviate $Aut_\k(R)$ to $Aut_X$.


\item  Denote by $Der_\k(R)\sset End_\k(R)$ the set of $\k$-linear derivations. This is an $R$-module.
 Geometrically this is the module of germs of vector fields on $X=Spec(R)$, classically
    denoted by $\Theta_X$.

 For  $R=\k[[x]]$, $\k\{x\}$, $\k\bl x\br$ this module is   generated by the partials
  $Der_\k(R)=R\bl\frac{\di}{\di x_1},\dots,\frac{\di}{\di x_n}\br,$ see e.g. \cite{Matsumura} Theorem 30.6
  More generally, for $R=\quots{S}{J}$ with $S=\k[[x]], \k\{x\},\k\bl x\br$, we have the presentation
\beq
J\cdot Der_\k(S)\to DerLog(J)\to Der_\k(R)\to0.
\eeq
Here $DerLog(J)=\{\sum c_j \frac{\di}{\di x_j}|\ \sum c_j \frac{\di}{\di x_j}(J)\sseteq J\}\sset Der_\k(S)$
 is the module of log-derivations.
 Its elements  correspond to
  the vector fields tangent to the subgerm $V(J)\sset (\k^n,o)$.

Denoting $X=Spec(R)$  we   abbreviate $Der_\k(R)$ to $Der_X$ and $End_\k(R)$ to $End_X$.

\noindent Applying  derivations to an element \!$f\!\in  \! R^{\oplus p}$\! we get the \!(finitely generated)\! submodule \!$Der_X(f)\!\sseteq\! R^{\oplus p}.$

\item
In $\R/\C$-analytic (or $C^r$-) geometry  one integrates vector fields to flows. Taking a flow  at time $t=1$ one gets the map $\Theta_X\to Aut_X$.
 Over an arbitrary ring $\k$ (without any topology/integration) the maps $Der_X\rightleftarrows Aut_X$ are constructed in \S\ref{Sec.Preliminaries.Exp.Ln},
  using the following filtration-topology.
 (For more detail  see   \cite[\S2.2.3]{BGK.20}.)

Fix a descending filtration by ideals, $R=I_0\supset I_1\supset\cdots$,
 with $\cap I_j=(0)$. This defines the filtration topology on $R$.
  Accordingly, for each $j\ge0$, one considers the coordinate changes that are identity modulo higher order terms:
\beq
Aut^{(j)}_X:=\{\Phi\in Aut_X| \quad \Phi(I_d)=I_d \quad\quad \text{ and }\quad\quad \Phi|_{\quots{I_d}{I_{d+j}}}=Id|_{\quots{I_d}{I_{d+j}}},\ \forall\ d\ge1\}.
\eeq
 Thus $Aut^{(0)}_X$ consists of automorphisms that preserve the filtration, while    $Aut^{(1)}_X$ consists of
 ``topologically unipotent" automorphism.
 For the particular filtration $\cm^\bullet\sset R$
  one has  $Aut_X=Aut^{(0)}_X$.
   The sets $Aut^{(j)}_X$ are normal subgroups, and we get the filtration
    $ Aut^{(0)}_X\trianglerighteq Aut^{(1)}_X\triangleright\cdots$.
\\Similarly  we consider (for every $j\!\in \!\Z$) ``topologically nilpotent" $\k$-linear endomorphisms and derivations.
\bei
\item
 $End^{(j)}_\k(R)\!:=\!\{\phi\!\in \!End_\k(R)|\ \phi(I_d)\!\sseteq\! I_{d+j},\ \forall\ d\ge1\}$
\\Below we often abbreviate $End_X\!:=\!End_\k(R)$ and $End^{(j)}_X\!:=\!End^{(j)}_\k(R).$
\item
$Der^{(j)}_X:=\{\xi\in Der_X|\ \xi(I_d)\sseteq I_{j+d},\ \forall\ d\ge1\}=Der_X\cap End^{(j)}_\k(R).$


We get the filtration by $R$-submodules, $Der_X\supseteq\cdots\supseteq  Der^{(-1)}_X\supseteq  Der^{(0)}_X\supseteq Der^{(1)}_X\supset  \cdots$.

E.g. for the filtration by powers of ideals, $I^\bullet\sset R$, we have  $Der_X=Der^{(-1)}_X$ and
   $Der^{(j)}_X\supseteq I^{j+1}\cdot Der_X$.

\eei

\eee

\subsection{The maps $exp,ln: Der_X\leftrightarrows Aut_X$ and their approximations.  Vector field integration}\label{Sec.Preliminaries.Exp.Ln}

Let $(R ,\cm , I^\bullet)$ be a  local filtered $\k$-algebra of \S\ref{Sec.Preliminaries.Rings.Germs}.ii. Assume   $\k\supseteq\Q$ and $R$ is complete \wrt the filtration.
 Then the derivations and automorphisms of $R$ admit the exponential and logarithmic maps, see Theorem 3.6 of \cite{BGK.20}:
 \vspace{-0.4cm}
\beq\label{Eq.Exp.Ln.full.series}
\forall\ j\ge1: \quad\quad  Der^{(j)}_X\stackrel{exp}{\to} Aut^{(j)}_X,\quad   exp(\xi ):=e^{\xi }=\sum\nolimits^\infty_{i=0}\frac{\xi^i}{i!},
\eeq
 \vspace{-0.4cm}\[
\forall\ j\ge1:\quad\quad Aut^{(j)}_X\stackrel{ln}{\to}Der^{(j)}_X ,\quad   ln(\Phi ):=-\sum\nolimits^\infty_{i=1}\frac{(\one-\Phi )^i}{i}.
\]

When $\k\not\supseteq\Q$ or the ring $R$ is not $I^\bullet$-complete, one cannot use the classical exponential/logarithmic maps. But in many cases
  one can use the following ``$N$'th jet
 approximations" to these maps.  Fix $1\le N\le \infty.$
  For our filtration criteria, \S\ref{Sec.R.Filtration}, \S\ref{Sec.K.Filtration}, \S\ref{Sec.A.Filtration},
   we impose the conditions (for all $j\ge1$):
\\$\bullet$ $\boldsymbol{jet_N(Exp)\!:}$ every derivation $\xi\!\in\! Der^{(j)}_X$  defines  an
automorphism  $\!jet_N(e^\xi)\!:=\!\sum^N_{i=0}\!\!\frac{\xi^i}{i!}\!+\!\phi_\xi\!\!\in\! Aut^{(j)}_X\!.$

 (Here $\phi_\xi\in End^{(jN+1)}_X$ is an  auxiliary $\k$-linear endomorphism, see \S\ref{Sec.Preliminaries.Aut.Der}.iii.)
\\$\bullet$ $\!\!\boldsymbol{jet_N(Ln)\!:}$ every automorphism $\!\phi\!\in \! Aut^{(j)}_X$ defines a
 derivation  $\!jet_N(ln(\phi))\!:=\!-\!\sum^N_{i=1}\! \frac{(\one-\phi )^i}{i}\!+\!\xi_\phi\!\!\in\!\! Der^{(j)}_X\!.$

(Here $\xi_\phi\in End^{(jN+1)}_X$ is  an  auxiliary $\k$-linear endomorphism.)
\\These conditions imply, in particular, $2,\dots,N\in \k^\times$, see \S\ref{Sec.Preliminaries.Rings.Germs}.i. Thus $char(\k)>N$ or $char(\k)=0.$

\

For our ``$\cG$-implicit function theorems", \S\ref{Sec.R.IFT},\S\ref{Sec.K.IFT},\S\ref{Sec.A.IFT}, we impose a  variation of  condition $jet_1$.
 Take a derivation $\xi\!\in\! Der_X$ with $\xi(\cm)\sseteq\cm^2.$
  If \!$J\!=\!0$\! then  the map $f(x)\!\to\! f(x\!+  \xi(x)),$\!  defines an automorphism  $\Phi\!\in\! Aut_X.$\!
 If  $J\!\neq\!0$  then $f(x)\!\to\! f(x\!+  \xi(x))$ (defined on the representatives)
   is not a well-defined self-map of $R ,$ as it does not preserve \!$J.$  Therefore we impose the condition:

\bei
\item $\boldsymbol{jet_0:}$ Any map $x\to x+\xi(x)$,  with $\xi\in Der_X$ and $\xi(\cm)\sseteq \ca^2\cdot \cm^j$ (for some $\ca\sseteq\cm$ and  $j\ge0$),
    extends to an automorphism $\Phi\in Aut_X$ that satisfies: $\Phi(x)-x-\xi(x)\in \ca^3\cdot\cm^j$.
     Moreover, presenting $\Phi=Id+\xi+\varphi_\xi$, one has (for $c_1,c_2\in \k$): $\varphi_{c_1\xi_1+c_2\xi_2}$ is a power series in $c_1,c_2$.
\eei


\bex\label{Ex.jet_N.conditions.hold}\bee[\!\!\bf i.\!]
\item  The  conditions $jet_N(Exp),jet_N(Ln)$ hold, e.g. in the following cases.
\bei
\item $R=\quots{ \k[[x]]}{J} $, with $\k\supseteq \Q$, and any $0\le N\le\infty$. Here one takes just the full power series in \eqref{Eq.Exp.Ln.full.series}.
 \item $R=\quots{\k\{x\}}{J}$ or $R=\quots{\k\bl x\br}{J}$, with $\k\supseteq \Q$, and any $0
 \le N<\infty$. Here one shows that such $\phi_\xi$, $\xi_\phi$ do exists,
 under very weak assumptions on the filtration: $Der^{(1)}_X(\cm)\sseteq\cm^2  .$
 For more detail  see \cite[Theorem 3.6]{BGK.20}.
\item  $R= \k[[x]]$, $\k\{x\}$, $ \k\bl x\br $, (and any $\k$), with  $2,\dots,N\in \k^\times$ or $N=0$,
 see \S3.2.2 of \cite{BGK.20}.

\eei
 \item
 The condition   $  jet_N(Exp)$  (resp.  $  jet_N(Ln)$) obviously  implies   $  jet_l(Exp)$  (resp. $  jet_l(Ln)$)   for $l<N$.
\item
The condition $jet_\infty(Exp)$ means: every vector field $\xi\in T_X$ integrates to the flow $\Phi_t\circlearrowright X$.

\item   For  $char(\k)\!>0$  and  $J\!\neq\! 0$  already the condition $jet_1$  is  non-trivial, see \S3.2 of \cite{BGK.20}.
 For example, let $\k$ be a field of characteristic $p$ and take $f(x,y)=x^p+y^{3p}\in \k[[x,y]]$.
  Take the ring $R=\quots{\k[[x,y]]}{(f)}$ filtered by the ideals $(x,y)^\bullet$.
  Take  $\xi:=y^2\frac{\di}{\di x}\in Der^{(1)}_X$. The conditions $jet_0$ and $ jet_1(Exp)$   imply: there exist $\phi_x,\phi_y\in (x,y)^3$ such that
   the map $(x,y)\to (x+y^2+\phi_x,y+\phi_y)$ defines an automorphism of $R$. But then we must have
   \beq
f(x,y)\to f(x+y^2+\phi_x,y+\phi_y)=  x^p+y^{2p}+(\phi_x)^p+y^{3p}+(\phi_y)^{3p} \in (f(x,y))=(x^p+y^{3p}).
   \eeq
And this condition is clearly non-resolvable because of the term $y^{2p}$.
\eee
\eex



\section{Group  actions  $\cG\!\circlearrowright\!  \Maps$  and  the  corresponding  tangent  spaces  $T_\cG f$}\label{Sec.3.Maps.of.Germs.Group.Actions}
Classically one studies maps of germs $(\k^n,o) \to (\k^p,o)$, for $\k=\R,\C$, up to the right ($\cR$), left ($\cL$),
 contact ($\cK$) or left-right ($\cA$)
 equivalences.  We set up the corresponding general notions.

Let $X$ and \!$Y\!=\!(\k^p,o)$ be the formal/$\k$-analytic/$\k$-Nash germs, i.e. $(R_X,\!R_Y)$ is one of the   pairs,
  $(\quots{\k[[x]]}{J},\k[[y]]),(\quots{\k\{x\}}{J},\k\{y\}),(\quots{\k\bl x\br}{J},\k\bl y\br),$ see \S\ref{Sec.Preliminaries.Rings.Germs}.ii.

\subsection{The basic group-actions of Singularity Theory,  $\mathbf{\cR, \cL, \cC, \cK, \cA\circlearrowright \Maps}$}\label{Sec.3.Groups.R.K.A}

\bee[\bf\!i.\!]
\item The elements of $R_X$ are presentable by power series in $x.$
 Therefore the ``coordinate changes", $Aut_X\!\circlearrowright\! X,$ induce the  action $Aut_X\!\circlearrowright\! \Maps\!=\!\cm\!\cdot\! \RpX$ by composition on the right,
  \cite[\S3.1]{Mond-Nuno}

\beq
  \Phi_X(f)=\Phi(f_1,\dots,f_p)=f\circ\Phi^{-1}_X=(f_1\circ\Phi^{-1}_X,\dots,f_p\circ\Phi^{-1}_X),\hspace{2cm}
  f\in \cm\cdot \RpX.
  \eeq
 Traditionally one  denotes this group action by $\cR:=Aut_X.$
The action $\cR\circlearrowright \cm\cdot \RpX$ is $\k$-linear.

  \item
 Present the automorphisms $Aut_Y\circlearrowright R_Y$ explicitly, $(y_1,\dots,y_p)=y\to \Phi_Y(y)=\big(\Phi_1(y),\dots,\Phi_p(y)\big)$.
 Accordingly define the left action $Aut_Y\circlearrowright \Maps=\cm\cdot \RpX$,   by $f\to \Phi_Y(y)|_f:=(\Phi_1(f),\dots,\Phi_p(f)).$
  See  \cite[\S3.2]{Mond-Nuno}. Traditionally one denotes this group action by $\cL:=Aut_Y.$

Note that the group $\cL$ does not act on the whole module $\RpX$,
 as the composition $\Phi_Y(y)|_f$ is not defined for $f\not\in \cm\cdot \RpX$.
  (See \S\ref{Sec.Preliminaries.Rings.Germs}.iv.)

The action $Aut_Y\circlearrowright R_Y$ is $\k$-linear, but the action $\cL\circlearrowright \cm\cdot\RpX$ is neither additive nor multiplicative.

Occasionally we write also $Aut_X,Aut_Y\circlearrowright \Maps.$
  \item
 The contact transformations  $X\times Y\to X\times Y $ are defined by $(x,y)\to(\Phi(x),\Psi(x,y))$, see e.g. \cite[pg. 108]{Mond-Nuno}.
 Here
  $\Phi\in Aut_X$ and
 $\Psi\in Maps(X\times Y,Y)$ with $\Psi(o,-)\in Aut_Y$. Algebraically, in the formal case:
  $\Psi(\cm_X,\cm_Y)\in \cm_Y\cdot R_X[[y]]^{\oplus p}$, with $\Psi(o,\cm_Y)\in \k[[y]]^{\oplus p}.$ 
   (And similarly in the $\k$-analytic/$\k$-Nash cases.)
\\
Accordingly we  define the contact group action   $\cK\circlearrowright \Maps= \cm\cdot \RpX$   by $f(x)\to \Psi(x,f(\Phi^{-1}(x))$.

Taking $\Phi\!=\!Id_X\in Aut_X$ we get the subgroup $\cC\sset \cK$. It acts on the fibres of the projection $X\!\times\! Y\!\to\! X.$
 This presents the contact group as the semi-direct product  $\cK\!=\!\cC\!\rtimes\! \cR,$ see   pg. 156 of \cite{AGLV} or page 109 of \cite{Mond-Nuno}.

  The action $\cK\circlearrowright  \Maps=\cm\cdot \RpX $ is neither additive nor $\k$-multiplicative.
 But the $\cK$-orbits coincide with the orbits of the much smaller group $\cK^{lin}:=GL(p,R_X)\rtimes Aut_X\circlearrowright \cm\cdot \RpX .$
 For $R\!=\!\k\{x\}$, $\k\in \R,\C$ this is well known,
  see e.g. \cite[pg.110]{Mond-Nuno}. For  the  general  case  see  \cite[pg.123]{B.K.motor}.
 This action $\cK^{lin}\circlearrowright\Maps$ is  $\k$-linear.

Two maps $f,\tf\in \cm\cdot \RpX$ are $\cK^{lin}$-equivalent iff the corresponding local $\k$-algebras are isomorphic,
   $\quots{R_X}{(f)}\cong\quots{R_X}{(\tf)}.$ And the later is an isomorphism of scheme-germs, $Spec(\quots{R_X}{(f)})\cong  Spec(\quots{R_X}{(\tf)}).$
   Therefore the $\cK^{lin}$, $\cK$-equivalences coincide with the $V$-equivalence of \cite{Martinet.1976}.

\item
 The action of the left-right group,
    $\cA:=\cL\times\cR \circlearrowright \Maps=\cm\cdot \RpX$, is defined by
\beq\label{Eq.A.action}
(\Phi_Y,\Phi_X)(f):=\Phi_Y\circ f\circ \Phi_X^{-1}:=\Phi_Y(y)|_{f\circ \Phi_X^{-1}}.
\eeq
This action is compatible with the product structure:
  $(\Phi_Y,\Phi_X)\circ (\tilde\Phi_Y,\tilde\Phi_X)=(\Phi_Y\circ \tilde\Phi_Y, \Phi_X\circ \tilde\Phi_X)$. In particular,
   $(\Phi_Y,\Phi_X)= (\Phi_Y,Id_X)\circ (Id_Y,\Phi_X)= (Id_Y,\Phi_X) \circ (\Phi_Y,Id_X)$.
\eee

While the actions $\cR,\cK^{lin}\circlearrowright \cm\cdot \RpX $ are $\k$-linear,
 the actions  $\cL,\cA\circlearrowright \cm\cdot \RpX$  are neither additive nor $\k$-multiplicative.
   To my knowledge it is not known whether/how the $\cA$-action can be ``reduced"
  to a $\k$-linear one. See \cite{Damon}, \cite{Mond-Montaldi}, \cite{Houston.Wik.Atique} for the related works.

\subsection{The (extended) tangent spaces  $\mathbf{T_\cR, T_\cL, T_\cA,  T_\cK}$}\label{Sec.3.TG}
 The natural  (extended)  tangent spaces to the  groups $\cR,\cL,\cK,\cA$ are defined via the modules of derivations
  and endomorphisms, see \S\ref{Sec.Preliminaries.Aut.Der} and \cite{Mond-Nuno}.
 \beq
 T_\cR:=Der_X,\quad\quad\quad    T_\cL:=Der_Y,\quad\quad\quad  T_\cK:= End_{R_X}(\RpX)\oplus T_\cR,\quad\quad\quad T_\cA:=T_\cL\oplus T_\cR.
 \eeq
As $R_Y\!=\!\k[[y]] ,  \k\{y\} ,
   \k \bl y\br,$\!   one can write explicitly $Der_Y\!=\!R_Y\bl \frac{\di}{\di y_1},\dots,\frac{\di}{\di y_p}\br,$\! i.e.
   \!$\xi_Y\!=\!\sum \xi_{Y,i}\frac{\di}{\di y_i},$\! for
   $ \xi_{Y,i}\!\in\! R_Y.$

The name ``tangent space" is justified by the relation $T_\cG\rightleftarrows\cG$, see \S\ref{Sec.Preliminaries.Exp.Ln}, and \cite[\S3]{BGK.20}.

These tangent spaces act on the (extended) space of maps as in the classical case:
\beq
T_\cR\circlearrowright \RpX \quad\quad\quad{\rm by}\quad\quad\quad \xi_X(f):=\xi_X(f_1,\dots,f_p)=(\xi_X(f_1),\dots,\xi_X(f_1))\in \RpX,
\eeq
\[
T_\cL\!\!: \  \cm\cdot \RpX\to \RpX\quad\quad\quad{\rm by}\quad\quad\quad \xi_Y(f):=\xi_Y(y)|_f=(\xi_{Y,1}|_f,\dots,\xi_{Y,p}|_f)\in \RpX.
\]
The action $T_\cR\circlearrowright \RpX$ is $\k$-linear. This embeds $T_\cR=Der_X\sset End_\k(\RpX)$  as an $R_X$-submodule.
 The map $T_\cL:  \cm\cdot \RpX\to \RpX$ is by power series, it is neither additive nor $\k$-multiplicative.

Given a map $f\in \cm\cdot \RpX,$ one gets the (extended) image tangent space
 $T_\cG f\sseteq \RpX$. For our groups:
\beq
T_\cR f=Der_X(f),\quad\quad
T_\cL f=f^\#(R_Y)^{\oplus p},\quad\quad
T_\cA f=T_\cR f+T_\cL f,\quad\quad
T_\cK f=T_\cR f+(f)\cdot \RpX.
\eeq
\beR
\bee[\bf i.]
\item Classically these $T_\cG f$ are called ``the extended tangent spaces",  with the notations from differential geometry, \cite[pg.72]{Mond-Nuno}:
\beq\label{Eq.Classical.Notations.for.TG}
 R_X^{\oplus p}=\theta(f),\quad\quad\quad T_{\cR^e} f=tf(\theta_n),\quad\quad \quad T_{\cL^e} f=\omega f(\theta_p) , \quad\quad\quad
   T_{\cK^e} f=tf(\theta_n)+f^*\cm_p\cdot\theta(f).
\eeq
We avoid the notation  $T_{\cG^e}$ to prevent any confusion with the filtered tangent spaces,  $T_{\cG^{(j)}},$ \S\ref{Sec.3.Filtrations.on.G.TG}.
\item
 In zero characteristic, $\k\supseteq\Q,$ these ``image tangent spaces" coincide with the
 ``tangent spaces to the orbit", $T_\cG f=T_{\cG f}$. In positive characteristic the two notions differ, see \S1.5.2 of \cite{BGK.20}.

\item Note that $T_\cR f, T_\cK f\sset \RpX$ are $R_X$-submodules, while $T_\cL f, T_\cA f\sset \RpX$
 are only $R_Y$-submodules.
 Here $R_Y$ acts by substitution,
$q(y_1,\dots,y_p)\cdot T_\cA f:=q(f_1,\dots,f_p)\cdot T_\cA f\sset\RpX$, see \eqref{Eq.3.R_Y.action.on.R_X}.

\item The actions $\Phi_X,\!\Phi_Y\!\circlearrowright\!\Maps$ commute, see \eqref{Eq.A.action}. The derivations do not commute:
\begin{multline}\label{Eq.xi_x.xi_y.don't.commute}
\xi_X(\xi_Y(y)|_f)=\xi_X(\xi_{Y,1}|_f,\dots,\xi_{Y,p}|_f)=(\sum\nolimits_i \di_{y_i} \xi_{Y,1}|_f\cdot \xi_X(f_i), \dots, \sum\nolimits_i \di_{y_i} \xi_{Y,p}|_f\cdot \xi_X(f_i))
=\\=
\sum\nolimits_i \di_{y_i} \xi_{Y}(y)|_f\cdot \xi_X(f_i)\neq \xi_Y(y)|_{\xi_X(f)}.
\end{multline}
To pass from $\xi_X\cdot\xi_Y$ to $\xi_Y\cdot\xi_X$ we aply lemma \ref{Thm.Thom.Levine.alg}. This is used in \S\ref{Sec.A.Filtration}.
\eee
\eeR

\subsection{Filtrations on the group $\cG$ and on the tangent space $T_\cG$}\label{Sec.3.Filtrations.on.G.TG}
 Fix an ideal $I\sseteq \cm\sset R_X$ and take the filtration by powers of this ideal, $M_\bullet:=I^\bullet\cdot \RpX\sseteq \RpX$.
 Thus $M_\bullet\sset \Maps$ consists of maps that `vanish to a given order" on the locus $V(I)\sset X$.

\subsubsection{}\label{Ex.L.action.Filtration}
For $\cG= \cR$ or $\cL$ we get the corresponding filtrations on the groups/tangent spaces, as in \S\ref{Sec.Preliminaries.Aut.Der}.iii:
\beq\label{Eq.filtration.on.G.TG}
\cG^{(j)}:=\{g\in \cG|\ g\cdot M_d=M_d\ \text{ and } \ g|_{\quots{M_d}{M_{d+j}}}=Id|_{\quots{M_d}{M_{d+j}}},\ \forall d\ge1\}.
\eeq
\[
T_{\cG^{(j)}}:=\{\xi\in T_\cG|\quad \xi(M_d)\sseteq M_{d+j},\quad  \forall d\ge1\}.
\]
For the right equivalence   we have $T_{\cR^{(j)}}=T_{\cR}\cap End^{(j)}_\k(\RpX)$ and  $ T_{\cR}= T_{\cR^{(-1)}}\supsetneq T_{\cR^{(0)}}$,
  cf. \S\ref{Sec.Preliminaries.Aut.Der}. 

For the left equivalence we have
 $\cL^{(j)}=\{\Phi_Y\in \cL|\ \Phi_Y(y)-y|_{I^d\cdot \RpX}\in I^{d+j}\cdot \RpX\}$,
 i.e. $\Phi_Y(y)=y+\phi(y)$, where $\phi(I^d)\sseteq I^{d+j}$.
 In particular $\cL^{(0)}=\cL$.

 The tangent space filtration is $T_{\cL^{(j)}}=\{\xi_Y\in Der_Y|\ \xi_Y(y)|_{I^d}\sseteq I^{d+j}\cdot \RpX,\ \forall\ d\ge1\}$.
 In particular:
 \beq
T_{\cL}\supsetneq  T_{\cL^{(0)}}=(y)\cdot T_\cL\supsetneq T_{\cL^{(1)}}=(y)^2\cdot T_\cL.
 \eeq
  More generally,    $T_{\cL^{(j)}}=\cb_{j}\cdot T_{\cL}$, for the largest ideal $\cb_j\sset R_Y$ satisfying: $f^\#(\cb_{j})\sseteq I^{d+j}$ for all $d\in \N$ and
  all $f\in I^d\cdot \RpX$. In particular, the tangent image is $T_{\cL^{(j)}}f=f^\#(\cb_{j}\cdot R^{\oplus p}_Y)\sseteq I^{d+j}\cdot\RpX$.

We have the obvious inclusions:  $(y)^{j+1}\sseteq \cb_{j }\sset R_Y$, and $(f)^{j+1}\sseteq f^\#(\cb_{j})\cdot R_X$ when $(f)\sseteq I$.
 These are equalities when $\k\supseteq\Q$ is a field, but can be proper if $\k$ is a finite field or a ring.
 \bex\label{Ex.Pathology.of.TL}
\bee[\bf i.]
\item
Suppose $R_X$ is one of $\k[[x]]$, $\k\{x\}$, $\k\bl x\br$, where $\k$ is an infinite field. Suppose $I=\cm=(x)\sset R_X$. Take a map $f\in (x)\cdot \RpX$.
 Then:
 \[
 T_{\cR^{(j)}}=(x)^{j+1}\cdot T_\cR,\quad\quad
T_{\cR^{(j)}}f=(x)^{j+1}\cdot Der_X(f),\quad\quad
T_{\cL^{(j)}}=(y)^{j+1}\cdot T_\cL,\quad\quad
T_{\cL^{(j)}}f= f^\#(y)^{j+1}\cdot T_\cL f\sset \RpX.
 \]
In the classical notations, cf. \eqref{Eq.Classical.Notations.for.TG}: $T_{\cR^{(j)}}f=tf(\cm^{j+1}_n\cdot \Theta_n)$
 and $T_{\cL^{(j)}}f=\omega f(\cm^{j+1}_p\cdot \Theta_p)$.

\item
 Consider $Maps((\k,o),(\k^2,o))$, where $\k$ is a field with two elements, and the rings are complete. Take the filtration $(x)^\bullet\cdot R^{\oplus 2}_X$.
  Take the derivation  $\xi_Y=q(y)\di_{y_1}+0\cdot \di_{y_2}\in Der_Y$, where $q(y)=y_1y_2(y_1+y_2)\in \k[[y_1,y_2]]$.
   Then $q((x)^i,(x)^i)\sseteq (x)^{3i+1}$. Thus $\xi_Y\in T_{\cL^{(3)}}$, even though $q(y)\not\in (y)^4\sset \k[[y_1,y_2]].$

\item Let $\k=\tilde\k[[t]]$  for a field $\tilde\k\supseteq\Q$. Take $I=(t,\cm)$, then $T_{\cL^{(j)}}=(t,y)^{j+1}T_\cL\supsetneq (y)^{j+1}T_\cL$.
\eee \eex

\beR\label{Rem.Filtrations.on.TR.TL}
\bee[\bf i.]
\item In the definition \eqref{Eq.filtration.on.G.TG} it is important to assume $d\ge1$, rather than $d\ge0$.
 In fact the group $\cL$ does not act on the whole $\RpX$, thus the condition ``$g|_{\quots{M_0}{M_{j}}}=Id|_{\quots{M_0}{M_{j}}}$" makes no sense.
 One could restrict  to $\quots{\cm\cdot M_0}{M_{j}}$, but then, for $\sqrt{I}\ssetneq\cm$,    the condition
 $g|_{\quots{\cm\cdot M_0}{M_{j}}}=Id|_{\quots{\cm\cdot M_0}{M_{j}}}$  would imply $\cL^{(j)}=\{e\}$. Similarly, the condition
 $\xi(\cm\cdot M_0)\sseteq M_{j}$ would imply $T_{\cL^{(j)}}=(0)$.
\item
The relation $(y)\cdot T_\cL=T_{\cL^{(0)}}$ motivates the notation $T_{\cL^{(-1)}}:=T_\cL$.
 This should be used   with caution, as $T_\cL(I^d)\not\sseteq I$ for any $d$.
\eee
\eeR

\subsubsection{}\label{Sec.3.Filtrations.on.A.TA}
For the groups $ \cA,\cK$ we define the filtration  in a special way:
\bei
\item $\cA^{(j)}:=\cL^{(j)}\times \cR^{(j)}$ and $T_{\cA^{(j)}}:=T_{\cL^{(j)}}\oplus T_{\cR^{(j)}} $.

\item $\cK^{(j)}:=\cC^{(j)}\rtimes \cR^{(j)}$ and $T_{\cK^{(j)}}:=T_{\cC^{(j)}}\oplus T_{\cR^{(j)}}$. Here
 $T_{\cC^{(j)}}:=Mat_{p\times p}(I^j)$.
\eei
In particular $T_\cK=T_{\cK^{(-1)}}$. As in remark \ref{Rem.Filtrations.on.TR.TL}.ii we denote $T_{\cA^{(-1)}}:=T_\cA$. Again care is needed, as
 $T_\cA (I^d\cdot \RpX)\not\sseteq I\cdot \RpX$ for any $d$.

\

One would like to use the general definition \eqref{Eq.filtration.on.G.TG}, for the groups $\cA,\cK$, i.e. to define
\beq
\cA^{(j)} =(\cL\times\cR)^{(j)} \quad \text{ and  } \quad T_{\cA^{(j)}}=(T_{\cR}\oplus T_{\cL})^{(j)},\quad
\text{ and so on.}   
\eeq
 In this definition $\cA^{(j)},T_{\cA^{(j)}}$ become non-explicit/complicated.
But the two versions often coincide.   

\bel\label{Thm.TA.TK.distinct.filtrations} Suppose  $j\ge1$,  and   $\k$ is an infinite field, and $I^\bullet\neq(0)$ for each $\bullet\in \N$. Then:
\bee[\bf 1.]
\item $(\cL\times\cR)^{(j)}=\cL^{(j)}\times \cR^{(j)}$, i.e. $(Aut_Y\times Aut_X )^{(j)}=Aut^{(j)}_Y\times Aut^{(j)}_X$.
\item
  $(T_{\cL}\oplus T_{\cR})^{(j)}=T_{\cL^{(j)}}\oplus T_{\cR^{(j)}}$, i.e.
  $(Der_Y\oplus Der_X)^{(j)}=Der^{(j)}_Y\oplus  Der^{(j)}_X$.
\eee
\eel
\bpr The inclusions $(\cL\times\cR)^{(j)}\ge \cL^{(j)} \times \cR^{(j)}$ and  $(T_{\cL}\oplus T_{\cR})^{(j)}\supseteq T_{\cL^{(j)}}\oplus T_{\cR^{(j)}}$ are
  obvious. We prove the parts $\le$ and $\sseteq$.

Below we use the property: there exists $f\in I$  such that $f^d\in I^d\smin I^{d+1}$ for each $d\ge1$. Indeed, otherwise we
 would get $I^d\sseteq I^{d+1}$ for $d\gg1$,   as $I$ is finitely generated. But then $I^d=0$ by Nakayama.
 \bee[\bf 1.]
 \item
Suppose   $(\Phi_Y,\Phi_X)\in (\cL\times\cR)^{(j)}$. As $\Phi_Y\in \cL=\cL^{(0)}$, we get: $\Phi_X\in \cR^{(0)}$.
 Suppose $\Phi_Y\in  \cL^{(l)}\smin \cL^{(l+1)}$, then, after a $\k$-linear reshuffling of $y$-coordinates, we can assume:
 $\Phi_Y(y_1,0,\dots,0)=(q(y_1),\dots).$ Here:
 \beq
 q(y_1)=y_1+a\cdot y^{\tl+1}_1+h.o.t. \quad \text{(for some $\tl\le l$)},  \quad   a\in \k^\times, \quad h.o.t.\in (y_1)^{l+2} \quad \text{and}  \quad
  a\cdot I^{\tl+1}\sseteq I^{l+1}.
  \eeq
   As $I^\bullet\neq0$ for each $\bullet,$ we get  $\tl=l.$

  For each $f\in I$ define the difference $\De_f:=[\Phi_Y(y_1,0,\dots,0)|_f]_1-\Phi_X(f)=q(f)-\Phi_X(f).$
  By our assumption $\De_{f\circ \Phi^{-1}_X} \in I^{ord(f)+j}.$ As $\Phi _X\in \cR^{(0)},$ we get: $\De_{f } \in I^{ord(f)+j}.$

\bei
\item The case $l=0$. Then $\Phi_Y\in  \cL^{(0)}\smin \cL^{(1)}$ and  $\Phi_X\in  \cR^{(0)}\smin \cR^{(1)}.$
 For each $f\in I$ and $d\in \N$ we get
\beq
I^{d\cdot ord(f)+j}\ni f^{d-1}\cdot \De_f-\De_{f^d}=f^{d-1}\cdot \Phi_X(f)-\Phi_X(f)^d+(h.o.t.),\quad \text{where}\quad  (h.o.t)\in I^{(d+1)\cdot ord(f)}.
\eeq
As this holds for each $f\in I,$ and the field $\k$ is infinite, we get: $\Phi_X\in \cR^{(1)},$ i.e. a contradiction.


\item The case $l\ge1$. For any $f\in I$, $c\in \k^\times$ we have:
\beq
\De_{c\cdot f}-c\cdot \De_f=q(c\cdot f)-c\cdot q(f)=a\cdot (c^{l+1}-c)\cdot f^{l+1}+h.o.t.\in I^{ord(f)+j}.
\eeq
As $\k$ is infinite there exists $c\in \k^\times$ for which $c^{l+1}-c\in \k^\times.$ Thus for each $f\in I$ and $l>0$ we get: $f^{l+1}\in I^{ord(f)+j}$.
 Therefore $l\ge j$, hence $\Phi_Y\in \cL^{(j)}$, and thus $\Phi_X\in \cR^{(j)}$.
\eei

\item
Suppose  $(\xi_Y,\xi_X)\in (T_{\cL}\oplus T_{\cR})^{(j)}$, thus $(\xi_Y,\xi_X)(I\cdot \RpX)\sseteq I^{1+j}\cdot \RpX$.
 If $\xi_Y\in T_{\cL}\smin T_{\cL^{(0)}}$ then $\xi_Y(I)\not\sseteq \cm$.
 Then for $f\in (I\cap \cm^2)\cdot \RpX$ we get: $(\xi_Y,\xi_X)(f)\not\in\cm\cdot\RpX.$
 Therefore necessarily $\xi_Y\in T_{\cL^{(0)}},$ and then  $\xi_X\in T_{\cR^{(0)}}$.

 Suppose $\xi_Y\in T_{\cL^{(l)}}\smin T_{\cL^{(l+1)}}$, then after a $\k$-linear reshuffling of $y$-coordinates we can assume $\xi_Y(y_1,0,\dots,0)=(q(y_1),\dots).$
  Here $q(y_1)=a\cdot y^{\tl+1}_1+h.o.t.$  (for some $\tl\le l$), \quad  $a\in \k^\times,$ \quad $h.o.t.\in (y_1)^{l+2}$ and $a\cdot I^{\tl+1}\sseteq I^{l+1}.$
 Therefore $\tl=l.$
\bei
\item The case $l=0$. For any $f\in I$ and $d\in \N$ we have
\beq
[(\xi_Y,\xi_X)|_{f^d,0,\dots,0}]_1-d\cdot f^{d-1}\cdot [(\xi_Y,\xi_X)|_{f,0,\dots,0}]_1=q(f^d)-d\cdot f^{d-1}\cdot q(f)=a(1-d)f^d+h.o.t\in I^{ord(f^d)+j}.
\eeq
Choose $f\in I\smin I^2$ satisfying $f^2\in I^2\smin I^3$ to get the contradiction.

\item The case $l\ge1$. For any $f\in I$ and $c\in \k$ we have
\beq
[(\xi_Y,\xi_X)|_{c\cdot f,0,\dots,0}]_1-c\cdot [(\xi_Y,\xi_X)|_{f,0,\dots,0}]_1=q(c\cdot f)-c\cdot   q(f)=a(c^{l+1}-c)f^{l+1}+h.o.t\in I^{ord(f)+j}.
\eeq
As $\k$ is infinite there exists $c\in \k$ such that $c^{l+1}-c\neq0$. Therefore we get  $l\ge j$, hence $\xi_Y\in T_{\cL^{(j)}}$,
 and thus $\xi_X\in T_{\cR^{(j)}}$. \epr
\eei
\eee

\beR\!\!By  similar arguments  one \!can prove, for the contact group $\cK\!=\!\cC\!\rtimes\!\cR,$ and    \!an \!infinite field $\k:$
\[
(\cC\rtimes \cR)^{(j)}=\cC^{(j)}\rtimes \cR^{(j)},\quad
 (GL(p,R_X)\rtimes Aut_X)^{(j)}=GL^{(j)}(p,R_X)\rtimes Aut^{(j)}_X, \quad
(T_{\cC}\oplus T_{\cR})^{(j)}=T_{\cC^{(j)}}\oplus T_{\cR^{(j)}}.
\]
\eeR

\subsection{The closures in the filtration topology}\label{Sec.3.Filtration.Top.Closures}
For our groups $\cG= \cR,\cK,\cA$ of \S\ref{Sec.3.Groups.R.K.A} we have the closures of the orbits in the filtration topology (for $M_\bullet=I^\bullet\cdot \RpX$):
 \beq
\overline{T_{\cG^{(j)}}f}=\cap_{\bullet\ge1} (T_{\cG^{(j)}}f+M_\bullet),\quad\quad\quad\quad
\overline{\cG^{(j)} f}=\cap_{\bullet\ge1} (\cG^{(j)} f+M_\bullet).
 \eeq

\bel\label{Thm.Filtration.Closures.TR.TK}
\bee
\item   Any $R_X$-submodule of $\RpX$ is closed in $I^\bullet$-topology.
\\ In particular, the tangent spaces $T_{\cR^{(j)}} f$, $T_{\cK^{(j)}} f\sseteq \RpX$ are closed.
\item  The orbits $\cR^{(j)}f$, $\cK^{(j)}f$ are closed.
\eee
\eel
\bpr\bee
\item
 Let $M\sseteq\RpX $ then $M+I^d\cdot \RpX\supseteq \overline{M}$ for any $d\ge1$. Therefore  $M+(I^d\cdot \RpX)\cap \overline{M} \supseteq \overline{M}$.
 Applying the Artin-Rees lemma to the (finitely generated, \S\ref{Sec.Preliminaries.Rings.Germs}.viii) submodule $\overline{M}\sset \RpX$
   we get: $M+I\cdot \overline{M}\supseteq \overline{M}$. Then  Nakayama gives $M=\overline{M}$.

\item (The $\cK$-case.) Suppose $f+g\in \overline{\cK^{(j)} f}$ for  a perturbation $g\in \cm\cdot \RpX.$
 Thus $f+g\in \overline{(\cK^{lin})^{(j)} f}.$
  Thus the system of implicit
 function equations $f(\Phi(x))=U\cdot (f+g)$ has an order-by-order solution.
  By the theorem of Pfister-Popescu, \S\ref{Sec.Preliminaries.Basic.Tools}.v., it has a formal solutions,
  $(\hat\Phi,\hat{U})\in\widehat{(\cK^{lin})^{(j)}}$.
\\
For the rings  $\quots{\k\{x\}}{J},\quots{\k\bl x\br}{J}$ apply the Artin approximation  to
   get an ordinary solution,   $(\Phi,{U})\!\in\!(\cK^{lin})^{(j)}.$
\epr
\eee

\beR
\bee[\bf i.]
\item Part 1 is a pure commutative algebra, it holds for Noetherian local rings, provided all the objects are defined.
\item
Below we use the simple observation: to verify the inclusion $\overline{\cG^{(j)} f}\supseteq \{f\}+M_d$
  it is enough to verify $\cG^{(j)} f+M_{\bullet+1}\supseteq \{f\}+M_\bullet$ for every $\bullet\ge d$.
 And similarly for $\overline{T_{\cG^{(j)}}f}\supseteq M_d$.
\item
 The case of $T_\cA f$ is more delicate, see lemma \ref{Thm.Filtration.Closures.TA}.
\eee
\eeR

\subsection{Related derivations and related automorphisms. The ``naturality" of  vector fields integration}\label{Sec.3.Related.Derivations.Autos}
  Given a map of ($C^\infty$ or real/complex-analytic)
 manifolds, $f: X \to{Y}$, take some vector fields $\xi_X\in \Theta_X$, $\xi_Y\in \Theta_Y$,
   and the corresponding time-one flows, i.e. the automorphisms $\Phi_X\circlearrowright X$,
  $\Phi_Y\circlearrowright Y$.
\bel \
[Thom-Levine,   e.g.  \cite[Lemma 2.6]{Mond-Nuno}] If the vector fields are related, i.e.  $\xi_Y\circ f=df\circ \xi_X$, then so are the corresponding automorphisms,
 $\Phi_Y\circ f=f\circ \Phi_X$.
\eel

 We need the algebraic version of this  and other results on vector fields.
Below we use  the exponential and logarithmic maps of  \S\ref{Sec.Preliminaries.Exp.Ln}.

Take  a  local filtered $\k$-algebra  $(R _X,\cm , I^\bullet).$ Assume   $\k\supseteq\Q$ and $R_X$ is complete \wrt the filtration.
 Take the filtration of $\RpX$ by $M_\bullet:=I^\bullet\cdot \RpX$. Put $M_\infty:=0$. The map $exp$ of \eqref{Eq.Exp.Ln.full.series} is well defined.

\bel\label{Thm.Thom.Levine.alg}
 Fix some elements $(\xi_Y,\xi_X)\in T_{\cA^{(l)}},$        $f\in \cm\cdot \RpX$,  for some    $l\ge1,$ $1\le d\le\infty,$ and
 $w_d \in M_{d}$.
\bee[\bf 1.]
\item $\xi_Y (e^{\xi_X}f)=e^{\xi_X}\cdot \xi_Y(y)|_f$.
\item   If \quad$\xi_X(f)\stackrel{mod\ M_{d}}{\equiv}\xi_Y(y)|_f $ \quad
 then \quad $\xi^i_X(\xi^j_Y(y)|_f)\stackrel{mod\ M_{d+l}}{\equiv}\xi^{i+j}_Y(y)|_f$ \quad for any $i+j\ge2$.
\\ In particular, $e^{\xi_X}f-e^{\xi_Y}f\in M_{d}$.
\item(algebraic ``naturality")    $\xi_Y(y)|_f-\xi_X(f)\stackrel{mod\ M_{d+l}}{\equiv}   w_d  $ \ if and only
 if \ $e^{\xi_Y}f-e^{ \xi_X}f\stackrel{mod\ M_{d+l}}{\equiv}  w_d$.
\eee
\eel
\bpr
\bee[\bf 1.]
\item The direct check:  $\xi_Y (e^{\xi_X}f)=\xi_Y(y)|_{e^{\xi_X}f}=\xi_Y(y)|_{f\circ e^{\xi_X}}=e^{\xi_X}\cdot \xi_Y(y)|_f$.
\item Induction on $i$. The case $i=0$ is trivial. Assuming the statement for $i$ we verify it for $i+1$:
\begin{multline}
\xi^{i+1}_X(\xi^j_Y(y)|_f)=\xi_X\Big(\xi^i_X(\xi^j_Y(y)|_f)\Big)=\xi_X\Big( \xi^{i+j}_Y(y)|_f\ mod\ M_{d}\Big)\stackrel{mod\ M_{d+l}}{\equiv}
\xi_X\Big( \xi^{i+j}_Y(y)|_f\Big)=\\
\sum\nolimits_q \di_{y_q}\Big( \xi^{i+j}_Y(y)\Big)|_f\cdot \xi_X(f_q)
\stackrel{mod\ M_{d+l}}{\equiv}
\sum\nolimits_q \di_{y_q}\Big( \xi^{i+j}_Y(y)\Big)|_f\cdot \xi_Y(y)_q|_f=\xi^{i+j+1}_Y(y)|_f.
\end{multline}
The statement  $e^{\xi_X}f-e^{\xi_Y}f\in M_{d}$ is now obtained by Taylor expansion of the exponents.
\item {\bf The part $"\Rrightarrow"$.} \ $e^{\xi_Y}f- e^{\xi_X}f=w_d+\sum_{j\ge 2}\frac{\xi_Y^j(y)|_f-\xi_X^j(f)}{j!}\equiv w_d\ mod\ M_{d+l}$.
 For the last equality we use part 2.
\\
{\bf The part $"\Lleftarrow"$.}  \   $\xi_Y(y)|_f-\xi_X(f)=ln(e^{\xi_Y})(y)|_f-ln(e^{\xi_X})f\stackrel{!}{=}ln(e^{-\xi_X}e^{\xi_Y})(y)|_f=
 -\sum_{j\ge 1}\frac{(1-e^{-\xi_X}e^{\xi_Y})^j}{j}f= \sum_{j\ge 1}\frac{(1-e^{-\xi_X}e^{\xi_Y})^{j-1}}{j}e^{-\xi_X}(w_d\ mod\ M_{d+l})
  \stackrel{mod\ M_{d+l}}{\equiv}  w_d$.
 For the transition "!" we use the commutativity $e^{\xi_X}e^{\xi_Y}=e^{\xi_Y}e^{\xi_X}$.
\eee
\beR
For the subsequent applications (\S\ref{Sec.R.Filtration}, \S\ref{Sec.K.Filtration}, \S\ref{Sec.A.Filtration}) one needs the strongest Thom-Levine statement.
 In Part 3 the conclusion cannot be strengthened to the statement ``$   (mod\ M_{d+l+1})$."
  For example, define the map $(\k^1,o)\to(\k^1,o)$ by $f(x)=x+x^{d-1}\in\k[[x]]$. Here $\k\supseteq\Q$, the ring $R_X=\k[[x]]$ is filtered by $(x)^\bullet$,
  the ring $R_Y=\k[[y]]$ is filtered by $(y)^\bullet$.
\\ Take the derivations $\xi_X=x^2\di_x$,  $\xi_Y=y^2\di_y$. One has $\xi_X(f)-\xi_Y(y)|_f\stackrel{mod\ x^{2d-2}}{\equiv} (d-3)x^d$.
 But   $e^{\xi_X}(x)=\frac{x}{1-x}$ and therefore $e^{\xi_X}(f)-e^{\xi_Y}(f)=f(\frac{x}{1-x})-\frac{f(x)}{1-f(x)}=(d-3)x^d\big(1+\frac{3d-2}{2}x+\dots\big)$.
  Thus $e^{\xi_X}(f)-e^{\xi_Y}(f)\not\equiv (d-3)x^d\ mod\ (x)^{d+2}$.


\eeR
\beR\label{Rem.Thom-Levine.without.Q}
In the proof of part 3 we needed the exponential/logarithmic expansions only up to order $\lceil \frac{d-ord(f)}{l}\rceil+1$. Therefore this lemma can be
 used even when the full exp/log maps do not exist, e.g. when $\k\not\supseteq\Q$ or $\RpX$ is not $M_\bullet$-complete.
  Namely, fix some $(\xi_Y,\xi_X)\in T_{\cA^{(l)}}$, $w_d\in M_{d}$,   and suppose the assumptions $jet_N Exp$,  $jet_N Ln$ of
   \S\ref{Sec.Preliminaries.Exp.Ln} hold for $N=\lceil \frac{d-ord(f)}{l}\rceil+1$. Then we get the statement:
  \beq
  \xi_Y(y)|_f-\xi_X(f)\stackrel{mod\ M_{d+l}}{\equiv}   w_d \quad\quad  \text{ if and only  if } \quad\quad
   jet_N(e^{-\xi_X}) jet_N(e^{\xi_Y})f\stackrel{mod\ M_{d+l}}{\equiv} f+w_d.
  \eeq
\eeR

\subsection{Critical/singular/instability loci and the annihilator  of $\pmb{ T^1_\cG f}$}\label{Sec.3.Annihilators}
Fix a group  $\cG= \cR,\cA,\cK$ and a map $f\in  \Maps=\cm\cdot \RpX.$ Take the image tangent space $T_\cG f\sseteq \RpX$, \S\ref{Sec.3.TG}.
 Take the quotient module $T^1_\cG f:=\quots{\RpX}{T_\cG f}$. (In \cite{Kerner.Unfoldings} we show that this is the tangent
  space to the $\cG$-miniversal unfolding, when $f$ is $\cG$-finite.)
The standard way to measure
 how large is this tangent space goes via the support of $T^1_\cG f$, i.e. the annihilator of this quotient module
 \beq\label{Eq.Annihilators.Definition}
 \ca_\cG:=Ann[ T^1_\cG f]=\{q\in R_X|\ q\cdot \RpX\sseteq T_\cG f\}\sseteq R_X.
 \eeq
  This is the largest ideal satisfying: $T_\cG f\supseteq \ca_\cG \cdot \RpX$.

Fix a filtration $I^\bullet\sset R_X$   and the corresponding filtrations $\cG^{(\bullet)}$, $T_{\cG^{(\bullet)}}$, see \S\ref{Sec.3.Filtrations.on.G.TG}.
 One gets the filtered annihilators $\ca_{\cG^{(j)}}:=Ann [T^1_{\cG^{(j)}} f]\sseteq R_X$ for $j\ge0$.

  \bed
  \bee
\item  For $\ca_\cG\sset R_X$  the scheme-germ $V(\ca_\cG)\sseteq X$ is called the critical ($Crit(f)$),
  resp. the singular ($Sing(f)$), resp. the instability locus of $f$.
\item   The map $f$ is called $\cG$-finite if $\sqrt{\ca_\cG}\supseteq\cm$, i.e. geometrically $V(\ca_\cG)=o\in X$ or $V(\ca_\cG)=\empty.$
\eee
\eed

If $X$ is smooth, $\k=\bk,$ and $f_1,\dots,f_p\in R_X$ is a regular sequence then (set-theoretically) $Sing(f)=Sing(V(f))$, the singular locus of the complete intersection.
 Otherwise the two objects can differ.

\beR
Classically  the Fitting (determinantal) ideals of the modules $T^1_\cR f$, $T^1_\cK f$ were used. We use the annihilator ideals in \eqref{Eq.Annihilators.Definition}
 for two reasons:
 \bei
 \item These annihilators appear naturally in the criteria for group-orbits, \S\ref{Sec.R.IFT}, \S\ref{Sec.K.IFT},  \S\ref{Sec.A.IFT}.
 \item
 The ideals $ Fitt_\bullet(T^1_\cA f)$ are
 not immediately defined, as $T^1_\cA f$ is not an $R_X$ module, while as  $R_Y$-module it is not finitely generated.
 \eei
\eeR

\subsubsection{} Below we discuss the annihilators $\ca_\cR,\ca_\cK$. The annihilator $\ca_\cA$ is much more delicate, see \S\ref{Sec.TA}.
\bex\label{Ex.Annihilators.R.K.simple}
 Let  $R$  be one of $\quots{\k[[x]]}{J}$, $\quots{\k\{x\}}{J}$, $\quots{\k\bl x\br}{J}$,   see \S\ref{Sec.Preliminaries.Rings.Germs}.ii.
\bee[\bf i.]
\item For $p\!=\!1$ and $R$   one of $\k[[x]],\k\{x\}, \k\bl x\br$ the annihilators are just the Jacobian  and the Tjurina ideals,
\[
T^1_\cR f=\quots{R}{Jac(f)},\quad\quad
\ca_\cR=Jac(f),\quad\quad
T^1_\cK f=\quots{R}{(f)+Jac(f)}, \quad\quad
\ca_\cK=Jac(f)+(f).
\]
Here (when $\k$ is not a field) the derivatives are taken \wrt $x$-variables only.

For $J\neq 0$ one gets the Bruce-Roberts versions of the Jacobian/Tjurina ideals, \cite{Bruce-Roberts}.
\item
For $p\ge2$ the annihilators are not computed easily, but one has immediate bounds via the determinantal ideals of the generating matrix $[\xi_i(f_j)]$ of $T_\cR f$.
 (For the regular rings $\k[[x]]$, $\k\{x\}$, $\k\bl x\br$, this is just the matrix of partials,  $[\di_{x_i}f_j]$.)
  One has (see e.g. \cite[\S20]{Eisenbud-book}):
\beq
I_p[\xi_i(f_j)]\sseteq \ca_\cR\sseteq \sqrt{I_p[\xi_i(f_j)]},\quad\quad\quad\quad
I_p[\xi_i(f_j)]+(f)\sseteq \ca_\cK\sseteq \sqrt{I_p[\xi_i(f_j)]+(f)}.
\eeq
\item In particular one gets: $\ca_\cR+(f)\sseteq \ca_\cK\sseteq \sqrt{\ca_\cR+(f)}$. Hence $Sing(f)=Crit(f)\cap V(f)$.
\item Let $p\ge2$.  In many cases one has $\ca_\cR\cdot \RpX\sseteq \cm\cdot T_\cR f$. For example, this holds if
 the determinantal ideal $I_p[\xi_i(f_j)] \sset R_X$ is radical, as then $\ca_\cR=I_p[\xi_i(f_j)]$.

\item Take the filtration $I^\bullet\sset R_X$. One has $I^{j+1}\cdot\ca_\cG\sseteq \ca_{\cG^{(j)}}\sseteq \ca_{\cG}$.

\item (The case of unfolding) Take a map  $f:(X\times\k^r_u,o)\to (\k^p_y\times\k^r_u,o)$, $(x,u)\to (f_o(x)+g(x,u),u)$, where $(g(x,u))\sseteq (x)\cdot(u)$.
 Then $\ca_\cK (f)=\ca_\cK(f_o)+(u)$ and $\ca_{\cK^{(0)}} (f)=\ca_{\cK^{(0)}}(f_o)+(u)$.

\eee
\eex

\subsubsection{Maps of finite singularity type}\label{Sec.3.Finite.Sing.Type}
The following criterion is well known for $\k=\C$ and $X$-ICIS, see e.g. \cite{Mond-Montaldi}, \cite[pg.224]{Mond-Nuno}.
\bel\label{Thm.K.finite.maps} Let the pair $(R_X,R_Y)$ be one of $(\quots{\k[[x]]}{J},\k[[y]])$, $(\quots{\k\{x\}}{J},\k\{y\})$, $(\quots{\k\bl x\br}{J},\k\bl y\br),$ see \S\ref{Sec.Preliminaries.Rings.Germs}.ii.
The following conditions are equivalent for a  map $f:X\to (\k^p,o)$.

\bee[\!\!\bf 1.\!]
\item The restriction of $f$ to $Crit(f)\!\to \!(\k^p,o)$ is a finite morphism, i.e. the $R_X$-module $\quots{R_X}{\ca_\cR}$ is a f.g. $R_Y$-module;

\item The $R_X$-module $T^1_\cR f$ is a f.g. $R_Y$-module;

\item  $T^1_\cK f$  is a f.g. $\k$-module;  \quad\quad  (If $\k$ is a field then $f$ is $\cK$-finite.)
\item $f^{-1}(o)\cap Crit(f)$ is a (fat) point over $Spec(\k)$ (i.e. $\sqrt{\ca_\cR+(f)}\supseteq(x)$).
\eee
\noindent
 Moreover, if   $T^1_\cK f\in mod\text{-}\k$ is generated by $\{v_\bullet\}\sset \RpX$, then $T^1_\cR f\in mod(R_Y)$ is generated  by $\{v_\bullet\}$.
 \eel
In this case $f$ is called \underline{a map of finite singularity type}.
\bpr
$\bf 1\Rrightarrow 2.$ $T^1_\cR f$ is a f.g. module over  $\quots{R_X}{\ca_\cR}$, and thus a f.g. module over $R_Y$.
\\
$\bf 2\Rrightarrow 3.$ If \!$T^1_\cR f$ is   f.g.   over  \!$R_Y,$ then $T^1_\cR f\!\otimes\!\quots{R_Y}{(y)}$ is f.g. over \!$\quots{R_Y}{(y)}\!=\!\k.$
 Now observe: \!$T^1_\cR f\!\otimes_{R_Y}\!\quots{R_Y}{(y)}\!=\!T^1_\cK f.$
\\
$\bf 3\Rrightarrow 4.$ As $T^1_\cK f$ is f.g. over $\k,$ it is a module over $\quots{R_X}{(x)^d}$ for some $d\gg1.$
 Therefore  $\sqrt{\ca_\cK}=\sqrt{\ca_\cR+(f)}\supseteq(x)$. Thus $f^{-1}(o)\cap Crit(f)$ is a fat point over $Spec(\k).$
\\
$\bf 4\Rrightarrow 1.$ By the assumption $\sqrt{\ca_\cR+(f)}\supseteq(x).$ Hence the map $Crit(f)\to (\k^p,o)$ is quasi-finite.
 The finiteness follows by the Weierstra\ss-finiteness, \ref{Sec.Preliminaries.Basic.Tools}.iii.
\epr

\section{Criteria for right orbits, ``$\cR f$ vs $T_\cR f$"}\label{Sec.R.Orbits}
Let $R_X$ be one of the rings $\quots{\k[[x]]}{J}$, $\quots{\k\{x\}}{J}$, $\quots{\k\bl x\br}{J}$, with $\k$ a local ring,  see \S\ref{Sec.Preliminaries.Rings.Germs}.ii.

\subsection{The $\cR$-implicit function theorem}\label{Sec.R.IFT}

     Take  two ideals $\ca\sseteq I\sseteq  R_X$, with $\ca\sseteq\cm^2$, and a map-germ $f\in I\cdot \RpX$.
   (The usual choice is: $I=R_X$ or $I=\cm$ or $I=\sqrt{(f)}$.)
    For the notation $I^{ord(f)-2}$ see \S\ref{Sec.Preliminaries.Rings.Germs}.v.

\bthe\label{Thm.IFT.R.case} Suppose the  ideal $\ca$ satisfies  $\ca^2\cdot I^{ord(f)-2}\!\cdot \RpX\sseteq\cm\cdot  \ca\cdot T_{\cR} f.$
\bee
\item (The case $J=0$)
 Then    $\{f\} + \ca\cdot T_{\cR} f \sseteq \cR f .$
 \item  (The case $J\neq0,$  assume the $jet_0$-condition of \S\ref{Sec.Preliminaries.Exp.Ln}.)
    Then    $\{f\} + \ca^2\cdot T_{\cR} f \sseteq \cR f .$
 \eee
\ethe
A remark: the assumption    $\ca^2\cdot I^{ord(f)-2}\!\cdot \RpX\sseteq\cm\cdot  \ca\cdot T_{\cR} f$ is preserved by $\cR$-equivalence.
\\
The automorphism $\Phi_X\!\in \!Aut_X$ (constructed in the proof) satisfies: $\Phi(x)\!-\!x\!\in\! \ca,$ resp. \!(for part 2) $\Phi(x)\!-\!x\!\in\! \ca^2.$
\bpr
\bee[\bf 1.]
\item  Take a perturbation $g\in \ca \cdot T_\cR f$. We want to resolve the condition $\Phi_X(f)= f+g $, where $\Phi_X\in \cR=Aut_X$.
 Fix some (finite sets of) generators  $\{\xi_i(f)\}$ of the $R_X$-module $\ca\cdot T_\cR f=\ca\cdot Der_X(f)$.
 Expand $g=\sum c^g_{i}\xi_j(f)$, here $c^g_{i}\in R_X$.

    We look for the coordinate change $\Phi_X\in Aut_X$ in the form $x\to x+\sum c_{i}\xi_i(x)$.
 Thus we want to resolve the equation
\beq\label{Eq.proof.of.R.TR.eq.to.resolve}
f(x+\sum   c_{i}\xi_i(x))=f(x)+\sum   c^g_i\cdot \xi_i(f), \quad  \text{ with the unknowns $\{c_*\}$}.
\eeq
 Taylor-expand the left hand side:
   $f(x)+\sum  c_{i}\xi_i(f)+H_{\ge2}(\{ c_{i}\xi_i(x)\})$.
   Here $H_{\ge 2}(z)\in \sum^{ord(f)}_{l\ge 2} I^{ord(f)-l}\cdot (z)^l\cdot \RpX$.
   As $\ca\sseteq I$ one has: $H_{\ge2}(\ca\cdot z)\in I^{ord(f)-2}\cdot (\ca\cdot z)^2\cdot R_X[[\ca\cdot z]][z]$.
    These are polynomials in $z$ whose coefficients are (formal/analytic/algebraic)
    power series in     $\ca\cdot z$, see \S\ref{Sec.Preliminaries.Rings.Germs}.iv. Therefore:
\beq\label{Eq.bound.2nd.order.terms}
H_{\ge2}(\{  c_{i }\xi_i(x)\})\in
 (\{c_*\})^2\cdot\ca^2\cdot I^{ord(f)-2}\cdot \RpX\sseteq
  (\{c_*\})^2 \cdot\cm\cdot  \ca\cdot T_\cR f.
\eeq
Thus we can present $H_{\ge2}( c_{i }\xi_i(x))=\sum  h_{i }(c_*)\xi_i(f)$, for some (fixed) power series
$h_{i }(z)\in (z)^2\cdot\cm\cdot R_X[[\ca\cdot z]][z]$.
 These (formal/analytic/algebraic) series $h_{i }$ satisfy: $h_{i }(R_X)\sset \cm$.

 Thus to resolve the equation \eqref{Eq.proof.of.R.TR.eq.to.resolve} it is enough to resolve the   (finite) system of equations:
   $\{c_{i }+h_{i }(c_*)=c^g_{i }\}_{i }$.
 And now apply the $IFT_\one$, see \S\ref{Sec.Preliminaries.Basic.Tools}.i, to get the solutions $c_*(x)\in R_X$.

 Finally we take the coordinate change  $\Phi_X: x\to x+\sum  c_{i }(x)\xi_i(x)$.
 It is invertible (as  $\ca\sseteq\cm^2$), and
   defines the needed automorphism $\Phi_X\in Aut_X$.

\item Take a perturbation $g\in \ca^2\cdot T_\cR f.$ By part 1
  we have a  map  $\Phi_X:x\to x+ \xi(x)$, with $\xi\in \ca^2\cdot T_\cR$, satisfying $\Phi^{-1}_X(f+g)=f$. As $J\neq 0$, this $\Phi_X$
 is not necessarily a ``coordinate change" on $X$, i.e. not an automorphism of $R_X$. By the $jet_0$-assumption  we can extend $\Phi_X$ to an automorphism $\Psi_X\in Aut_X$, satisfying
   $\Psi_X(x)- x -\xi(x)\in \ca^3$. Then $\Psi^{-1}_X(f+g)-f\in \ca^3\cdot I^{ord(f)-1}\cdot\RpX\sseteq \ca^2\cdot\cm\cdot T_\cR f.$
 (Note that we cannot assume   $\Psi^{-1}_X(f+g)-f\in \ca^3\cdot T_\cR f,$ because $\Psi_X-Id-\xi$ is not a derivation of $R_X$.)

Iterate  this argument  to get: $\Psi_{X,d}^{-1}(f+g)-f\in \ca^3\cdot\cm^d\cdot \RpX$ for each $d\ge1$.
 Hence $f+g\in \overline{\cR f}$,  the orbit-closure in the $\cm^\bullet$-filtration topology.
 The statement follows now by lemma \ref{Thm.Filtration.Closures.TR.TK}.
  \epr\eee

\bex Let $p\!\ge\!2$ and $R_X$ be one of $\k[[x]], \k\{x\}, \k\bl x\br$, \!see \S\ref{Sec.Preliminaries.Rings.Germs}.ii.
  \!Take $I\!=\cm$ and $(\k^n,o)\!\stackrel{f}{\to} \!(\k^p,o).$
\bee[\bf i.]
\item Suppose $f$ is a submersion, i.e. $T_\cR f=\RpX$.  (Thus $p\le n$.)
Then for $\ca=\cm^2$ one gets $\cR f\supseteq\{f\}+\cm^2\cdot \RpX$.
 In particular one gets the normal form of submersion, $f\stackrel{\cR}{\sim}(x_1,\dots,x_p)$.
\item More generally, the theorem gives:
    If    $\ca^2\! \sseteq\cm\cdot  \ca\cdot \ca_\cR$
 then    $\cR f \supseteq \{f\} \!+\ca\cdot T_{\cR} f  $.
\\(Here $\ca_\cR$ is the annihilator of the critical locus,  \S\ref{Sec.3.Annihilators}.)
E.g. one has $\cR f \supseteq \{f\} \!+ \cm\cdot \ca^2_\cR\cdot   \RpX   $.
\item In many cases one has $\ca_\cR\cdot \RpX\sseteq \cm\cdot T_\cR f$, see example \ref{Ex.Annihilators.R.K.simple}.iv.
 Then theorem \ref{Thm.IFT.R.case} gives:  $\cR f\supseteq \{f\}+\ca^2_\cR\cdot \RpX$.
  Geometrically: $f$ is determined up to $\cR$-equivalence by its 2-jet on the critical locus $Crit(f)=V(\ca_\cR)$, taken with its annihilator structure.
\eee
\eex

\subsection{More examples for the case of one power series, i.e. $\pmb{p=1}$}\label{Sec.R.Examples}
 Let $R_X$ be one of $\k[[x]]$, $\k\{x\}$, $\k\bl x\br$. Take $f\in \cm ^2 $,      $I=\cm $ and $\ca\sseteq\cm^2$ .  Theorem \ref{Thm.IFT.R.case}    reads:
\beq
\text{If   \quad  $\ca^2\cdot \cm ^{ord(f)-2} \sseteq\cm\cdot  \ca\cdot Jac(f)$
 \quad  then  \quad $\{f\}+\ca\cdot Jac(f) \sseteq \cR f $.}
\eeq
 Below we use a simpler  though more restrictive condition:  $\ca\cdot \cm ^{ord(f)-2} \sseteq\cm\cdot   Jac(f)$.
 Then we get:
\bcor
\bee[\bf 1.]
\item $\cR f\supseteq \{f\}+Jac(f)\cdot \Big((\cm\cdot Jac(f)):\cm ^{ord(f)-2}\Big)$.
\item Take the smallest $d\in \N$ satisfying $\cm ^{2d+ord(f)-2}\sseteq \cm ^{d+1}\cdot Jac(f)$. Then $\cR f\supseteq \{f\}+\cm ^d\cdot Jac(f)$.
\eee
\ecor
\bex\label{Ex.R.IFT.p=1}
\bee[\bf i.]
\item
 Obviously $(\cm\cdot Jac(f)):\cm ^{ord(f)-2}\supseteq \cm\cdot Jac(f)$. Therefore: $\{f\}+\cm\cdot Jac(f)^2\sseteq \cR f $.
 This is well known for $R_X=\C\{x\}$, e.g. \cite[lemma 2.2, pg.91]{Ruiz}, or for $char(\k)=0$, \cite{Kucharz}.

 In particular, if $\sqrt{Jac(f)}=\cm$ then $f$ is finitely-determined.

More generally, in this way one can extend many other results of   \cite{Kucharz}. We omit the details.
 \item (Morse lemma)
{\em Let $\k$ be a field, $f\in \cm^2$, and $rank[f''|_o]=r$. Then  $f\stackrel{\cR}{\sim}Q_2(x)+\tf(w)$, where $x=(x_1,\dots,x_r)$, $w=(w_1,\dots,w_{n-r})$,
 $\tf(w)\in (w)^3$ (is independent of $x$)
 and $Q_2(x)$ is   a homogeneous quadratic polynomial.}

  Moreover, for $\k=\bar\k$ one can diagonalize the quadratic form $Q_2$ to get: $f \stackrel{\cR}{\sim}\sum x^2_i+\tf(w)$.
\bpr
 Apply a $GL(n,\k)$ transformation to get $f\!-\!Q_2(x)\!\in  \! (x,w)^3$. Then $Jac(f)\!+\!(x,w)^2\!\supset\! (x_1,\dots,x_r)$.
 And thus we have $f -Q_2(x)\in (x,w)\cdot Jac(f)^2+ (w)^3$. The statement follows now by part i.
\epr

\item Assuming $ord(f)\ge3$ one has $(\cm\cdot Jac(f)):\cm ^{ord(f)-2}\supseteq Jac(f)$.
 Then we get: $\{f\}+ Jac(f)^2\sseteq \cR f $. For $R_X=\C\{x\}$ this is Lemma 3.20 in \cite{Ebeling}. Its proof is based on vector field integration,
  thus cannot be used when $\k$ is an arbitrary field.

\item If $ord(f)=2$ then by Morse lemma we split the variables,  $f=Q_2(x)+\tf(w)$  with $\tf(w)\in(w)^3 $. Then one has:
 $\cR f\supseteq \{f\}+(x)^2\cdot (x,w)+(x)\cdot Jac(\tf(w))+Jac(\tf(w))^2$.
\item To show how this corollary strengthens the classical bounds, let $\k$  be a field of characteristic zero, and take
$f$ with an ordinary multiple point at $o\in (\k^n,o)$.
 This means: the partials  $\di_1 f,\dots,\di_n f$ form a regular sequence, and their leading terms  (of order $ord(f)-1$) form a regular sequence as well.
 Let $q$ be the smallest integer satisfying
 $\cm ^q\sseteq \cm ^2\cdot Jac(f)$. Then  $q= n\cdot (ord(f)-2)+1$, see exercise 2.2.7 in \cite{Gr.Lo.Sh}. The known determinacy bound is:
    $\{f\}+\cm^q\sset \cR f$, see \S5 of \cite{BGK.20}.

{\em We claim:} $\cR f\supseteq \{f\}+\cm ^{\lceil\frac{q-ord(f)}{2}\rceil}\cdot Jac(f)$. (This bound is better, even asymptotically.)
\bpr
As $f$ has an ordinary multiple point we have: $\cm ^q\!=\!Jac(f)\!\cdot\! \cm ^{q-ord(f)+1}$. Assume $ord(f)\!\ge\! 2+\frac{3}{n-1}$ and denote
  $d\!:=\lceil\frac{q-ord(f)}{2}\rceil\!+\!1$. Take $\ca\!=\!\cm ^d\!\sset\! I\!=\!\cm $. To use part 2 of the corollary we verify:
  \beq
  \ca^2\cdot \cm^{ord(f)-2}\sseteq \cm\cdot\ca\cdot  Jac(f), \quad \quad \text{   i.e. } \quad \quad
   \cm ^{2d+ord(f)-2}\sseteq   \cm ^{d+1}\cdot Jac(f).
   \eeq
     Indeed,     by our assumption:
   \beq
   2d+ord(f)-2\ge q \quad \text{ and } \quad (d+1)+(ord(f)-1)=\lceil\frac{q-ord(f)}{2}\rceil+1\le q.
 \hspace{2cm}
 \epr
   \eeq

\item
Take the $E_6$ singularity, $f(x,y)=x^3+y^4$. Assume $\k$ is a field, $char(\k)\neq 2,3$. Denote $\cm=(x,y)$,   $\ca=\cm^2\sset I=\cm$.
 We get:  $\ca^2\cdot \cm^{ord(f)-1}\sseteq \cm\cdot \ca\cdot Jac(f) $. Therefore $\cR f\supseteq \{f\}+\cm^2\cdot (x^2,y^3) $.
  The ideal $\cm^2\cdot (x^2,y^3) $ contains all the monomials lying above the Newton diagram of $x^3+y^4$. Compare this result to the
   standard determinacy statement: $\cR f\supseteq \{f\}+\cm^5$.
\eee
\eex

\subsection{The filtration criterion for $\cR$-orbits}\label{Sec.R.Filtration}
  Take  the  filtration \!$M_\bullet\!\!:=\!I^\bullet\!\cdot\!  \RpX.$  Fix some integers $1\!\le\! j\!<\!d.$
\bthe\label{Thm.Orbits.R.Filtration}
 Suppose      the conditions $jet_N (Exp)$, $jet_N(Ln)$ of \S\ref{Sec.Preliminaries.Exp.Ln}
  hold for $N=\lceil \frac{d-ord(f)}{j}\rceil$. Then:

\hspace{2.5cm}
$ {\cR^{(j)} f}\supseteq \{f\}+I^d\cdot  \RpX$      \quad\quad\quad \text{ if and only if }\quad\quad\quad
   ${T_{\cR^{(j)}} f}\supseteq I^d\cdot  \RpX$.
\ethe
\bpr
{\bf The part $\boldsymbol{``\Rrightarrow"}$.} First we prove: $  T_{\cR^{(j)}} f+ M_{d+1}\supseteq  M_{d}$.
  Let $w\in M_{d}$, then $gf\equiv f+w\ mod\ M_{d+1}$ for some $g\in \cR^{(j)}$. Use the assumption $jet_N (Ln)$ to get the derivation
 $-\sum^{N}_{i=1}\frac{(\one-g)^i}{i}+\xi_g\in T_{\cR^{(j)}}$, for some
   $\xi_g\in End^{(Nj+1)}_\k(R)$. And then
\beq
T_{\cR^{(j)}}f\ni \Big(-\sum^{N}_{i=1}\frac{(\one-g)^i}{i}+\xi_g\Big)f\equiv w\ mod(M_{d+1}+M_{Nj+1+ord(f)})
\equiv w\ mod\ M_{d+1}.
\eeq
Here $Nj+ ord(f)\ge d$ because of the assumption $N=\lceil \frac{d-ord(f)}{j}\rceil$.

Therefore we have $  T_{\cR^{(j)}} f+ M_{d+1}\supseteq  M_{d}$. As $M_{d+1}=I\cdot M_d$, we get (by Nakayama) $  T_{\cR^{(j)}} f\supseteq  M_{d}$.

{\bf The part $\boldsymbol{``\Lleftarrow"}$.}
First we prove:  $  \cR^{(j+k)}  f+ M_{d+k+1}\supseteq \{f\}+M_{d+k}$ for every $k\ge 0$.  Indeed, $T_{\cR^{(j)}}  f\supseteq M_{d}$ implies
 $T_{\cR^{(j+k)}}  f\supseteq M_{d+k}$ for every $k\ge 0$. (Using $M_{\bullet}=I^\bullet\cdot \RpX$.)
 Accordingly we present an  element $w\in M_{d+k}$ as $w=\xi (f)$, for some $\xi\in  T_{\cR^{(j+k)}} $.
 The condition $jet_N (Exp)$ ensures an element $\Phi:=\sum^{N}_{i=0}\frac{\xi^i}{i!}+\phi_\xi\in \cR^{(j+k)}$, where
  $\phi_\xi\in End^{(N(j+k)+1)}_\k(R_X)$. And then
\[
\Phi(f)-f-w\in  M_{d+k+1}+M_{ (j+k)N+ord(f)+1}.
\]
 Finally use $(j+k)N+ord(f)\ge d+k$ to get
   $\Phi(f)-f-w\in  M_{d+k+1}$.

\

Iterating $  \cR^{(j+k)}  f+ M_{d+k+1}\supseteq \{f\}+M_{d+k}$ for $k\ge 0$ we get: $  \overline{\cR^{(j)}  f} \supseteq \{f\}+M_{d}.$
 Now apply part 2 of Lemma \ref{Thm.Filtration.Closures.TR.TK}.
\epr

\bcor\label{Thm.Orbits.Rj.vs.TRj.for.nice.rings}\!($R_X$ is one of  $ {\k[[x]]},{\k\{x\}},{\k\bl x\br} $, with $\k$ a \!local ring,  \!see  \!\S\ref{Sec.Preliminaries.Rings.Germs}.ii.)  Take $I\!\sseteq\!(x)$  \!and     \!$f\!\in \!(x)\!\cdot \!\RpX.$
\bee[\bf \!1.]
\item
 Suppose  $2,3,\dots,N\in \k^\times$ for $N\!=\!\lceil\frac{d-ord(f)}{j}\rceil.$
 Then:
$\cR^{(j)} f \supseteq \{f\}+I^{d}\cdot  \RpX$   
 if and only if 
$ T_{\cR^{(j)}} f \supseteq I^{d}\cdot  \RpX$.
\item
Suppose   $ 2\in \k^\times$. If $ T_{\cR^{(j)}} f \supseteq I^{d}\cdot  \RpX$
 then $\cR^{(d-j-ord(f))} f \supseteq \{f\}+I^{2d-2j-ord(f)}\cdot  \RpX$.
\eee
\ecor
\bpr
{\bf 1.}
Combine theorem \ref{Thm.Orbits.R.Filtration} with lemma \ref{Thm.Filtration.Closures.TR.TK}.
 The $jet_N$-condition holds by example \ref{Ex.jet_N.conditions.hold}.

{\bf 2.} We get $ T_{\cR^{(d-j-ord(f))}} f \supseteq I^{2d-2j-ord(f)}\cdot  \RpX$.
 Now apply part one with  $N:=\lceil\frac{2d-2j-2ord(f)}{d-j-ord(f)}\rceil=2$.
 \epr
\bex
\bee[\bf i.]
\item
Take $p=1$, $I=\cm$, $j=1$ and suppose $\k$ is a field, with $char(\k)=0$ or $char(\k)> d-ord(f)$. Then:
 $\cR^{(1)} f \supseteq \{f\}+\cm^d$        if and only if $\cm^2\cdot Jac(f)\supseteq \cm^d$.

In particular, the $\cR$-order of determinacy of $f$ is $\le \mu(f)+1$.
In $char(\k)=0$ this is well known, e.g. \cite[Corollary 2.24]{Gr.Lo.Sh}.
  For $char(\k)> d-ord(f)$ this improves the known bound, \cite{Boubakri.Gre.Mark}, \cite{Greuel-Pham.2017}:
  the $\cR$-order of determinacy of $f$ is $\le 2\mu(f)-ord(f)+2$.

\item In part two of corollary \ref{Thm.Orbits.Rj.vs.TRj.for.nice.rings} the conclusion is weaker, but also the assumption is weaker.
 For $j=1$ this part two is known, e.g. see corollary 5.2 of \cite{BGK.20}.
\eee
\eex

\section{Criteria for   contact orbits, ``$\cK f$  vs $T_\cK f$"}\label{Sec.Orbits.K}
Let $R_X$ be one of the rings $\quots{\k[[x]]}{J},\! \quots{\k\{x\}}{J}, \!\quots{\k\bl x\br}{J},$ with $\k$ a local ring, see \S\ref{Sec.Preliminaries.Rings.Germs}.ii.
\subsection{The  $\cK$-implicit function theorem}\label{Sec.K.IFT}
     Take  two ideals $\ca\sseteq I\sseteq  R_X$, with $\ca\sseteq\cm^2$, and a map-germ $f\in I\cdot \RpX$.
   (The usual choice is   $I=R_X,$ $I=\cm$ or $I=\sqrt{(f)}$.) For the notation $I^{ord(f)-2}$ see \S\ref{Sec.Preliminaries.Rings.Germs}.v.
\bthe\label{Thm.IFT.K.case}
Suppose
 $\ca^2\cdot I^{ord(f)-2}\cdot \RpX \sseteq\cm\cdot  \ca\cdot T_\cR f+\cm\cdot(f)\cdot \RpX.$
\bee
\item (The case $J=0$.) Then
      $\{f\}+  (\ca ^2\cdot I^{ord(f)-2} +\cm\cdot(f))\cdot \RpX\sseteq \cK f $.
\item (The case $J\neq0,$ assume the $jet_0$-condition of \S\ref{Sec.Preliminaries.Exp.Ln}.)
       Then
      $\{f\}+  (\ca^3\cdot I^{ord(f)-2} +\cm\cdot(f))\cdot \RpX\sseteq \cK f $.
\eee
\ethe
A remark: the assumption  $\ca^2\cdot I^{ord(f)-2}\cdot \RpX \sseteq\cm\cdot  \ca\cdot T_\cR f+\cm\cdot(f)\cdot \RpX$ is preserved by $\cK$-equivalence.
\bpr
\bee[\bf 1.]
\item Take a perturbation $g\in(\ca^2\cdot I^{ord(f)-2}  +\cm\cdot(f))\cdot \RpX$. We want to resolve the condition $\Phi_X(f)=(\one+U)\cdot (f+g)$, where $\Phi_X\in Aut_X$
 and $U\in Mat_{p\times p}(\cm)$.
   As in the proof of Theorem \ref{Thm.IFT.R.case} we fix some (finite) sets of generators
       $\{\xi_i(f)\}$ of $\ca\cdot T_\cR f\sseteq \RpX$.

  Present the coordinate change in the form $\Phi_X:x\to x+\sum   c_{i }\xi_i(x)$, here $\{c_{i }\}$ are unknowns.
\bei
\item
    As in the $\cR$-case we get:
 $\Phi_X(f)-f-\sum  c_{i }\xi_i(f)\in (\{c_*\})^2\cdot\ca^2\cdot I^{ord(f)-2}\cdot \RpX$.

 Therefore, as in the $\cR$-case, we can present
\beq\label{Eq.Proof.of.K.current.situation}
\Phi_X(f)=f+\sum   \big( c_{i }+H^{(1)}_{i }(\{c_*\})\big)\xi_i(f)+H^{(2)}(\{c_*\})\cdot f.
\eeq
Here  $H^{(1)}_{i } (z)\in \cm\cdot R_X[[\cm\cdot z]][z]$ and $H^{(1)}_i(\cm)\sseteq \cm\sset R_X$.
 Similarly, $H^{(2)} (z)\in Mat_{p\times p}(\cm\cdot   R_X[[\cm\cdot z]][z])$ and $H^{(2)} (\cm)\sset Mat_{p\times p}( \cm)$.

\item
Expand $g=g_\ca+U^g\cdot f$, here $g_\ca\in \ca^2\cdot I^{ord(f)-2}\cdot\RpX$ and $U^g\in Mat_{p\times p}(\cm)$.
  Thus $(\one+U)(f+g)=f+(\one+U)g_\ca+(U+U^g+U\cdot U^g)f.$
    We can present $(\one+U)g_\ca=\sum_j l_{j}(U)\cdot \xi_j(f)+\tU^g\cdot f$,
    where
   $l_{j}(U)$ and $\tU^g$ are linear functions in the entries of $U$.
\eei

   Altogether the condition $\Phi_X(f)=(\one+U)\cdot (f+g)$ becomes:
\beq
f+\sum   \big( c_{i }+H^{(1)}_{i }(\{c_*\})\big)\xi_i(f)+H^{(2)}(\{c_*\})\cdot f=f+\sum_{i}   l_{i }(U)\cdot \xi_i(f)+(U+U^g+U\cdot U^g+\tU^g)\cdot f.
\eeq
Therefore we have transformed the condition $\Phi_X(f)=(\one+U)(f+g)$ into  the system of equations:
\beq
c_{i }+H^{(1)}_{i }(\{c_*\})=    l_{i }(U),\quad
H^{(2)}(\{c_*\})= U+U^g+U\cdot U^g+\tU^g, \quad \text{with the unknowns $\{c_*\},U$} .
\eeq
Apply the $IFT_\one$, \S\ref{Sec.Preliminaries.Basic.Tools}.i, to determine $c_*(x)\in \cm$ and $U\in Mat_{p\times p}(\cm)$.
 We get the transformation $(\Phi_X,\one+U)\in \cR\ltimes GL(p,R_X)$ that satisfies  the condition $\Phi_X(f)=(\one+U)\cdot (f+g)$.

\item
It is enough to prove: $\cK f\supseteq \{f\}+\ca^3\cdot I^{ord(f)-2}\cdot \RpX$ (using the $GL(p,R_X)$-transformation).
Take a perturbation $g\in \ca^3\cdot I^{ord(f)-2} \cdot \RpX.$  The previous argument gives
 a map  $\Phi_X:x\to x+\xi(x)$, with $\xi(\cm)\sseteq \cm\cdot \ca^2$, and a matrix $U$, satisfying $\Phi_X(f)=(\one+U)(f+g)$.
  As in the proof of theorem \ref{Thm.IFT.R.case} this $\Phi_X$
 is not necessarily ``a coordinate change" on $X$. By the $jet_0$-assumption we can extend $\Phi_X$ to a coordinate change
  $\Psi_X\in Aut_X$, satisfying:
   $\Psi_X( x)-x-\sum   c_{i }\xi_i(x)\in  \ca^3\cdot\cm$. Then
    $\Psi_X^{-1}\big((\one+U)(f+g)\big)-f\in  \cm \cdot \ca^3\cdot  I^{ord(f)-2}\RpX$.
 (As in the $\cR$-case we remark: one cannot assume  $\Psi_X^{-1}\big((\one+U)(f+g)\big)\in\{f\}+  \ca^3\cdot\cm\cdot   T_\cR f$.)

 Iterate this argument  to get $ f+g \stackrel{\cK}{\sim}   \{f\}+ \cm^d\cdot\ca^3\cdot I^{ord(f)-2}\cdot  \RpX$ for each $d\ge1$.
 Thus $f+g\in \overline{\cK f}$,  the orbit-closure in the $\cm^\bullet$-filtration topology.
   The statement follows now by part 2 of lemma \ref{Thm.Filtration.Closures.TR.TK}.
\epr
\eee

\bex\label{Ex.K.IFT}
Let $R_X$ be one of the rings $\k[[x]]$, $\k\{x\}$, $\k\bl x\br$, with $\k$ a local ring, see \S\ref{Sec.Preliminaries.Rings.Germs}.ii.
 Take $I=\cm$ and assume $(f)\sseteq\cm^2$. We get the following bounds on the orbit $\cK f.$
\bee[\bf i.]
\item  $\cK f\supseteq \{f\}+  (\cm^2\cdot \ca^2_\cR  +\cm\cdot(f))\cdot \RpX $.
 For $p=1$ this reads: $\cK f\supseteq \cm^2\cdot Jac(f)^2+f\cdot \cm$.
\bei
\item
Suppose $f$  has an isolated singularity, i.e.  $Crit(f)\cap V(f)=\{o\}\in \k^n$, i.e.
 $\sqrt{\ca_\cR+(f)}=\cm.$  Then $f$ is finitely-$\cK$-determined.
\item The inclusion  $\cK f\supseteq \{f\}+  (\cm^2\cdot \ca^2_\cR  +\cm\cdot(f))\cdot \RpX $ strengthens Theorem 2.1 of \cite{Cutkosky-Srinivasan.1997},
 which reads:
 {\em If $R_X=\k[[x]]$ and $\ca=(f)+Fitt_0(T^1_\cR f)$ then $\cK f\supseteq \{f\}+\ca^3\cdot\RpX$}.

Note that $\ca_\cR\supseteq Fitt_0(T^1_\cR f)$, and for $p\ge2$ one often has $\cm\cdot \ca_\cR\supseteq Fitt_0(T^1_\cR f)$, see example
 \ref{Ex.Annihilators.R.K.simple}.
\eei
\item If $(f)\sseteq\cm^3$ then $\cK f\supseteq \{f\}+  (\ca^2_\cR  +\cm\cdot(f))\cdot \RpX $.
For $p=1$ this reads: $\cK f\supseteq  \{f\}+ Jac(f)^2+f\cdot \cm$.

\item For $p=1$ and $f\in \cm^2$ use Morse lemma, example \ref{Ex.R.IFT.p=1}.ii, to present $f=Q_2(x)+\tf(w)$, with $\tf(w)\in (w)^3$.
 Then we get:   $\cK f\supseteq \{f\}+  \cm\cdot (x)\cdot((x)+Jac_w(\tf(w)))+(Jac_w(\tf(w)))^2  +\cm\cdot(f) $.

 \item For $p\ge2$ the assumption on the modules $\ca^2\cdot I^{ord(f)-2}\cdot \RpX\sseteq\cm\cdot  \ca\cdot T_{\cR} f+\cm\cdot(f)\cdot \RpX$ is
 weaker than the assumption on the  ideals $\ca^2\cdot I^{ord(f)-2} \sseteq\cm\cdot  \ca\cdot \ca_\cR+\cm\cdot(f).$
   But the later assumption is often simpler to verify.

\item Below we assume  $p=1$.
\bei
\item Suppose $\ca^2\cdot \cm ^{ord(f)-2}\sseteq \cm\cdot \ca\cdot Jac(f)+f\cdot \cm$. Then $\cK f\supseteq \{f\}+\ca\cdot Jac(f)+f\cdot \cm$

\item Suppose $\cm^e\cdot Jac(f)+(f)\cm\supseteq \cm^d$ for $e=\lceil\frac{d-ord(f)}{2}\rceil+1$.

Take $\ca=\cm^e$ to get:
 $\cK f\supset\{f\}+(\cm^e\cdot  Jac(f) +\cm(f))\cdot \RpX$.
\item
  Take the smallest $d\in \N$ satisfying $\cm ^{2d+ord(f)-2}\sseteq \cm ^{d+1}\cdot Jac(f)+\cm \cdot(f)$.

  Then $\cK f\supseteq \{f\}+\cm ^d\cdot Jac(f)+\cm \cdot f$.
\eei
\eee
\eex

\subsection{An application: the $\cK$-orbit of a reduced complete intersection curve germ via the  semigroup of values}\label{Sec.K.application.curve.germ}
Let  $R_X$ be one of $\k[[x]]$, $\k\{x\}$, $\k\bl x\br$, where $\k=\bar\k$ is a  field, see \S\ref{Sec.Preliminaries.Rings.Germs}.ii.
 Suppose $(C,o):=V(f)\sset (\k^n,o)$ is a reduced complete intersection curve germ. (Here $f=(f_1,\dots,f_{n-1})$.)
 Take its branch decomposition, $(C,o)=\cup^r_{i=1}(C_i,o) $. We have the normalization morphism
 $\cO_{(C,o)}\into\prod \cO_{(\tC_i,o)}$. Here each ring $\cO_{(\tC_i,o)}$ is a DVR.
  It has the valuation map  $val_i: \cO_{(\tC_i,o)}\stackrel{}{\to}\bN:=\N\cup\{\infty\} $.
 The ``full" valuation map is the product of morphisms, $val:R_X\to \cO_{(C,o)}\into \prod \cO_{\tC_i,o}\to \oplus_i \bN$.
 Accordingly we have two sub-semigroups: $val(\cO_{(C,o)})\sseteq \oplus_i \bN $ and
  $val(\ca^2_\cR)\sseteq \oplus_i \bN $. Note that $val(\ca^2_\cR)=val(\ca^2_\cR+(f))$.

Below we consider the ``$\N$-submodules" of $\N^r$, i.e. subsets $S\sset\N^r$ satisfying: $S+\N^r=S$.

\bprop
Suppose a submodule $\cS\sset \bN^r$ satisfies: $\cS\cap val (R_X)=\cS\cap val(\cm\cdot\ca^2_\cR)$.
 Then $\cK f\supseteq \{f\}+\cm\cdot val^{-1}(\cS)\cdot \RpX$.
\eprop
 If $(C,o)$ is not a plane curve germ, i.e. $p\ge2$, then
 in many cases one has  $\ca_\cR\cdot \RpX\sseteq \cm\cdot T_\cR f,$ rather than just $\ca_\cR\cdot \RpX\sseteq   T_\cR f,$ see example \ref{Ex.Annihilators.R.K.simple}.iv.
  Then the bound is strengthened to:   $\cK f\supseteq \{f\}+  val^{-1}(\cS)\cdot \RpX$.

\bpr
Let $g\in val^{-1}(\cS )\sset R_X$. Then $val(g)=val(q)$ for an element $q\in  \cm\cdot\ca^2_\cR$. Thus  $val(g-c\cdot q)>val(g)$
 for some constant $c\in \k$.
 Namely, $val_i(g-c\cdot q)\ge val_i(g)$ for each $i=1,\dots,r$, and this inequality is strict for at least  one $i$.

 Apply this procedure for each coordinate of the vector $val(g)\in\bN^r$ to get an element $q\in \ca^2_\cR$ that satisfies:
  $val(g-c\cdot q)\ge val(g)+(1,\dots,1)$.
 Iterating this we get $val^{-1}(\cS )\sset    \ca^2_\cR+(f) +val^{-1}(\cS+ d(1,\dots,1) )$ for each $d\ge1$.

As $(C,o)$ is a reduced germ, the ideal $\ca^2_\cR+(f)$ is $\cm$-primary. Therefore we get: $val^{-1}(\cS )\sset    \cm\cdot\ca^2_\cR+(f)$.
 And hence $\{f\}+\cm\cdot val^{-1}(\cS )\cdot \RpX\sset \{f\}+\cm\cdot ( \cm\cdot\ca^2_\cR+(f))\cdot \RpX$.
 Finally, apply Example \ref{Ex.K.IFT}.i.
\epr

\subsection{The filtration criterion for   $\cK$-orbits}\label{Sec.K.Filtration}
  Take the filtration \!$M_\bullet\!\!:=\!I^\bullet\!\cdot \! \RpX.$  Fix some integers $1\!\le\! j\!<\!d.$
\bthe\label{Thm.Orbits.K.Filtration}
 Suppose      the conditions $jet_N (Exp)$, $jet_N(Ln)$ of \S\ref{Sec.Preliminaries.Exp.Ln}
  hold for $N=\lceil \frac{d-ord(f)}{j}\rceil$. Then:

  \quad\quad\quad
$ {\cK^{(j)} f}\supseteq \{f\}+I^d\cdot  \RpX$      \quad\quad\quad \text{ if and only if }\quad\quad\quad
   ${T_{\cK^{(j)}} f}\supseteq I^d\cdot  \RpX$.
\ethe
\bpr
\bee[\bf Step 1.]
\item In the proof of $\cR$-case, Theorem \ref{Thm.Orbits.R.Filtration},
  we have used the transitions $T_{\cR^{(j)}}\rightleftarrows \cR^{(j)}$, ensured by the $jet_N$ assumptions.
 Now we need the transitions $T_{\cK^{(j)}}\rightleftarrows \cK^{(j)}$.
 These are ensured by lemma 3.22 of \cite{BGK.20}:
\li{\em For any  element $(u,\xi)\in T_{\cK^{(j)}}$ there exists an element $(U,\Phi)\in \cK^{(j)}$ satisfying:}
\[ord_I[(U,\Phi)f-f-(u,\xi)f]>ord_I [(u,\xi)f].\]
 $\bullet$ {\em For any  element $(U,\Phi)\in \cK^{(j)}$  there exists an element $(u,\xi)\in T_{\cK^{(j)}}$ satisfying:}
\[ord_I[(U,\Phi)f-f-(u,\xi)f]>ord_I [(U,\Phi)f-f].\]

The statement of that lemma in \cite{BGK.20} assumes the conditions $jet_N(Exp)$,  $jet_N(Ln)$ for each $N$.
 In particular this implies:  $\k\supseteq \Q$.
 However, the proof in \cite{BGK.20} uses these conditions only up to $N=\lceil\frac{d-ord(f)}{j}\rceil$,   as we explain now.

 The only part of the proof involving the conditions $jet_N(Exp)$,  $jet_N(Ln)$ is the Baker-Campbell-Hausdorff formula,
  \S\ref{Sec.Preliminaries.Basic.Tools}.iv.
 For the transitions $T_{\cK^{(j)}}\rightleftarrows \cK^{(j)}$ we need the expansion only up to order $l$ that satisfies: $l\cdot j+ord(f)\ge d$.
  Thus it suffices to expand up to order $l=\lceil\frac{d-ord(f)}{j}\rceil$, which
  coincides with $N$. Our $jet_N(\dots)$ assumptions imply $2,\dots,N\in \k^\times$. Therefore the factor $(l!)^4$ is invertible.
  Therefore the Baker-Campbell-Hausdorff formula, \eqref{Eq.Baker.Campbell.Hausdorff},
   holds in $R_X$ up to order $l=N$. Thus  lemma 3.22 of \cite{BGK.20} is applicable to our case.

\item
 We restate lemma 3.22 of \cite{BGK.20} in the form:
 \beq\label{Eq.inside.proof}
 w_d\in T_{\cK^{(j)}}f+M_{d+1} \quad\quad \text{ iff } \quad \quad w_d\in (\cK^{(j)}f-f)+M_{d+1}.
 \eeq
Using this equivalence the proof goes as follows.
\bei
\item If \!$\overline{\cK^{(j)}}f\!\supseteq\! \{f\}\!+\!M_d$ then \!$T_{\cK^{(j)}}f\!+\!M_{d+1}\!\supseteq\! M_{d}$. As \!$M_{d+1}\!=\!I\!\cdot\! M_d$ we get (by Nakayama):
  \!$T_{\cK^{(j)}}f \!\supseteq \!M_{d}.$
\item If $T_{\cK^{(j)}}f \supseteq M_{d}$ then $T_{\cK^{(j+l)}}f \supseteq M_{d+l}$ for every $l\ge0$.
 For each $l$ we use \eqref{Eq.inside.proof} with $N=\frac{d+l-ord(f)}{j+l}\le \frac{d-ord(f)}{j}$.
 And then
   $ \cK^{(j+l)} f+M_{d+l+1}\supseteq \{f\}+M_{d+l}$ for every $l\ge0$. Thus $\overline{\cK^{(j)}}f\supseteq \{f\}+M_d.$
    Finally, apply part 2 of Lemma \ref{Thm.Filtration.Closures.TR.TK}.
  \epr
  \eei
\eee

\bcor\label{Thm.Orbits.Kj.vs.TKj.for.nice.rings}
 ($R_X$ \!is \!one \!of  $ \k[[x]]  ,  {\k\{x\}} ,   { \k\bl x\br}{}$, \!with \!$\k$ a \!local \!ring, \!see\! \S\ref{Sec.Preliminaries.Rings.Germs}.ii.)
  Take $I\!\!\sseteq\!\!(x)$  and     $f\!\in\! (x)\!  \RpX.$
  \bee[\bf 1.\!\!]
  \item
 Suppose $2,3,\!\dots,\!N\!\in\! \k^\times$  for $N\!=\!\lceil\frac{d-ord(f)}{j}\rceil.$
 Then:
$\cK^{(j)} f\! \supseteq \!\{f\}\!+ I^{d}\!\cdot \! \RpX$        if and only if
${T_{\cK^{(j)}} f}\!\supseteq\! I^{d}\!\cdot   \RpX.$

\item
Suppose   $2\in \k^\times.$ If $ T_{\cK^{(j)}} f \supseteq I^{d}\cdot  \RpX$
 then $\cK^{(d-j-ord(f))} f \supseteq \{f\}+I^{2d-2j-ord(f)}\cdot  \RpX$.
\eee
\ecor
 The proofs of both parts are the same as in Corollary \ref{Thm.Orbits.Rj.vs.TRj.for.nice.rings}.

\bex
\bee[\bf i.]
\item
Take $p\!=\!1$, $I\!=\!\cm$, $j\!=\!1$ and suppose $\k$ is a field, $char(\k)\!=\!0$ or $char(\k)\!>\! d-ord(f).$
\beq
\text{ Then: \quad $\cK^{(1)} f \supseteq \{f\}+\cm^d$     \quad    if and only if \quad $\cm^2\cdot Jac(f)+\cm\cdot (f)\supseteq \cm^d$.}
\eeq
In particular, the $\cK$-order of determinacy of $f$ is $\le \tau(f)+1.$
For $char(\k)=0$ this is well known, e.g. \cite[Corollary 2.24]{Gr.Lo.Sh}.
 For  $char(\k)> d-ord(f)$ this strengthens the bound $(2\tau(f)-ord(f)+2)$ of \cite{Boubakri.Gre.Mark2011}, \cite{Greuel-Pham.2017}.
\item For $j=1$ (and $\k$ a field) part two of Corollary \ref{Thm.Orbits.Kj.vs.TKj.for.nice.rings} is known, e.g. see Corollary 5.10 of \cite{BGK.20}.
\eee
\eex

\section{Properties of the (extended) tangent space $T_\cA f$}\label{Sec.TA}
The tangent space $T_\cA f$ is not an $R_X$-module. And $T_\cA f$ is not finitely generated as a module over $R_Y$.
 As is mentioned in \S\ref{Sec.Intro.Contents}, numerous simple/immediate properties of $T_\cR f, T_\cK f$  become cumbersome (and non-trivial) for $T_\cA f$.
   E.g. even the simplest properties of $T_\cA f$ in the classical $C^\infty$-case rely heavily on the deep Mather-Malgrange preparation theorem.

 Below we study the annihilator $\ca_\cA$,
 the extended module structure of $T_\cA f$,  the Nakayama\&Artin-Rees properties, and the topological closure $\overline{T_\cA f}\sseteq \RpX$.

 Let the pair $(R_X,R_Y)$ be one of
  $(\quots{\k[[x]]}{J},\k[[y]])$, $(\quots{\k\{x\}}{J},\k\{y\})$, $(\quots{\k\bl x\br}{J},\k\bl y\br)$,  \S\ref{Sec.Preliminaries.Rings.Germs}.ii.
Fix the filtration $I^\bullet\cdot \RpX$ of $\RpX,$ we always assume $(f)\sseteq I$.
  Accordingly one has the filtrations $\cA^{(\bullet)}$ and $T_{\cA^{(\bullet)}}$, with $T_{\cA^{(-1)}}:=T_\cA$,
  see \S\ref{Sec.3.Filtrations.on.A.TA}. Fix a map $f:X\to (\k^p,o)$, i.e. $f\in \cm\cdot\RpX$.
   Take the annihilator ideals $\ca_\cA,\ca_{\cA^{(j)}}\sset R_X$ defining the instability locus, see \S\ref{Sec.3.Annihilators}.
     In general $\ca_\cA,\ca_{\cA^{(j)}}\not\sseteq f^\#(R_Y)$.

The map $f$ is called infinitesimally $\cA$-stable if  $T_\cA f= \RpX$, i.e. $\ca_\cA=R_X$.

\subsection{The annihilator $\ca_\cA$}
\subsubsection{} The case $p=1,$ i.e. $f:X\to (\k^1,o),$ is elementary. (Probably well-known.)
 \bel\label{Thm.Annihilator.a_A.p=1.case} ($p=1$)
 Assuming $f\in \cm^2$, one has  $\ca_\cR\sseteq \ca_\cA\sseteq \ca_\cR:\cm$ and $\ca_\cA=\ca_\cR+\ca_\cA\cap f^\#(R_Y)$.
  \eel
In particular, for $R_X$ one of $\k[[x]]$, $\k\{x\}$, $\k\bl x\br$, one has: $Jac(f)\sseteq \ca_\cA\sseteq Jac(f):\cm$.
\bpr (The part $\ca_\cR\sseteq \ca_\cA$ is obvious.)
Take any $g\in \ca_\cA$, then $g\in \ca_\cR+f^\#(y)^d$.
 For any $q\in \cm$ one has  $q\cdot g\in \ca_\cR+f^\#(y)^{d'}$. (Here we chose $d,d'\le \infty$ as the largest possible.)
 \bei
 \item  If $d'\le d$ then $(1-\tq\cdot q)g\in \ca_\cR$ for some $\tq\in R_Y\sset R_X$. And therefore $g\in \ca_\cR$.
\item Suppose  $d'\!>\! d$ for each $q\!\in\! \cm$. Then  $\cm \cdot g\sseteq \ca_\cR+g\cdot(f)\!\sseteq \ca_\cR+g\cdot\cm^2$.
 Nakayama gives: $\cm\cdot g\sseteq \ca_\cR$.
\eei
In both cases we get  $\cm\!\cdot\! \ca_\cA\!\sseteq \!\ca_\cR$. Hence $\ca_\cA\!\sseteq\! \ca_\cR\!:\!\cm$.
The statement $\ca_\cA=\ca_\cR+\ca_\cA\cap f^\#(R_Y)$ follows. \epr

\subsubsection{}\!\!\!For  \!$p\!\ge\!2$ the annihilator  \!$\ca_\cA$ is much more complicated. \!We  list  several basic properties.
 \!Most of them
 \\are immediate/well known.  \!Let the pair \!$(R_X,R_Y)$ be as in \!\S\ref{Sec.3.Maps.of.Germs.Identification}.
 \!Take the filtration \!$I^\bullet\!\cdot\! \RpX,$ \!with $(f)\!\sseteq\! I.$


\begin{Properties}\label{Thm.Annihilators.Basic.Properties}
\bee[\bf i.]
\item $\ca_\cA\cdot \RpX\sseteq \ca_\cK\cdot \RpX+\k^{\oplus p}$.
\item If $(f)\sseteq\ca_\cA$ then $\ca_\cK\sseteq\ca_\cA $.
\\If $f$ is infinitesimally stable, i.e. $\ca_\cA=R_X$, then $T_{\cA^{(0)}}f=T_{\cK^{(0)}}f$.
\\{\rm(In Lemma 7.1 of \cite{Mond-Nuno} this is stated for $I=\cm$,
 but the proof is the same.)}
\\More generally, $T_{\cA^{(0)}}f=T_{\cK^{(0)}}f$ \ iff \ $(f)\cdot \RpX\sseteq T_{\cA^{(0)}}f$ \ iff \ $\ca_{\cA^{(0)}}=\ca_{\cK^{(0)}}$.
\item If $(f)\sseteq \cm^2$ and $\ca_\cA\sseteq \cm$ then $\ca_\cA+(f)\sseteq\ca_\cK$.
 \item If $\ca_\cA\sseteq\cm$ then $\sqrt{\ca_\cA+(f)}=\sqrt{\ca_\cK}\sset R_X$.
\quad \quad Moreover, $\sqrt{\ca_{\cA^{(j)}}+(f)}=\sqrt{\ca_{\cK^{(j)}}}$ for $j\ge0$.
\\{\rm(Geometrically: if $f$ is $\cA$-unstable then the singular locus $Sing(f)$ coincides with the instability locus restricted to the germ $V(f)$.)
\bpr  {\bf The part $\mathbf{\sseteq}$.} For $d\!\gg\!1$ we have:
 $(\ca_\cA+(f))^d\!\cdot\! \RpX\!\sseteq\! T_\cR f+ T_{\cL^{(0)}} f\!\sseteq\! T_\cR f+(f)\!\cdot\! \RpX=T_\cK f.$

The case $j\ge0$ is immediate.

{\bf   The part $\mathbf{\supseteq}$.} Take the filtration $\cm^\bullet\sset R_X$.
Then $\sqrt{\ca_\cA+(f)}\supseteq \sqrt{\ca_{\cA^{(0)}}+(f)}\supseteq \sqrt{\ca_{\cR^{(0)}}+(f)}=\sqrt{\ca_\cR+(f)}=\sqrt{\ca_\cK}$.
 The last transition is by example \ref{Ex.Annihilators.R.K.simple}.

The statement with $j\ge0$ is similar:
  $\sqrt{\ca_{\cA^{(j)}}+(f)}\supseteq \sqrt{\ca_{\cR^{(j)}}+(f)}=\sqrt{\ca_{\cK^{(j)}}}$.
\epr}

\item If $\ca_\cK\sseteq \sqrt{I}$ and $\ca_\cA\sseteq\cm$ then $\sqrt{\ca_{\cA^{(j)}}+(f)}=\sqrt{\ca_\cA +(f)}$ for each $j\ge0$.
\bpr{\rm $\sqrt{\ca_{\cA^{(j)}}+(f)}=\sqrt{\ca_{\cK^{(j)}}}\supseteq\sqrt{I\cdot \ca_{\cK}}=\sqrt{\ca_\cK}=\sqrt{\ca_\cA+(f)}$.}
\epr

\item  $\ca_\cA\cdot\ca_\cK\cdot  \RpX\sseteq (\ca_\cA+(f))\cdot T_\cR f+ T_{\cL^{(0)}} f$.
\\If $(f)\sseteq \cm^2$ and $\ca_\cA\sseteq \cm$ then $\ca_\cA\cdot\ca_\cK\cdot  \RpX\sseteq (\ca_\cA +(f))\cdot T_\cR f+ T_{\cL^{(1)}} f$.

\item    $\ca_{\cA^{(0)}}\cdot\ca_{\cK}\cdot  \RpX\sseteq  \ca_{\cA^{(0)}}\cdot T_{\cR} f+ (f) \cdot T_{\cR^{(0)}} f+ T_{\cL^{(1)}} f$.
\\$\ca_{\cA^{(0)}}\cdot\ca_{\cK^{(0)}}\cdot  \RpX\sseteq ( \ca_{\cA^{(0)}}+ (f))\cdot T_{\cR^{(0)}} f+ T_{\cL^{(1)}} f$.
\\$\ca_\cA\cdot\ca^2_\cK\cdot  \RpX\sseteq (\ca_\cA\cdot \ca_\cK+(f)^2)\cdot T_\cR f+ T_{\cL^{(1)}} f$.
 \item
$\sqrt{\ca_{\cA^{(j)}}}=\sqrt{I\cdot \ca_{\cA}}$ for every $j\ge0$.

\bpr{\rm {\bf The part $\sseteq$.} By the assumption $(f)\sseteq I$. Therefore $ T_{\cA^{(j)}} f\sseteq I\cdot \RpX$. Thus $\ca_{\cA^{(j)}}\sseteq I\cap \ca_\cA$.
 Therefore (by Artin-Rees or via the primary decomposition) $\ca_{\cA^{(j)}}\sseteq \sqrt{I\cdot\ca_\cA}$.

{\bf The part $\supseteq$.} The quotient $\quots{T_\cA f}{T_\cR f}$ is a finitely-generated module over $R_Y$.
  Thus the submodule $\quots{\ca_\cA\cdot \RpX+T_\cR f}{T_\cR f}\sseteq \quots{T_\cA f}{T_\cR f}$  is finitely generated,
    see \S\ref{Sec.Preliminaries.Rings.Germs}.vi. Take the filtration $(I\cdot \ca_\cA)^\bullet\cdot \quots{\ca_\cA\cdot \RpX+T_\cR f}{T_\cR f}$.
    It is strictly decreasing.
     Therefore for each $i\ge1$ there exists $d_i<\infty$ satisfying:
\[(I\cdot \ca_\cA)^{d_i}\cdot \quot{\ca_\cA\cdot \RpX+T_\cR f}{T_\cR f}\sset f^\#(y)^{i+1}\cdot \quot{T_\cA f}{T_\cR f}.
\]
      Hence $(I\cdot \ca_\cA)^{d_i+1}\cdot \RpX\sset T_\cR f+T_{\cL^{(i)}}f$. Then
  $(I\cdot \ca_\cA)^{d_i+1}\cdot \RpX\sset T_\cR f\cap I^{i+1}\cdot \RpX+T_{\cL^{(i)}}f$, because
   $(f)\sseteq I$.   Artin-Rees over $R_X$ gives:
  $(I\cdot \ca_\cA)^{d}\cdot \RpX\sset T_{\cR^{(j)}} f+T_{\cL^{(j)}}f$ for each $j\ge0$ and $d\gg j$. Hence
  $ \sqrt{\ca_{\cA^{(j)}}}\supseteq I\cdot \ca_{\cA}$.
  \epr
 }
\eee
\end{Properties}

\subsubsection{}
\noindent The following technical lemma is used to apply theorem \ref{Thm.IFT.A.case} in particular cases.
\bel
\bee
\item  $(\ca_\cR+\ca_\cA\cap f^\#(y))^2\cdot \RpX\sseteq (\ca_\cR+\ca_\cA\cap f^\#(y))\cdot T_\cR f+(\ca_\cA\cap f^\#(y))\cdot T_\cL f$.
\item   $(\ca_{\cR^{(0)}}+\ca_{\cA^{(0)}}\cap f^\#(y))^2\cdot \RpX\sseteq (\ca_{\cR^{(0)}}+\ca_{\cA^{(0)}}\cap f^\#(y))\cdot T_{\cR^{(0)}} f+
 (\ca_{\cA^{(0)}}\cap f^\#(y))\cdot  T_{\cL^{(0)}} f$.
\eee
\eel
\bpr (of the second statement) Take the largest $R_Y$-submodule $M\sseteq T_{\cA^{(0)}}f$ such that $M\sseteq \cm\cdot \RpX$ is an $R_X$-submodule.
 (Thus $M$ is uniquely defined.) In particular $M\supseteq T_{\cR^{(0)}}f+\ca_{\cA^{(0)}}\cdot \RpX$.
   Therefore we can present $M=T_{\cR^{(0)}}f+\Lambda$ for some $R_Y$-submodule $\Lambda \sseteq T_{\cL^{(0)}}f$. We take $\Lambda$ the largest possible,
 thus $\Lambda\supseteq  \ca_{\cA^{(0)}}\RpX\cap T_{\cL^{(0)}}f.$ Therefore $Ann(\quots{R^{\oplus p}_Y}{\Lambda})\supseteq \ca_{\cA^{(0)}}\cap f^\#(y).$

 We claim: $Ann(\quots{R^{\oplus p}_Y}{\Lambda})=\ca_{\cA^{(0)}}\cap f^\#(y)$. Indeed, suppose $\Lambda\supseteq \cb_Y\cdot R^{\oplus p}_Y$
   for an ideal  $\cb_Y\sset f^\#(R_Y)$. Then $M\supseteq R_X\cdot (\cb_Y\cdot R^{\oplus p}_Y)$. And thus $\ca_{\cA^{(0)}}\supseteq R_X\cdot \cb_y$.
    Hence $\cb_Y\sseteq \ca_{\cA^{(0)}}\cap f^\#(R_Y)$.

Altogether: $(\ca_{\cA^{(0)}}\cap f^\#(y))\cdot \RpX\sseteq M= T_{\cR^{(0)}}f+\Lambda$. Therefore
\beq
(\ca_{\cR^{(0)}}+\ca_{\cA^{(0)}}\cap f^\#(y))^2\cdot \RpX=(\ca_{\cR^{(0)}}+\ca_{\cA^{(0)}}\cap f^\#(y))\cdot \ca_{\cR^{(0)}}\cdot \RpX+
 (\ca_{\cA^{(0)}}\cap f^\#(y))^2\cdot \RpX\sseteq
\eeq
\[\hspace{8cm}
\sseteq  (\ca_{\cR^{(0)}}+\ca_{\cA^{(0)}}\cap f^\#(y))\cdot T_{\cR^{(0)}}+(\ca_{\cA^{(0)}}\cap f^\#(y))\cdot \Lambda.
\hspace{1cm} \epr\]

\subsection{Extending the module structure of $T_{\cA^{(\bullet)}}f$}
The tangent spaces $T_{\cA^{(\bullet)}}f\sset \RpX$ are $R_Y$-submodules. In fact they are modules over a larger ring, and are often finitely generated.
\bel\label{Thm.TA.module.over.Ry.a}
\bee[\bf 1.]
\item    $\ca_\cA\cdot T_{\cA} f\sseteq T_{\cA}f$.
\item         $\ca_\cA\cdot T_{\cA^{(j)}}f\sseteq T_{\cA^{(j)}}f$ for each $j\ge0$.
\\Therefore all the   spaces $T_{\cA^{(j)}}f$ are modules over the  subring $ f^\#(R_Y)+\ca_\cA\sseteq R_X$.
\item The map $f$ is of finite singularity type (i.e. the  $\k$-module $T^1_\cK f$ is  finitely generated, \S\ref{Sec.3.Finite.Sing.Type})
  \iff the $\k$-module $\quots{R_X}{(f)+\ca_\cA}$ is finitely generated.
\item If $f$  is of finite singularity type, then  any  $(f^\#(R_Y)+\ca_\cA)$-submodule of $\RpX$ is finitely generated.
\\In particular, $\RpX$ and $T_\cA f$ are f.g. over
  the local Noetherian ring $f^\#(R_Y)+\ca_\cA$.
 \eee
\eel
\bpr
\bee[\bf 1.]
\item Immediate verification: $\ca_\cA\cdot T_{\cA} f\sseteq \ca_\cA\cdot \RpX\sseteq T_{\cA}f$.

\item We have:
 $\ca_\cA\cdot T_{\cA^{(j)}}f\stackrel{\S\ref{Sec.3.Filtrations.on.G.TG}}{=}\ca_\cA\big( T_{\cR^{(j)}}f+f^\#(\cb_{j})\cdot T_{\cL}f\big)\sseteq
   T_{\cR^{(j)}}f+f^\#(\cb_{j})\cdot T_{\cA}f$. Note that $f^\#(\cb_{j })\sseteq I^{j+1}$, as $(f)\!\sseteq\! I$.
 Thus   \!$f^\#(\cb_{j })\!\cdot\! T_{\cR}f\!\sseteq \! T_{\cR^{(j)}}f.$
  In addition \!$f^\#(\cb_{j })\!\cdot\! T_{\cL}f\!\sseteq \!  T_{\cL^{(j)}}f.$
  Therefore $ f^\#(\cb_{j })\!\cdot\! T_{\cA}f\!\sseteq\! T_{\cA^{(j)}}f.$

\item
The subring $f^\#(R_Y)+\ca_\cA\sseteq R_X$  is local as the rings $R_Y$, $R_X$ are local and $\ca_\cA\sseteq R_X$ is an ideal.

We can assume $\ca_\cA\sseteq\cm$. Then, by property \ref{Thm.Annihilators.Basic.Properties}.iv we have $\sqrt{(f)+\ca_\cA}=\sqrt{\ca_\cK}$.
 Hence $\quots{R_X}{(f)+\ca_\cA}$ is $\k$-finite iff $\sqrt{(f)+\ca_\cA}\supseteq \cm$ iff $T^1_\cK f$ is $\k$-finite.

\item
By part 3:  $\quots{R_X}{\ca_\cA+(f)}$ is finitely generated over $\k$.

 Applying Weierstra\ss\ finiteness to $\quots{R_X}{\ca_\cA}\otimes\quots{R_Y}{(y)}=\quots{R_X}{\ca_\cA+(f)}$ we get: $\quots{R_X}{\ca_\cA}$ is f.g. over $ R_Y $.
 And then $R_X$  is a f.g. module over $f^\#(R_Y)+\ca_\cA$.
  (Fix a finite set of $R_X$-preimages of the generators of  $\quots{R_X}{\ca_\cA}\in mod\text{-} f^\#(R_Y)$. They generate $R_X$ over $f^\#(R_Y)+\ca_\cA$.)

Therefore the subring $f^\#(R_Y)+\ca_\cA\sset R_X$ is Noetherian. (This is the Eakin-Nagata theorem, \cite[Exercise A.3.7,  pg.625]{Eisenbud-book}.)

Finally, apply the remark of \S\ref{Sec.Preliminaries.Rings.Germs}.viii to $\RpX$ as a f.g. module over  $f^\#(R_Y)+\ca_\cA$.
\epr
\eee

\subsection{Nakayama-type results for $T_\cA f$}\label{Sec.TA.Nakayama}
 Suppose a submodule $M\sseteq \RpX$ satisfies  $T_\cG f+\cm\cdot M\supseteq M$. For $\cG= \cR,\cK$ one
 concludes (by Nakayama): $T_\cG f\supseteq M$. But $T_\cA f$ is not an $R_X$-module.
  If considered as $R_Y$-modules,   $T_\cA f,M$
  are not always  finitely generated. Recall the ``adjusted Nakayama lemma".
\bel\label{Thm.Nakayama.mixed.modules} \cite{Gaffney},\cite[\!lemma \!1.6]{Wall-1981}
 Take a homomorphism of (local, Noetherian) rings $f\!:\!\! (R_Y,\cm_Y)\!\! \to\!\! (R_X,\cm)$,
 with Weierstra\ss\ finiteness condition. \!\!Take f.g. \!\!modules $M_X\!\!\in\! mod \text{-} R_X$, $M_Y\!\!\in\! mod \text{-} R_Y.$
\bee
\item   If $M_X\sseteq M_Y+(f)\cdot M_X$   then $M_X\sseteq M_Y$.
  \item Suppose the inclusion $  M_X\sseteq M_Y+ \cb\cdot M_X$ holds for an ideal  $\cb\sseteq\cm$ satisfying:
   $\cb\cdot M_Y\sseteq  (f)\cdot M_X$.
   Then $  M_X\sseteq M_Y$.
\eee\eel
\bpr
\bee
\item One has $\quots{R_Y}{(y)}\otimes M_X\cong \quots{M_X}{(f)M_X}\sseteq\quots{M_Y+(f)M_X}{(f)M_X}$. Here $\quots{M_Y+(f)M_X}{(f)M_X}$ is f.g. over $R_Y$.
 Thus $\quots{R_Y}{(y)}\otimes M_X$ is f.g. over $R_Y$, see \S\ref{Sec.Preliminaries.Rings.Germs}.vi.
  By Weierstra\ss\ finiteness $M_X$ is f.g. over $R_Y$. Finally, apply Nakayama over $R_Y$ to the inclusion $M_X\sseteq M_Y+(f)\cdot M_X$.
\item
One gets  $  \cb\cdot M_X\sseteq  (f)\cdot M_X+  \cb^2\cdot M_X$. Nakayama (over $R_X$) gives
 $  \cb\cdot M_X\!\sseteq\!  (f)\cdot M_X.$ The initial assumption gives  $  M_X\!\sseteq\! M_Y\!+\! \!(f)\!\cdot\! M_X$. Apply part one.
\epr\eee

\bcor\label{Thm.Nakayama.For.TA}
 Fix some $j\ge-1$, and $M\sseteq\RpX$, and $(f)\sseteq I$.
\bee[\bf 1.]
\item If $T_{\cA^{(j)}}f+ (f) \cdot M\supseteq M$  then $T_{\cA^{(j)}}f \supseteq M$.
\item Suppose $\ca_\cA\sseteq\cm$. If  $T_{\cA^{(j)}}f+(\ca_\cA+(f))\cdot M\supseteq M$ then $T_{\cA^{(j)}}f \supseteq M$.
\eee
\ecor
\bpr As we need an upper bound on $M$, we can assume $M\supseteq T_{\cR^{(j)}}f$.
\bee[\bf 1.]
\item Pass to the quotients, $\quots{T_{\cA^{(j)}}f}{T_{\cR^{(j)}}f}+(f)\cdot\quots{M}{T_{\cR^{(j)}}f}\supseteq \quots{M}{T_{\cR^{(j)}}f}$.
 Here $\quots{M}{T_{\cR^{(j)}}f}\in mod\text{-}R_X$ and $\quots{T_{\cA^{(j)}}f}{T_{\cR^{(j)}}f}\in mod\text{-}R_Y$.
  Therefore (by part one of  lemma \ref{Thm.Nakayama.mixed.modules}) $\quots{T_{\cA^{(j)}}f}{T_{\cR^{(j)}}f}\supseteq \quots{M}{T_{\cR^{(j)}}f}$.

\item
 Multiply the assumption by $\ca_\cA$ to get $\ca_\cA\cdot T_{\cA^{(j)}}f+\ca_\cA\cdot (\ca_\cA+(f))\cdot M\!\supseteq\! \ca_\cA\!\cdot\! M.$
  Nakayama over $R_X$ gives
  $\ca_\cA\!\cdot\! T_{\cA^{(j)}}f\!\supseteq \!\ca_\cA\!\cdot\! M.$
 Part two of lemma \ref{Thm.TA.module.over.Ry.a} gives:   $ T_{\cA^{(j)}}f\!\supseteq \!\ca_\cA\!\cdot\! M.$
  The initial assumption becomes:     $T_{\cA^{(j)}}f\!+\! (f)\! \cdot\! M\!\supseteq\! M.$
 Now invoke part 1.
\epr
\eee

\bex Let $(R_X,R_Y)$ be one of the   pairs  $(\quots{\k[[x]]}{J},\k[[y]])$, $(\quots{\k\{x\}}{J},\k\{y\})$, $(\quots{\k\bl x\br}{J},\k\bl y\br)$,
 see \S\ref{Sec.Preliminaries.Rings.Germs}.ii.
\\Let $I=\cm$ and suppose $ {\ca_\cA}\sseteq\cm$.  Take $M=\cm^j\cdot (f)\cdot\RpX$ for some $j\ge0$.
 Part 2 of corollary \ref{Thm.Nakayama.For.TA} gives:
 \beq
 \text{If $T_{\cA^{(j)}}f+(\ca_\cA+(f))\cdot \cm^j\cdot(f)\cdot \RpX\supseteq \cm^j(f)\RpX$ then
  $T_{\cA^{(j)}}f=T_{\cK^{(j)}}f$.}
  \eeq
   For $j=0$ this is well known, e.g. theorem 5.1 of \cite{Ruas}.
\eex

\beR
One would like a stronger bound:
 \!``If  \!$T_{\cA}f\!+\!(\sqrt{\ca_\cA}\!+\!(f))\!\cdot\! M\! \!\supseteq\!\! M$ then $T_{\cA}f \!\!\supseteq \!\!M$". This fails already for the case \!$\sqrt{\ca_\cA}\!=\!\cm.$
For example, one always has $T_\cA f\!+\!\cm\!\cdot \!\RpX\!=\!\RpX,$\! though of course $T_\cA f\!\neq\! \RpX.$
\eeR

\bcor
 Suppose $\cb\!\cdot\! \RpX\!\sseteq\! T_{\cA^{(j)}}f+\cb\cdot(\cb(f):(f)^{j+1})\RpX$, for an ideal $\cb\sset R_X$.
 Then  $\cb\cdot \RpX\!\sseteq\! T_{\cA^{(j)}}f.$
\ecor
 \bpr Pass to the quotients, $\cb\!\cdot\!\quots{ \RpX }{T_{\cR^{(j)}}f}\!\sseteq\! \quots{T_{\cA^{(j)}}f}{T_{\cR^{(j)}}f}
  +\cb\!\cdot\!(\cb(f):(f)^{j+1})\!\cdot\!\quots{\RpX}{T_{\cR^{(j)}}f}.$
   Note:  $(\cb(f):(f)^{j+1})\!\cdot\!\quots{T_{\cA^{(j)}}f}{T_{\cR^{(j)}}f}\!\sseteq\! \cb\!\cdot\!(f)\!\cdot\!\quots{\RpX}{T_{\cR^{(j)}}f}.$
    Then (by part 2 of  lemma \ref{Thm.Nakayama.mixed.modules}): $\cb\!\cdot\!\quots{ \RpX }{T_{\cR^{(j)}}f}\!\sseteq\! \quots{T_{\cA^{(j)}}f}{T_{\cR^{(j)}}f}.$
 \epr

\bex\bee[\bf i.]
\item Take  $\cb=I^r$ and $(f)\sseteq I,$ for some $0\le j<r.$ One gets the well-known bound: if
 $T_{\cA^{(j)}}f+I^{2r-j}\cdot\RpX\supseteq I^r\cdot \RpX$   then $ T_{\cA^{(j)}}f  \supseteq I^r\cdot \RpX$. \quad
 (See e.g. corollary 1.23, pg. 17  of \cite{Ruas}.)
\item Moreover, for any map $g\in \cm\cdot \RpX$ satisfying $g-f\in I^{2r-2j}\cdot \RpX$
 one gets (by the direct check) $ T_{\cA^{(j)}}g  \supseteq I^r\cdot \RpX$. For $I=\cm$ this is \cite[corollary 2.4]{du Plessis}.
\eee
\eex

\subsection{Artin-Rees type properties of $T_\cA f$}\label{Sec.TA.Artin.Rees.Properties}

   The image tangent spaces  $T_\cR f$, $T_\cK f\sseteq\RpX$ are $R_X$-submodules. For them the condition
  ``$T_{\cG^{(j)}}f\supseteq I^d\cdot \RpX$" implies ``$T_{\cG^{(j+\bullet)}}f\supseteq I^{ d+\bullet}\cdot \RpX$", for every $\bullet\ge 0$.
 This is heavily used in the study of $\cR$, $\cK$-orbits, theorems \ref{Thm.Orbits.R.Filtration}, \ref{Thm.Orbits.K.Filtration}.

The tangent space $T_\cA f$ is a module over $ R_Y $,   not over $R_X$.
 The condition $T_{\cA^{(\bullet)}}f\supseteq I^d\cdot \RpX$ does not imply ``$T_{\cA^{(j+\bullet)}}f\supseteq I^{d+\bullet}\cdot  \RpX$ for every $\bullet\ge0$".
 (For example, in the infinitesimally stable case, i.e. $T_\cA f=\RpX$, one does not have $T_{\cA^{(0)}}f\supseteq \cm\cdot \RpX$.)
Yet, some properties of Artin-Rees type do hold.

\

\noindent Take a map $ X\! \stackrel{f}{\to}\! Y=(\k^p,o),$ the corresponding map $ R_Y\!\stackrel{f^\# }{\to}\!R_X,$ and  a pair $mod\text{-}R_Y\ni\! M_Y\!\sset\! M_X\!\in\! mod\text{-}R_X. $
\bel\label{Thm.Artin-Rees.Mixed}
Then $M_Y\cap \cm^{d_j}\cdot M_X\sseteq  \cm^j_Y \cdot M_Y\cap \cm^{d_j}_X\cdot M_X$ for every $j$
 and a corresponding $d_j<\infty$.
\eel
\bpr
The filtration $\cm^\bullet_X\cdot M_X$ is strictly decreasing (by Nakayama over $R_X$). Take the induced filtration
 $\quots{M_Y}{\cm^j_Y\cdot M_Y}\supset \quots{M_Y\cap \cm^\bullet_X\cdot M_X+\cm^j_Y\cdot M_Y}{\cm^j_Y\cdot M_Y}$. Here
 $\quots{M_Y}{\cm^j_Y\cdot M_Y}$ is an Artinian module over the ring $\quots{R_Y}{\cm^j_Y}$. Thus the induced filtration stabilizes,
\beq
\quot{M_Y\cap \cm^d_X\cdot M_X+\cm^j_Y\cdot M_Y}{\cm^j_Y\cdot M_Y}=\cap_{\bullet\ge1} \quot{M_Y\cap \cm^\bullet_X\cdot M_X+\cm^j_Y\cdot M_Y}{\cm^j_Y\cdot M_Y}
\quad \text{ for } \quad d\gg1.
\eeq
Finally, $\cap_\bullet \quots{M_Y\cap \cm^\bullet_X\cdot M_X+\cm^j_Y\cdot M_Y}{\cm^j_Y\cdot M_Y}\sseteq
  \cap_\bullet \quots{ \cm^\bullet_X\cdot M_X+\cm^j_Y\cdot M_Y}{\cm^j_Y\cdot M_Y}\twoheadleftarrow \cap_\bullet (\cm^\bullet_X\cdot M_X)=0$.

   Thus $\quots{M_Y\cap \cm^{d_j}\cdot M_X+\cm_Y\cdot M_Y}{\cm^j_Y\cdot M_Y}=0$ for $d_j\gg1$. Hence the statement.
\epr

\

\bel\label{Thm.Artin.Rees.for.TA} Assume $(f)\sseteq I$.
\bee[\bf 1.]
\item Suppose $T_{\cA^{(j_l)}}f\supseteq \cc_l\cdot \RpX$ for $l=1,2$ and some ideals $\cc_l\sseteq I^{j_l+1}$, with $j_l\ge 0$.
 Then  $ T_{\cA^{(j_1+j_2+1)}}f\supseteq \cc_1\cdot \cc_2\cdot \RpX$.

\item Suppose $T_{\cA^{(j)}}f\supseteq I^{d_j}\cdot \RpX$ and  $T_{\cK^{(i)}}f\supseteq I^{d_i}\cdot \RpX$  for some $d_j\ge i\ge0$ and $d_j>j\ge -1$.
 Then $T_{\cA^{(j+1)}}f\supseteq I^{d_j+d_i-i}\cdot \RpX$.
\\By iterating this we get: $T_{\cA^{(j+k)}}f\supseteq I^{d_j+k(d_i-i)}\cdot \RpX$ for each $k\ge0$.
\item
 Suppose an ideal  $\ca\sset R_X$ and a subset $\cb_Y\sseteq f^\#(y)^2\cap \ca\sset R_X$ satisfy: $\ca^2\sseteq\ca_\cA\cdot \cb_Y$.
    Then $\ca\cdot  T_\cR f+ T_{\cL^{(1)}} f\supseteq  \cb_Y \cdot \ca_\cA\cdot \RpX\supseteq \ca^2\cdot\RpX$.

For a subset $\cb_Y\sseteq f^\#(y)^2\cap (\ca\cdot \cm)\sset R_X$ satisfying  $\ca^2\sseteq\ca_\cA\cdot C$ one has:
  $\ca\cdot \cm\cdot T_\cR f+ T_{\cL^{(1)}} f\supseteq \cb_Y \cdot \ca_\cA\cdot \RpX\supseteq \ca^2\cdot\RpX$.

\item Take a subset $C\sseteq (y)\sset R_Y$.
 Then $ C\cdot \cm\cdot T_\cR f+ T_{\cL^{(1)}} f\supseteq  C\cdot \ca_{\cA^{(o)}}\cdot  \RpX$.

In particular,  $ \ca_{\cA^{(o)}}\cdot T_{\cR^{(0)}} f+ T_{\cL^{(1)}} f\supseteq   ( \ca_{\cA^{(o)}}\cap f^\#(R_Y)+ \ca_{\cR^{(o)}}) \ca_{\cA^{(o)}}\cdot  \RpX$.

And also:
 $(\ca_{\cA}\cap f^\#(R_Y)+ \ca_{\cR^{(o)}})\cdot T_{\cR^{(0)}} f+(\ca_{\cA}\cap y)^2\cdot T_\cL f\supseteq
   ( \ca_{\cA}\cap f^\#(R_Y)+ \ca_{\cR^{(o)}}) ^2\cdot  \RpX$.

    $(\ca_{\cA}\cap f^\#(R_Y)+ \ca_\cR)\cdot T_{\cR^{(0)}} f+(\ca_{\cA}\cap y)^2\cdot T_\cL f\supseteq
   ( \ca_{\cA}\cap f^\#(R_Y)+ \ca_\cR) ^2\cdot  \RpX$.

  \item
Suppose $(f)\sseteq I$ then
$\ca_\cK\cdot(\ca_{\cA^{(j)}}\cap I^{j+2})\sseteq \ca_{\cA^{(j+1)}}$ for each $j\ge-1$.
\\
If $(f)\sseteq I$ and $\ca_{\cA^{(j)}}\cdot T_\cR f\sseteq  T_{\cR^{(j+1)}} f$  then
$\ca_\cK\cdot \ca_{\cA^{(j)}} \sseteq \ca_{\cA^{(j+1)}}$ for each $j\ge-1$.
\eee
\eel
We remark that in part 2 the sequence $\{d_i-i\}$  is non-increasing  and $d_i-i\ge ord(f)$.
\bpr  All the proofs are direct verifications.
\bee[\bf 1.]
\item
$\cc_1\cdot\cc_2\cdot \RpX\sseteq \cc_1\cdot \Big(T_{\cR^{(j_2)}} f+f^\#(\cb_{j_2})\cdot T_{\cL}f\Big)\sseteq T_{\cR^{(j_2+j_1+1)}} f
 + f^\#(\cb_{j_2})\cdot \Big(T_{\cR^{(j_1)}} f+ $
\\$ +f^\#(\cb_{j_1})\cdot T_{\cL}f\Big)
\sseteq T_{\cR^{(j_1+j_2+1)}} f+f^\#(\cb_{j_1+j_2+1})\cdot T_{\cL}f=T_{\cA^{(j_1+j_2+1)}} f.$

\item
$I^{d_j+d_i-i}\cdot \RpX\sseteq I^{d_j-i}\cdot \left(T_{\cR^{(i)}}f+I^i\cdot (f)\RpX\right) \sseteq  T_{\cR^{(d_j)}}f+I^{d_j}\cdot (f)\RpX
\sseteq  T_{\cR^{(d_j)}}f+(f)\cdot  T_{\cR^{(j)}}f+(y)\cdot  T_{\cL^{(j)}}f \sseteq T_{\cA^{(j+1)}}f$.

\item
Multiply the inclusion $T_\cA f\supseteq \ca_\cA\cdot \RpX$ by $\cb_Y$.

\item
Multiply the inclusion $T_{\cA^{(o)}} f\supseteq \ca_{\cA^{(o)}}\cdot \RpX$ by $C$.
\item
 $\ca_\cK\!\cdot\!(\ca_{\cA^{(j)}}\!\cap\! I^{j+2})\!\cdot\! \RpX\!\sseteq \!(\ca_{\cA^{(j)}}\!\cap\! I^{j+2})\!\cdot \!(T_\cR f\!+\! (f)\!\cdot\! \RpX)\!\sseteq\!
 T_{\cR^{(j+1)}}f\!+\!(f)\!\cdot\! T_{\cR^{(j)}}f\!+\!  (y)\!\cdot\! T_{\cL^{(j)}}f\!\sseteq\! T_{\cA^{(j+1)}}f.$
\epr
\eee

\bex
\bee[\bf i.]
\item Part 1 gives: $\ca^2_{\cA^{(0)}}\sseteq \ca_{\cA^{(1)}}$.
\item Let $I=\cm$ and  $\ca\sseteq T_{\cA^{(0)}}f$ and $(f)\sseteq \ca$. Then $\ca^2\cdot \RpX\sseteq I\cdot \cm\cdot T_\cR f+ T_{\cL^{(1)}} f$.
\eee
\eex

\subsection{Filtration closures, $\overline{T_\cA f}$ vs $T_\cA f$.}\label{Sec.TA.filtration.closure}
 The tangent spaces $T_\cR f, T_\cK f$ are filtration-closed, by  lemma \ref{Thm.Filtration.Closures.TR.TK}.
 As was explained in \S\ref{Sec.Intro.Contents}, the methods of that lemma are not directly applicable to $\cA$-equivalence.
    We give the $\cA$-version. Fix some $j\ge-1$.
\bprop\label{Thm.Filtration.Closures.TA}
\bee[\bf 1.]
\item
  $\ca_\cA\cdot \overline{T_{\cA^{(j)}} f}=\ca_\cA\cdot T_{\cA^{(j)}} f$.
\item
 Suppose either $R_X=\quots{\k[[x]]}{J} $ or  $f$ is of finite singularity type, \S\ref{Sec.3.Finite.Sing.Type}.
 Then $T_{\cA^{(j)}}f=\overline{T_{\cA^{(j)}}f}$.
\eee
\eprop
\bpr
\bee[\bf\!\!1.\!\!]
\item
Observe the obvious inclusions:  $\ca_\cA\!\cdot\! T_{\cA^{(j)}} f\!\sseteq \! \ca_\cA\!\cdot\!\overline{T_{\cA^{(j)}} f}\!\sseteq\! \overline{\ca_\cA\cdot T_{\cA^{(j)}} f}.$
Note that $\ca_\cA\!\cdot\! T_{\cA^{(j)}}f\!\sseteq\! \RpX$ is an $R_X$-module. Therefore it is closed (lemma \ref{Thm.Filtration.Closures.TR.TK}),
 $\ca_\cA\cdot T_{\cA^{(j)}}f=\overline{\ca_\cA\cdot T_{\cA^{(j)}}f}$. Hence
$\ca_\cA\cdot T_{\cA^{(j)}} f=  \ca_\cA\cdot\overline{T_{\cA^{(j)}} f}$.

\item
The case $R_X=\quots{\k[[x]]}{J}$ is trivial.  Suppose $f$ is of finite singularity type.

We have the $(f^\#(R_Y)+\ca_\cA)$-submodule $\overline{T_\cA f}\sset \RpX$, by lemma \ref{Thm.TA.module.over.Ry.a}.
   It satisfies: $T_\cA+I^d\cdot \RpX\supseteq \overline{T_\cA f}\supseteq T_\cA f$ for any $d\gg1$. Therefore $T_\cA+(I^d\cdot \RpX)\cap \overline{T_\cA f}\supseteq \overline{T_\cA f}$.

The ring $f^\#(R_Y)+\ca_\cA$ is local Noetherian (lemma \ref{Thm.TA.module.over.Ry.a}),
 therefore $\overline{T_\cA f}$ is f.g. over this ring (\S\ref{Sec.Preliminaries.Rings.Germs}.vi).

As $f$ is of finite singularity type, $\sqrt{\ca_\cK}=\cm\supseteq I.$
   For $d\gg1$ one has $I^d\sseteq (f)+\ca_\cA$, by property \ref{Thm.Annihilators.Basic.Properties}.iv.
     By lemma \ref{Thm.Artin-Rees.Mixed} one has:
    $ ((f^\#(y)+\ca_\cA)^d\cdot \RpX)\cap \overline{T_\cA f}\sseteq (f^\#(y)+\ca_\cA)\cdot \overline{T_\cA f}$.
    Thus we get $T_\cA f+(f^\#(y)+\ca_\cA)\cdot \overline{T_\cA f}\supseteq \overline{T_\cA f}$.
     Apply Nakayama to  $T_\cA f$ as a f.g. module over  $f^\#(R_Y)+\ca_\cA$.
\epr
\eee

\bex
Let $R_X$ be one of   $ \quots{\k[[x]]}{J} ,  \quots{\k\{x\}}{J} ,  \quots{\k\bl x\br}{J} $, see \S\ref{Sec.Preliminaries.Rings.Germs}.ii.
 If $f$ is $\cK$-finite then $T_{\cA^{(j)}} f=\overline{T_{\cA^{(j)}} f}$.
\eex

\section{Criteria for left-right orbits, ``$\cA f$ vs $T_\cA f$"}\label{Sec.A.Orbits}

\subsection{An auxiliary result:  orbits of unipotent groups are Zariski closed}\label{Sec.Preliminaries.Orbits.Unipotent.Group.Closed}
Let $G$ be an algebraic (connected)   unipotent  group
over a field $\k.$
Suppose $G$ acts on the    affine space  $\A^N_\k.$
If $\k$ is algebraically closed then the $G$-orbits are Zariski-closed.
 This is the Kostant-Rosenlicht theorem, see e.g. \cite[Theorem 2.11]{Ferrer Santos-Rittatore}.

For the $\cA$-implicit function theorem below (with $\k$ not algebraically closed) we need the following statement,  provided by
 Mikhail Borovoi.
\bel\label{Thm.Orbit.Unipotent.Group.Closed}
Let $G$ be a (algebraic, connected)   unipotent  group
over a field $\k$ of characteristic 0.
 Suppose a line  $L\subseteq \A^N_\k$ intersects an orbit
  $G x_o\sset  \A^N_\k$ in an infinite number of points.
Then $G x_o\supseteq L$.
\eel
\bpr
 Take the  algebraic closure  $\k\sseteq\bk$ and the corresponding orbit  $G_{\bar\k}\cdot x_o\sseteq \A^N_{\bar\k}.$
Thus $G_{\bar\k}\cdot x_o$ is Zariski-closed (by Kostant-Rosenlicht),   hence  contains $L$.
The orbit  $G_{\bar\k}\cdot x_o$ contains the $\k$-point $x_o$,
 hence is stable under the action of the Galois group ${\rm Gal}(\bar\k/\k)$.
We conclude (by Galois descent) that $G_{\bar\k}\cdot x_o$ is defined over $\k$,
see, e.g., \cite[Proposition 2.6]{Jahnel}.

  Denote by $G_0$ the stabilizer of $x_o$ in $G$. Then $G_0$ is a unipotent $\k$-group.
As $char(\k)=0$, $G_0$ is isomorphic, as a variety, to its Lie algebra ${\rm Lie}(G_0)$ and hence  is connected,
see, e.g., \cite[Lemma (3.1)]{Bruce.du-Plessis.Wall}.
By   \cite[Lemma 1.13]{Sansuc}, see also \cite[Proposition 3.1]{BDR},
 we have  $H^1(\k, G_0)=1$.
The set of orbits of $G_\k$ in $(G_{\bar\k}\cdot x_o)\cap \A^N_\k$
is in a canonical bijection with $\ker[H^1(\k,G_0)\to H^1(\k,G)]$,
see \cite[Section I.5.4, Proposition 36]{Serre}.
As $H^1(\k, G_0)=1$, we see that there is only one orbit in $(G_{\bar\k}\cdot x_o)\cap \A^N_\k$,
 that is, $(G_{\bar\k}\cdot x_o)\cap \A^N_\k=G_\k\cdot x_o$.
Therefore $G_\k\cdot x_o \supseteq L $, as required.
\epr

\

Lemma \ref{Thm.Orbit.Unipotent.Group.Closed} does not hold in positive characteristic, even if $\k$ is perfect.
\bex (by Zev Rosegarten) Let $\k$ be an infinite field of characteristic $p$, such that the Artin-Schreier map, $\k\!\ni\! c\!\to\! c^p\!-c\!\in\! \k$ is not surjective.
 Equivalently, $H^1(\k,\!\quots{\Z}{p\Z})\!\neq\!0.$\! One can even assume that $\k$ is perfect. (E.g. let $\k$ be the perfect closure of a local or global function field.) Take the action (of the additive group) \vspace{-0.3cm}
 \beq
 G_a:=\k\circlearrowright \A^2, \quad  \la\cdot \bbm x_1\\x_2\ebm:=\bbm 1&\la^p-\la\\0&1\ebm\cdot\bbm x_1\\x_2\ebm.
 \eeq
This action is algebraic and unipotent. It preserves the line $x_2=1$. And the orbit $G_a(0,1)$ is infinite.
 But this orbit does not contain this line, as the Artin-Schreier map is non-surjective.
\eex

\subsection{The  \!$\cA$-implicit function theorem}\label{Sec.A.IFT}
Let $R_X$ be one of $ \quots{\k[[x]]}{J} ,  \quots{\k\{x\}}{J} ,  \quots{\k\bl x\br}{J} $,
  \S\ref{Sec.Preliminaries.Rings.Germs}.ii.
 Denote by $\cm_\k$ the maximal ideal of $\k.$

  Take a map $f:X\to (\k^p,o)$. If its critical locus, $Crit(f)$, see \S\ref{Sec.3.Annihilators},  is a point, i.e. $\sqrt{\ca_\cR}=\cm$, then already theorem \ref{Thm.IFT.R.case}
  gives a good control on the orbit, e.g. $\cR f\supseteq \{f\}+  \ca^2_\cR\cdot\RpX$. For $p\ge2$ (and $f$ not a submersion)
   the critical locus is always of positive dimension. And ``in most cases" the instability locus, $V(\ca_\cA)$, is of
    dimension smaller than $Crit(f)$. We adopt the following assumptions.
\begin{Assumptions}\label{Assumptions}
\bee[\bf i.]
\item
 Take an ideal $\cm\cdot\ca_\cR\sseteq \ca\sseteq\cm$, so that $V(\ca)\sseteq Crit(f).$
  Suppose either $Crit(f)=V(\cm)\sset  X$   or the locus $V(\ca)\sset Crit(f)$ contains no  irreducible component of $Crit(f)\sset X$.
\item
If $J\!\neq\!0$ and $\k\not\supseteq\Q$ then we assume the $jet_0$-condition of \S\ref{Sec.Preliminaries.Exp.Ln}.
\eee
\end{Assumptions}
\noindent Moreover,  for the ring  $R_X\!=\!\quots{\k\{x\}}{J}$ (with $J\!\neq\!0$)  the statements   are for the  $\cm$-adic
closure,  $\overline{\cA f}\!\supseteq\! \{f\}\!+\!\cdots.$

   Recall the tangent space of the left equivalence, $T_{\cL^{(1)}} f=f^\sharp((y)^2\cdot \RpY)$, see \S\ref{Sec.3.Filtrations.on.A.TA}.

\bthe\label{Thm.IFT.A.case}
Suppose  $\ca^2\cdot \RpX\sseteq \ca\cdot \cm\cdot T_\cR f+ T_{\cL^{(1)}} f$.
\bee[\bf a.\!\!]
\item 
 Then $ {\cA f}\supseteq \{f\}+(\cm_\k+\ca^2+(f))\cdot\ca^2\cdot \RpX+ T_{\cL^{(1)}} f$.
\item  
   Suppose  the field $\quots{\k}{\cm_\k}$ is   of characteristic  zero or is  algebraically closed.
     \!\!Then
$\cA f \!\supseteq \!\{f\}\!+\!\ca^2\!\cdot\! \RpX\!+\! T_{\cL^{(1)}} f.$
\eee
\ethe

The trivial choice $\ca=\cm\cdot\ca_\cR$ satisfies the assumptions. Then the statement is implied by theorem \ref{Thm.IFT.R.case}.

We remark: the assumption   $\ca^2\cdot \RpX\sseteq \cm \cdot\ca\cdot T_\cR f+ T_{\cL^{(1)}} f$ is preserved by $\cA$-equivalence.
\bpr The proof is an essentially more  complicated version of the proofs of theorems \ref{Thm.IFT.R.case}, \ref{Thm.IFT.K.case}.
\bei
\item  The condition $f+g\in \cA f$ is not an  implicit function equation.
 Yet, in Step 2 we reduce it to a system of implicit function equations over the ring $R_Y$,  \eqref{Eq.IFT.A.to.prove.expanded.form},
\eqref{Eq.IFT.A.ift.to.resolve}, modulo higher order terms, $\ca^d_\cR$ for $d\gg1.$ Here the additional tool is the Weierstra\ss\ division.
\item Part a., for the case $J=0$,  is  proved  in Step 3  as in the $\cR,\cK$-cases, by applying $IFT_\one$.
\item Part b. (Step 4) is more involved. Unlike the $\cR,\cK$-cases we cannot resolve this system directly.
\bei
\item Step 4.i. gives  an ``arc-solution", $f+tg\in \cA_t f$. This arc is formal, with parametrization in $\k[\![t]\!]$.
 \item  Step 4.ii. gives: the condition
  $f\!+\!tg\!\in\! \cA_t f$ is resolvable for an infinite number of values of $t\!\in\! \quots{\k}{\cm_\k}.$
\item Step 4.iii.  To achieve ``$f+tg\in \cA_t f$ for any $t\in \k$"  we pass  to
 finite jets and use  the theorem  of Kostant-Rosenlicht or Lemma \ref{Thm.Orbit.Unipotent.Group.Closed}.
 \eei
 \item The case $J\neq0$ is then obtained iteratively, using $jet_o$-assumption,  in a way similar to theorems \ref{Thm.IFT.R.case}, \ref{Thm.IFT.K.case}. This gives the formal $\cA$-equivalence. In the Nash case, $\quots{\k\bl x\br}{J},$ one invokes the left-right Artin approximation, \S\ref{Sec.Preliminaries.Basic.Tools}.v. In the analytic case, $\quots{\k\{x\}}{J}$,   the statement is only for the orbit-closures.
\eei
\bee[\bf\!Step 1.\!\!]
\item (Preparations)
\bee[\!\!\bf i.\!\!]
\item Consider the \!$R_X$-module \!$M\!\!:=\!\!\quots{\ca^2\!\cdot\! \RpX\!+\!\ca\!\cdot \!\cm\!\cdot \!T_\cR f}{\ca\!\cdot\!\cm\!\cdot  \!T_\cR f}$.   \!Take  its annihilator \mbox{ideal \!$\ca_M\!\!:=\!Ann(M)\!\sset\! R_X.$}
  We have $\ca_\cR\!\cdot \!M\!=\!0$, thus $\ca_M\supseteq \ca_\cR.$  Geometrically: $Supp(M)\sseteq Crit(f)=V(\ca_\cR)\sset X.$

\quad
We claim: $\sqrt{\ca_M}=\sqrt{\ca_\cR}$. This is obvious for $\ca_\cR=\cm$ (i.e. $Crit(f)=V(\cm)=o\in X $). Thus we assume: $dim Crit(f)>0$ and
  $V(\ca)$ contains no  irreducible component of $Crit(f)$.
   Obviously   $Supp(M)\smin V(\ca)= V(\ca_\cR)\smin V(\ca)$. (In fact, over $X\smin V(\ca)$ one has: $M\cong T^1_\cR f.$)
 And $Supp(M)\sset X$ is closed. Therefore   $Supp(M)= Crit(f)$, i.e.
 $\sqrt{\ca_M}=\sqrt{\ca_\cR}$.

\item
Consider $M$ as an $R_Y$-module. By our assumption $M$ is a submodule of
 $\quots{  T_{\cL^{(1)}} f+\cm\ca \cdot T_\cR f}{\cm\ca\cdot  T_\cR f}$. The later module is finitely generated over $R_Y$. Therefore $M$ is f.g.
  over $R_Y$, \S\ref{Sec.Preliminaries.Rings.Germs}.vi.

\quad   Each $x_i\!\in\! R_X$ acts on $M$ as an $R_Y$-linear operator. And its action is filtration-nilpotent,
   i.e. $x_i^d\!\cdot\! M\!\sseteq\! f^\#(y)\!\cdot\! M$ for $d\!\gg\!1$. Therefore the
    characteristic polynomial of $x_i$ is of the form
    \beq
    \chi_{x_i}(t):=det[t\one-x_i]=t^{d_i}+y(\dots)\in R_Y[t].
\eeq
     One has $\chi_{x_i}(x_i)=0\circlearrowright M$,
     by Cayley-Hamilton theorem. (Hence $x^{d_i}_i\in (f)+\ca_M.$)
     Therefore we can apply the Weierstra\ss\ division, \S\ref{Sec.Preliminaries.Rings.Germs}.vii, and present
       $R_X=\ca_M+Span_{R_Y} \{v_\bullet\} $. Here $\{v_\bullet\}$ is a finite tuple whose image generates the $R_Y$-module
      $\quots{R_X}{\ca_M}$.

      We   take $v_0=1$ and $v_{l>0}\in \cm$.

\item We claim: $\{f\}+\cc\cdot \RpX+ T_{\cL^{(1)}} f\sseteq \cL(\{f\}+\cc\cdot \RpX)$, for any ideal $\cc\sset R_X$.
 Indeed, given an element $q(y)\in T_{\cL^{(1)}} f$, define the automorphism $\Phi_Y\in \cL^{(1)}$ by $\Phi_Y(y)=y+q(y)$.
  Then $\Phi^{-1}_Y(\{f\}+\cc\cdot \RpX+q(f))\sseteq \{f\}+\cc\cdot \RpX$.

In addition   $\cR f\supseteq \{f\}+\cm^2\cdot \ca^3_\cR\cdot \RpX$, by Theorem \ref{Thm.IFT.R.case}.

In view of the previous observations,
 it is enough to prove (for part b.):
 \beq\label{Eq.IFT.A.to.prove}
\quad\quad \{f\}+\ca^2\cdot \RpX\sseteq \cA\Big(\{f\}+\ca^d_M\cdot \RpX+ T_{\cL^{(1)}} f\Big)\quad\quad\quad \text{ for }\quad d\gg1.
\eeq

\eee

{\hspace{-1.8cm}\em From now and until Step 5, we assume $J=0$, i.e. $R_X$ is one of the rings $\k[[x]], \k\{x\},\k\bl x\br.$}
\item (Reduction of the condition \eqref{Eq.IFT.A.to.prove} to a system of implicit function equations)

As in the $\cR,\cK$-cases we fix some generators   $\{a^M_m\}$ of $\ca_M$  and $\{\xi_i(f)\}$ of $\ca\cdot T_\cR f$ (as $R_X$-modules).
 Recall the expansion $R_X=\ca_M+Span_{R_Y}\{v_\bullet\}$ of Step 1.ii.

\bee[\bf\!\!i.\!\!]
\item ({\em The first attempt.})
 Take a perturbation   $g\in \ca^2\cdot\RpX$. Expand it, $g-\sum   c^g_{i }\xi_i(f)\in T_{\cL^{(1)}} f$, here $c^g_{i}\in \cm$.
 Accordingly  define  the coordinate change $\Phi_X\in \cR$  by $x\to x+ \sum  c_{i}\xi_i(x)$, here $\{c_{i}\}\in R_X$ are unknowns.
  As in the proofs of $\cR,\cK$-cases we have:
\beq
 \Phi_X(f)-f-  \sum    c_{i}\cdot\xi_i(f)\in    \ca^2\!\cdot\!\{c_{*}\}^2\!\cdot\!\RpX,\quad \quad \quad
  \Phi_X(g)-g\in \big[\ca\!\cdot \!\{c_*\}\!\cdot\! (\ca^2)'+\ca^2\!\cdot\!(\ca^2)''\{c_*\}^2\big]\!\cdot\! \RpX.
\eeq
We cannot use the assumption $\ca^2 \RpX\sseteq \cm\cdot \ca\cdot  T_\cR f+T_{\cL^{(1)}}f$ directly, as  $\cm\cdot \ca\cdot  T_\cR f+T_{\cL^{(1)}}f$ is not an $R_X$-module. Instead we expand $c_{i}=\sum_l c_{i,l}v_l+\sum_m \tc_{i,m}\ca^M_m$, with the unknowns $c_{i,l}\in R_Y$, $\tc_{i,m}\in R_X$. (Using $R_X=\ca_M+ Span\{v_\bullet\}$.)
   Similarly one expands $\{c^g_{i}\}$.
  The total expression is:
\beq\label{Eq.IFT.A.to.prove.current}
\Phi_X(f+g)=f+f^\#(y)^2\cdot  H(\{c_{**}\})+ \hspace{9cm}
\eeq
\vspace{-0.4cm}\[\hspace{3cm}
+\sum_{i} \Big[ \sum_l v_l (c_{i, l}+c^g_{i,  l})+\sum_m a^M_m(\tc_{i, m}+\tc^g_{i,  m})+ Q_{i}(\{c_{*}\})+
 \ch_{i}(\{c_{*}\},\{\tc_{*}\})\Big]  \xi_i(f) .
\]
Here $Q_{i}(z)\in (z)^2\cdot R_Y[[z]]$, \ $\ch_{i}(z,w)\in (z,w)^2\cdot R_Y[[z,w]]$ and $H\in T_\cL f$.

We would like to get rid of the term  $[\dots]$, i.e. the coefficient of $ \xi_i(f)$. This cannot be done directly. Instead we will
 force this coefficient to belong to $\ca^d_M$ for $d\gg1$.
 This will ensure the condition \eqref{Eq.IFT.A.to.prove}.

\item ({\em The actual reduction.})
Iterate the expansion $R_X=\ca_M+Span_{R_Y}\{v_\bullet\}$ to  get $R_X=\ca^{d}_M+\sum^{d-1}_{k=0} \ca^k_M\cdot  Span_{R_Y}\{v_\bullet\}$,
  for  any $d\ge1$. Thus $\quots{R_X}{\ca^{d}_M}$ is a f.g.   $R_Y$-module.  Fix some generators,  $R_X=\ca^{d}_M+Span_{R_Y}\{v_\bullet\}$.
 As before, we take $v_0=1$ and $v_{l>0}\in \cm$. Moreover, we choose the generators in a filtered way: $v_0,\!\dots,\!v_{j_1}$\! generate $\quots{R_X}{\ca_M}$,\! then
  $v_0,\!\dots,\!v_{j_1},\!\dots,\!v_{j_2}$\! generate $\quots{R_X}{\ca^2_M}$,\! and so on.

 As in Step 1.ii., the action $\cm\!\circlearrowright\! \quots{R_X}{\ca^{d}_M}$ is filtration-nilpotent, i.e.
  \!$\cm^N\!\sseteq  \!(y)\!+\!\ca^d_M$\!  for \!$N\!\gg\!1.$

  Expand  $c_{i}= \sum_{l }  v_l c_{i, l}$, where   $c_{i,l}\in R_Y$  are the new unknowns.

Repeat the procedure of Step 2.i, then equation \eqref{Eq.IFT.A.to.prove.current} becomes:
\beq\label{Eq.IFT.A.to.prove.expanded.form}
\Phi_X(f+g)-f-
 \sum_{i,l}  v_l\cdot
 \big[ c_{i , l}+c^g_{i ,  l}+Q_{i ,l}(\{c_{**}\})\big]\cdot \xi_j(f)\in \ca^d_M\cdot \RpX+T_{\cL^{(1)}}f.
\eeq

Here
$Q_{i , l}(z)\in (z)\cdot R_Y[[y\cdot z]][z]$, resp. $(z)\cdot R_Y\{y\cdot z\}[z]$, resp. $(z)\cdot R_Y\bl y\cdot z\br[z]$, rather than just
 $Q_{i , l}(z)\in (z)\cdot R_Y[[z]]$.

As $v_0=1$ and $v_{l>0}\in \cm$ one has here: $c^g_{i ,0}\in (y)$ and $Q_{i ,0}(z)\in (y)\cdot (z)$.

 We get the (finite) system of
 implicit function equations over the ring $R_Y$:
\beq\label{Eq.IFT.A.ift.to.resolve}
 c_{i , l}+c^g_{i ,  l}+Q_{i ,l}(\{c_{**}\})=0,\quad \forall\ i, l.
\eeq
Once they are satisfied, equation \eqref{Eq.IFT.A.to.prove.expanded.form} will imply equation \eqref{Eq.IFT.A.to.prove}.

Unlike the $\cR,\cK$-cases we cannot immediately apply the $IFT_\one$ of \S\ref{Sec.Preliminaries.Basic.Tools}.i.
 In fact the coefficients $c^g_{i ,  l}$ do not necessarily belong to the maximal ideal $(y)+\cm_\k\sset R_Y$.
   Neither can one assume that  $Q_{i , l}(z)$ is nilpotent in some sense, e.g. $Q_{i , l}(z)\not\in (f^\#(y)+\cm_\k)\cdot (z)$.

\eee

\

 In Step 3 we prove part a., there the system \eqref{Eq.IFT.A.ift.to.resolve} is replaced by a simpler system.
  To prove part b. (in Step 4) we introduce the indeterminate $t$, and  get an ``arc-solution". Then we extend this arc-solution to an ordinary solution.

\item  (The proof of part a. in the case $J\!=\!0$.) The system \eqref{Eq.IFT.A.ift.to.resolve} was obtained for the perturbation $g\!\in\! \ca^2\!\cdot\! \RpX$.
  To prove part a. we start from the perturbation
  $g\!\in\! (\cm_\k+(f)+\ca^2)\!\cdot\!\ca^2\!\cdot\! \RpX$. Observe:
  \beq
  (\cm_\k+(f)+\ca^2)\cdot\ca^2\cdot \RpX\sseteq (\cm_\k+(f)+\ca^2)\cdot\ca\cdot\cm\cdot  T_\cR f+ T_{\cL^{(1)}}f.
  \eeq
  Thus we fix some generators $\{\xi_i\}$ of the $R_X$-module $(\cm_\k+(f)+\ca^2)\cdot\ca\cdot T_\cR f$.
  Define $\Phi_X\in \cR$ by $x\to x+\sum c_i \xi_i(x)$. We get:
  \beq
\Phi_X(f+g)-f-\sum c_i \xi_i(f)\in    (\cm_\k+(f)+\ca^2)^2\cdot\ca^2\cdot\{c_*\}^2+\hspace{5cm}
  \eeq
  \[
\hspace{3cm}  +(\cm_\k+(f)+\ca^2) \cdot\ca^2\cdot((\cm_\k+(f)+\ca^2)\ca^2 \cdot)' \{c_*\}\sseteq
(\cm_\k+(f)+\ca^4)\cdot\ca\cdot\cm\cdot  T_\cR f+ T_{\cL^{(1)}}f.
  \]
 Expand $c_i=\sum c_{i,l}v_l$ and note that  the ideal $(\cm_\k+(f)+\ca^4)$ acts nilpotently on the filtration of $\{v_*\}$.

  We get again   equation   \eqref{Eq.IFT.A.ift.to.resolve}. However, now  $Q_{i , l}(z)\in (\cm_k+(y))+Nilp$, where
   the (non-linear) operator $Nilp$ is filtration-nilpotent on $\{v_*\}$.
  And this system is indeed resolvable by the $IFT_\one$.

Once it is resolved, we have the transition $\{f\}+\ca^2\cdot\RpX\stackrel{\cR}{\rightsquigarrow}\{f\}+\ca^d_M\cdot \RpX+T_{\cL^{(1)}}f.$
 By Step 1.iii. this implies the statement a.

\item  (The proof of part b. in the case $J=0.$)

\bee[\bf i.]
\item Extend the base ring $\k$ to $\k[[t]]$, resp. $\k\{t\}$, resp. $\k\bl t\br$, for an indeterminate $t$. Instead of $f+g$ take $f+tg$.
 Accordingly define the coordinate change $\Phi_X$  by $x\to x+t\sum   c_{il }v_l\xi_i(x)$.
   Repeat Step 2.ii up to equation \eqref{Eq.IFT.A.to.prove.expanded.form}:
\beq
\Phi_X(f+g)-f-
 t\sum_{i,l}  v_l\cdot
 \big[ c_{i , l}+c^g_{i ,  l}+t\cdot Q _{i ,l}(\{c_{**}\})\big] \xi_i(f) \ \in \
 \ca^d_M\cdot \RpX+ t\cdot \{c_{**}\}^2\cdot T_{\cL^{(1)}}f  .
\eeq

Then the  system \eqref{Eq.IFT.A.ift.to.resolve} becomes
$c_{i , l}+c^g_{i ,  l}+t\cdot Q_{i ,l}(\{c_{**}\})=0 .$

 Applying $IFT_\one$ over the ring $\k[[y,t]]$, resp. $\k\{y,t\}$, resp. $\k\bl y,t\br$, we get the solution $\{c_{i , l}(y,t)\}$.
 As $c^g_{i,0}\in (y)$ and $Q_{i,0}(z)\in (y)\cdot (z)$, we get: $c_{i,0}(y,t)\in (y).$

For this solution we define the coordinate change $\Phi_X\in \cR_t$ by $x\to x+t\sum   c_{i ,l}\cdot v_l\cdot    \xi_i(x)$.
 Observe: this coordinate change is unipotent.
  Indeed,  $\xi(x)\in \ca\cdot T_\cR (x)=\ca\sseteq\cm$ and $c_{i ,0}\in (y)$ and $v_{l>0} \in \cm$.

\quad  Altogether, we have proved: $\Phi_X(f+tg)\in \{f\}+\ca^{d}_M\cdot \RpX+ T_{\cL^{(1)}} f$, for an indeterminate $t$ and any $d\gg1$.
   Now part iii. of Step 1 ensures:
   $f+t\cdot g\stackrel{\cA^{(1)}_t}{\sim}f$.

\

An aside remark,  taking $t=\cm_\k+(y)$ gives $\{f\}+\ca^2\cdot (\cm_\k+(y))\cdot \RpX\sseteq \cA f$, which is a weaker version of part a.

\

  From now on we assume:  the field $\quots{\k}{\cm_\k}$ is   of characteristic zero or  algebraically closed. Moreover, \!we \!assume $\{Q_{**}(z)\}\!\not\sset (\cm_k\!+\!(y))\cdot (z).$
 (Otherwise \!the \!system \eqref{Eq.IFT.A.ift.to.resolve} is \!directly resolvable.)

   \item We prove:    $f+tg\stackrel{\cA}{\sim}f$ for an infinite set of values of $t\in \k$. We still have to resolve the
    system $ c_{i , l}+c^g_{i ,  l}+t\cdot Q_{i ,l}(\{c_{**}\})=0$, now for a fixed $t\in \k$. Here $ Q _{i ,l}(z)$ are power series in $z$,
     not polynomials. However, it is enough to resolve this system over the quotient ring $\quots{R_Y}{(y)+\cm_\k}$. Indeed, $Q _{i ,l}\in R_Y[[y\cdot z]][z]$
      (resp. analytic, resp. algebraic).  And
     if $\{c^0_{**}\}$ is a solution $mod(y+\cm_\k)$,
      then by shifting the variables, $c_{**}\rightsquigarrow c_{**}-c^0_{**}$, we get the system of $IFT_\one$-type,
      $c_{i ,l}=\tQ_{i ,l}(c_{**})$, where now $\tQ _{i ,l}(z)\in ((y)+\cm_\k)\cdot(z)$.
       (This later system is readily resolvable.)

Over the ring $\quots{R_Y}{(y)+\cm_\k}$ we get a polynomial system in variables $\{c_{i ,l}\}$, over the (infinite) field $\quots{\k}{\cm_\k}$.
This system defines a closed algebraic subscheme of a finite dimensional affine space $\A^N_{\quots{\k}{\cm_\k}}$.
  In Step 4.i we have constructed an arc on this subscheme,  $\{c_{i , l}(t)\}$.
 This arc is non-trivial (i.e. not a point) as $\{Q _{**}(z)\}\not\sset \cm_k+(y)$.
   Thus the set of $\quots{\k}{\cm_\k}$-points of this (closed, affine, algebraic) subscheme is of positive dimension.

 Moreover, the projection of this scheme onto the line $\A^1_t$ has image of positive dimension.
 Hence, as the field $\quots{\k}{\cm_\k}$ is infinite, the system is solvable for an infinite number of $t$-values.
 Altogether,  $f+tg\stackrel{\cA}{\sim}f$ for an infinite set of values of $t\in \k$.

\item  Finally we deduce: $f+tg\stackrel{\cA}{\sim}f$ for any $t\in \k$. As in Step 4.ii, it is enough to resolve the system over the ring
 $\quots{R_Y}{(y)+\cm_\k}$.

Through the whole proof we have used only the unipotent subgroup $\cA^{(1)}<\cA$, for the filtration $\cm^\bullet\!\sset\! R_X.$
\!\!Passing to the quotient $\quots{R_X}{\cm_\k\!+\!\cm^d}$ for $d\!\!\gg\!\!1,$
 we get the action on the finite dimensional vector space, $\cA^{(1)}\!\circlearrowright\! (\cm\!\cdot\! \quots{R_X}{\cm_\k\!+\!\cm^d})^{\oplus p}$. This action is
 algebraic and unipotent.

If the field $\quots{\k}{\cm_\k}$ is algebraically closed then the group-orbits are Zariski-closed, \S\ref{Sec.Preliminaries.Orbits.Unipotent.Group.Closed}.
The intersection of such an orbit with a line, $\cA^{(1)}(f)\cap\{f+tg\}_t$, is either finite, or contains the whole line.
 In our case this intersection is infinite, thus $\cA^{(1)}(f)\supseteq f+tg$ for each $t\in \k$.

The case ``$\quots{\k}{\cm_\k}$ is of characteristic zero" is done similarly, using lemma
 \ref{Thm.Orbit.Unipotent.Group.Closed}.
\eee

\

\item(The case $J\neq0.$ We prove part b., as part a. is simpler.)
 In Steps 3,4 we have constructed the map $\Phi_X:x\to x+\sum c_i\xi_i(x)$ satisfying $\Phi_Y\circ f\circ\Phi_X=f+g$.
  For $J\neq0$ this map is not a ring-automorphism. As in the $\cR,\cK$-cases we adjust it by higher order terms.
  To use the $jet_0$-assumption we expand further, $c_{i,l}=\sum \uy^\um c_{il\um}$,
   now $c_{il\um}\in \k$ are unknowns.
   Then we extend $\Phi_X$ to the automorphism $\Psi: x\to x+t\sum \uy^\um c_{il\um} v_l\xi_i(x)+t^2\varphi(\{c_{***}\}^2)$.
    Here the entries of $\varphi(\{c_{***}\}^2)$ are power series in $\{c_{***}\}$. Resolve the corresponding system with $t$-indeterminate.
     This is an arc solution on the variety in $\A^N_\k$.
     Arguing as in Step 4.ii we get: $f+tg\in \cA^{(1)}_t(\{f\}+ \ca^d_M\cdot \RpX).$
     Proceed as in Step 4.iii to get:  $f+ g\in \cA^{(1)} (\{f\}+\ca^d_M\cdot \RpX).$ This holds for any $d\ge1$.
      Therefore  $f+ g\in   \overline{\cA^{(1)} f}.$
  This proves the statement for the ring $\quots{\k[[x]]}{J}.$ For the ring $\quots{\k\bl x\br}{J}$ one uses the
    left-right Artin approximation, \S\ref{Sec.Preliminaries.Basic.Tools}.v.
\epr

\eee

\beR\label{Rem.A.IFT.unipotence}
The factor $\cm$ in the assumption $\ca^2\cdot \RpX\sseteq \ca\cdot \cm\cdot T_\cR f+y^2\cdot T_\cL f$ is used
  only once in this whole proof, in Step 4.i.
 It is needed to ensure that the coordinate change  $\Phi_X: x\to x+t\sum c_{i ,l}v_l   \xi_i(x)$ is unipotent.
  But this unipotence
  holds trivially if one assumes $\ca\sseteq\cm^2$. Thus (almost) the same proof leads to the modified part b.:
  {\em If $\ca\sseteq\cm^2$ and $\ca^2\cdot \RpX\sseteq \ca\cdot   T_\cR f+T_{\cL^{(1)}} f$ then
   $\{f\}+\ca^2\cdot \RpX+ T_{\cL^{(1)}} f\sseteq \cA f$.}
\eeR

\subsection{Examples and corollaries}
 Let $R_X$ be one of $\quots{\k[[x]]}{J}, \quots{\k\{x\}}{J}, \quots{\k\bl x\br}{J}$, where $\k$ is a local ring (e.g. any field),  \S\ref{Sec.Preliminaries.Rings.Germs}.ii.
  Suppose the field $\quots{\k}{\cm_\k}$ is algebraically closed or of characteristic zero.
  Take the assumptions \ref{Assumptions}.
\\
We use the filtration $\cm^\bullet\cdot\RpX$ and the corresponding tangent spaces
$T_{\cR^{(j)}},T_{\cL^{(j)}},T_{\cK^{(j)}},T_{\cA^{(j)}},$ see \S\ref{Sec.3.Filtrations.on.G.TG}.

 \subsubsection{The case $p=1$, i.e. the scalar valued functions} Let $f:X\to (\k^1,o)$ and suppose $\ca_\cA\sseteq \cm^2.$
 Then
\[\cA f\supseteq \{f\}+\ca^2_\cA+Span_\k(f^2,f^3,f^4,\dots).\]
\vspace{-0.6cm}\bpr Combine lema \ref{Thm.Annihilator.a_A.p=1.case} with remark \ref{Rem.A.IFT.unipotence}.\epr

\subsubsection{Maps whose $\cA$ and $\cK$ orbits are ``close"}
\bee[\bf\!i.\!]
\item Suppose $\ca_{\cA^{(o)}}\supseteq \cm\cdot (f)$, equivalently  $T_ {\cA^{(o)}} f=T_ {\cK^{(o)}} f$, equivalently $\ca_{\cA^{(o)}}=\ca_{\cK^{(o)}}$.
(This holds, e.g. for infinitesimally-stable maps, see Property \ref{Thm.Annihilators.Basic.Properties}.ii.)
 Then
 \beq
 \ca_{\cA^{(o)}}^2\RpX\sseteq \ca_{\cA^{(o)}}\cdot T_{\cR^{(0)}}f+\ca_{\cA^{(o)}}\cdot  T_{\cL^{(0)}}f\sseteq  \ca_{\cA^{(o)}}\cdot T_{\cR^{(0)}}f+ T_{\cL^{(1)}}f.
 \eeq
Part b. of Theorem \ref{Thm.IFT.A.case} gives (for $\ca=\ca_{\cA^{(0)}}$): $\{f\}+\ca _{\cK^{(o)}}^2\cdot\RpX\sseteq \cA f$.
 This is close to Mather's statement ``Stable maps are determined by their local algebras".

 Note: it is much simpler to compute/to control the ideal $\ca_{\cK^{(o)}}$ than $\ca _{\cA^{(o)}}$.
 \item
 As a particular case, suppose $J=0$ and take a stable map.
  It is a versal unfolding of its genotype, i.e. $f: (\k^{n+\tau-1},o)\to (\k^{p+\tau-1},o)$ with $f=(f_o+\sum v_j u_j,u)$,
   where $f_o:(\k^n,o)\to(\k^p,o)$ and $Span_\k(\{v_\bullet\})=\cm\cdot T^1_\cK f_o$.
   (Thus $\tau=dim T^1_\cK f_o$.)
  For $\k=\R,\C$ this is the classical Mather's theorem, \cite[Theorem 7.2]{Mond-Nuno}, for an arbitrary field see \cite{Kerner.Unfoldings}.

 We have  $\ca^f_{\cK^{(o)}}\!=\!\ca^{f_o} _{\cK^{(o)}} \!+\!(u)$, see example \ref{Ex.Annihilators.R.K.simple}.vi.
 Thus we get $\cA f\!\supset\!\{f\}\!+\!(\ca^{f_o}_{\cK^{(o)}} +(u))^2\!\cdot\! \RpX$.
This is much more useful than the classical bound  $\cA f\supset\{f\}+\cm^{p+\tau}\cdot \RpX$.
    As we have seen in  examples \ref{Ex.R.IFT.p=1}, \ref{Ex.K.IFT}, in most cases the ideal $(\ca^{f_o}_{\cK^{(o)}}+(u))^2$ is  larger than $\cm^{p+\tau}$.
     (Even the asymptotics is better.)
\item More generally, for $J=0$, suppose $\ca_{\cA^{(j)}}\supseteq \cm\cdot (f)^{j+1}$ for some $j\ge0$.  (Thus $\sqrt{\ca_{\cA^{(j)}}}=\sqrt{\ca_{\cK^{(j)}}}.$)
  Then
  \beq
  \ca_{\cA^{(j)}}^2\cdot \RpX\sseteq (\ca_{\cA^{(j)}}+(f)^{j+1})\cdot T_{\cR^{(j)}}f+\ca_{\cA^{(j)}}\cdot  T_{\cL^{(1)}}f\sseteq
    \ca_{\cA^{(j)}}\cdot T_{\cR^{(j)}}f+ T_{\cL^{(1)}}f.
\eeq
  Therefore  $\{f\}+\ca _{\cA^{(j)}}^2\cdot \RpX\sseteq \cA f$.

Yet more generally, it is enough to assume  $(f)^{j+1}\!\cdot\!\ca_{\cA^{(j)}}\!\cdot\! \RpX\!\sseteq  \!\ca_{\cA^{(j)}} T_{\cR^{(j)}}f\!+\! T_{\cL^{(1)}}f$ for some $j\!\ge\!0.$

\eee

\subsubsection{}\!\!\!\!The main assumption of theorem \ref{Thm.IFT.A.case} is: $\ca^2\!\cdot\! \RpX\!\sseteq\! \ca\!\cdot\!\cm\!\cdot \!T_\cR f+T_{\cL^{(1)}}f.$ This can be difficult to verify.
\bcor\label{Thm.IFT.A.corollary} (A weaker but simpler version)
\bee
\item   Suppose  $ \cb^2_Y\!\sseteq \!  \cb_Y\!\cdot\! \ca_{\cA^{(o)}}$  for  an  ideal $\cb_Y\!\sseteq\! f^\#(y).$
  Then   $\{f\}\!+\!  (\cb_Y\!+\!\ca_{\cR^{(o)}} )^4\!\cdot\!\RpX\!+\! T_{\cL^{(1)}} f\!\sseteq\! \cA f.$
 \\Moreover, if the field $\quots{\k}{\cm_\k}$ is algebraically closed or of characteristic zero then
  \\$\{f\}+  (\cb_Y+\ca_{\cR^{(o)}} )^2\cdot\RpX+ T_{\cL^{(1)}} f\sseteq \cA f.$
\item Suppose $ \cb^2_Y\sseteq   \cb_Y\cdot \ca_{\cA}$ for an ideal $  \cb_Y\sseteq f^\#(y)^2$. Suppose $\ca_\cR\sseteq \cm^2$.
\\Then   $\{f\}+  (\cb_Y+\ca_{\cR} )^4\cdot\RpX+ T_{\cL^{(1)}} f\sseteq \cA f.$
\\Moreover, if the field $\quots{\k}{\cm_\k}$ is algebraically closed or of characteristic zero then
\\$\{f\}+  (\cb_Y+\ca_{\cR} )^2\cdot\RpX+ T_{\cL^{(1)}} f\sseteq \cA f.$
\eee
\ecor
\bpr
\bee
\item
It is enough to observe: $(\cb_Y+\ca_{\cR^{(o)}} )^2\cdot\RpX\sseteq (\cb_Y+\ca_{\cR^{(o)}} )\cdot T_{\cR^{(0)}}+\cb_Y^2\cdot \RpX\sseteq
(\cb_Y+\ca_{\cR^{(o)}} )\cdot T_{\cR^{(0)}}+\cb_Y\cdot T_{\cL^{(0)}}$. Then we apply theorem \ref{Thm.IFT.A.case}.

\item The same proof, just use   remark \ref{Rem.A.IFT.unipotence}.
\epr
\eee

\bex\label{Ex.Bounds.for.A.IFT} (Assuming the field $\quots{\k}{\cm_\k}$ is algebraically closed or of characteristic zero)
\bee[\bf\!i.\!]
\item
 Take $\cb_Y=\ca_{\cA^{(o)}}\cap f^\#(y)$ to get  $\cA f\supset\{f\}+  (\ca_{\cA^{(o)}}\cap f^\#(y)+\ca_{\cR^{(0)}})^2\cdot \RpX+ T_{\cL^{(1)}} f$.

Suppose $\ca_\cR\sseteq \cm^2$.
 Then part two gives: $\cA f\supset\{f\}+  (\ca_{\cA}\cap f^\#(y)^2+\ca_{\cR})^2\cdot \RpX+ T_{\cL^{(1)}} f$.

More generally, if  $ \cb^2_Y\sseteq   \cb_Y\cdot \ca_{\cA^{(o)}}$ then the ideal $\cb_Y+\ca_{\cA^{(o)}}\cap f^\#(y)$ satisfies this condition as well.

    \item Suppose $\ca_{\cA}\supseteq \cm\cdot (f)^{d}$,    then we have $\ca_{\cA}\cap f^\#(y)\supseteq f^\#(y)^{d}$.
     Assuming $\ca_{\cR }\sset \cm^2$ we get
 \\$\{f\}+(  (f)^d+\ca_{\cR})^2\cdot \RpX+ T_{\cL^{(1)}} f\sseteq \cA f$.

We remark: $\sqrt{ (f)^d+\ca_{\cR^{(o)}}}=\sqrt{\ca_\cK}$, thus in this case the ($\cA$) instability and the ($\cK$) singularity loci coincide.
\eee
\eex
See \S\ref{Sec.A.Geometry} for the geometric interpretation.

\subsubsection{Computational cases}
\bee[\bf i.]
 \item Define the map $f: (\k^2,o) \to (\k^3,o)$ by $f(x_1,x_2)=(x_1,x_2^2,x_2^3+x_1^{k+1}x_2)$.
   Suppose $2,3,k+1\in \k^\times$.
  By the direct check:  $\ca_\cA=(x_1^k,x_2^2)\sset R_X$ and  $\ca^2_\cA\cdot \RpX\sseteq\cm\cdot \ca_\cA\cdot T_\cR f+T_\cL f$.
   Thus $\cA f\supseteq\{f\}+ (x_1^k,x_2^2)^2\cdot \RpX+  T_{\cL^{(1)}} f$.
   Compare this to the classical criterion: $\cm^{k+1}\sseteq \ca_\cA$ and therefore $\cA f\supset \{f\}+\cm^{2k+2}\cdot \RpX$,
 \cite[Proposition 6.1, pg. 197]{Mond-Nuno}.

 \item Define the map $f: (\k^2,o) \to (\k^3,o)$ by $f(x_1,x_2)=(x_1,x_2^2,x_1x_2^k)$, with $k$-odd.
  Suppose $2, k \in \k^\times$. Note that $f^\#(y)\supset (x_1x_2^{k-1})\sset R_X$.
  By the direct check:  $\ca_\cA=(x_2^k,x_1x_2^{k-1})\sset R_X$. This is a non-isolated instability.

  By the direct check:   $\ca^2_\cA\cdot \RpX\sseteq\cm\cdot \ca_\cA\cdot T_\cR f+T_\cL f$.
   Thus $\cA f\supseteq\{f\}+ (x_2^k,x_1x_2^{k-1})^2\cdot \RpX+T_{\cL^{(1)}} f$.

\eee

\subsection{The filtration criterion for $\cA$-orbits}\label{Sec.A.Filtration}
  Take  an ideal $  I\sset R_X$ and the filtration $M_\bullet:=I^\bullet\cdot \RpX.$
\bthe\label{Thm.Orbits.Aj.vs.TAj}  Fix some integers $1\le j<d$.
       Assume one of the following:
\bee[\bf i.]
\item either the conditions $jet_N(Exp)$, $jet_N(Ln)$ of \S\ref{Sec.Preliminaries.Exp.Ln} hold for all $N\ge1$. (Thus in particular $\k\supseteq\Q$.)
\item or $(f)\sseteq I$ and the conditions $jet_N(Exp)$, $jet_N(Ln)$ hold for $N=\lceil\frac{2d-1-ord(f)}{j}\rceil+1$.
\eee
Then: \quad\quad\quad $\overline{\cA^{(j)}f}\supseteq \{f\}+I^d\cdot \RpX$ \quad \quad \quad if and only if \quad \quad \quad  $\overline{T_{\cA^{(j)}} f}\supseteq I^d\cdot \RpX$.
\ethe
\bpr
{\bf The part $\Lleftarrow$.}
\bee[{\rm\bf  The case }$i$.]
\item Let $w_d\in M_d$ then   $w_d\equiv \xi_Y(y)|_f+\xi_X (f) \  mod\ M_{d'}$  for some $(\xi_Y,\xi_X)\in T_{\cA^{(j)}}f$.
 We can assume $d'\gg d$.
 Using the $jet_N(Exp)$-assumption, for $N=\lceil\frac{d-ord(f)}{j}\rceil+1$,
 we take the corresponding automorphisms $jet_N(e^{\xi_X})\in \cR^{(j)}$, $jet_N(e^{\xi_Y})\in \cL^{(j)}$.

   By remark \ref{Rem.Thom-Levine.without.Q}, applied to   $M_d$ and $\cA^{(j)}$,
  we have:
  $jet_N(e^{\xi_X})\cdot jet_N(e^{\xi_Y})f-f-w_d\in M_{d+j}$.
Therefore $\cA^{(j)} f+M_{d+j}\supseteq \{f\}+M_d$. Now apply the same argument to the pair  $(M_{d+j},\cA^{(j)})$,
 with $jet_N(Exp)$-assumption  for $N=\lceil\frac{d+j-ord(f)}{j}\rceil+1$. Iterate this to get
 $\cA^{(j)}f +{M_{d'}}\supseteq \{f\}+M_d$ for every $d'\gg1$. Hence $\{f\}+M_d\sseteq \overline{\cA^{(j)}f}$.

\item
 It is enough to prove: $\cA^{(j)}f+M_{k+1}\supset \{f\}+M_k$ for each $k\ge d$.
 For this we prove:
\beq\label{Eq.to.prove.1}
 \cA^{(kj+k-1)}f+ M_{kd+l+1} \supseteq \{f\}+M_{kd+l} ,\quad \text{ for each }\ k\ge1\   \text{ and   each }\  \ l\in 0,\dots,d-1.
\eeq
At each step we work modulo $M_{d'}$ with $d'\gg1$, thus we can replace $\overline{T_{\cA^{(\bullet)}}}f$ by  $ T_{\cA^{(\bullet)}} f$.

 By part one of lemma \ref{Thm.Artin.Rees.for.TA} we get   $ T_{\cA^{(kj+k-1)}} f \supseteq M_{kd} $ for each $k\ge1$.

Take any element $w\in M_{kd+l} $ and present it
 as $w=\xi_X (f)+\xi_Y(y)|_f\in T_{\cA^{(kj+k-1)}}f$.
 Take the corresponding automorphisms $jet_N(e^{\xi_X})$, $jet_N(e^{\xi_Y})$, ensured by the condition  $jet_N(Exp)$.

Apply remark \ref{Rem.Thom-Levine.without.Q}   to   $M_{kd+l}$ and $\cA^{(kj+k-1)}$ to get:
$  jet_N(e^{\xi_X})\cdot jet_N(e^{\xi_Y})f-f-w\in M_{(kd+l)+(kj+k-1)} .$
 This verifies the condition \eqref{Eq.to.prove.1}.

We should only justify the use of remark \ref{Rem.Thom-Levine.without.Q}, i.e. to verify the condition:
 \beq
 \frac{kd+l-ord(f)}{kj+k-1}+1\le \lceil\frac{2d-1-ord(f)}{j}\rceil+1 \quad \text{ for all $l=0,\dots,d-1$, and all $k,j\ge1$.}
 \eeq
  Here the maximum
   of the left hand side is achieved for $l=d-1$ and $k=1$. And then we get the equality.
\eee
{\bf The part $\Rrightarrow$.}
\bee[{\rm\bf  The case} $i$.]
\item
Take $w_d\in M_d$ then $f+w_d \equiv (\Phi_Y,\Phi_X)f\ mod\ M_{d'}$ for some $(\Phi_Y,\Phi_X)\in \cA^{(j)}$. We can assume  $d'\gg d$.
 Using the $jet_N$-assumptions we   approximate $(\Phi_Y,\Phi_X)$ by
  the elements $jet_N (e^{\xi_X})\in Aut^{(j)}_X$,  $jet_N (e^{\xi_Y})\in Aut^{(j)}_Y$ satisfying:
 \beq
\Phi_X-jet_N (e^{\xi_X})\in End^{(jN+1)}_\k(R_X),\quad\quad\quad
\Phi_Y-jet_N (e^{\xi_Y})\in End^{(jN+1)}_\k(R_Y).
 \eeq
Therefore $(jet_N (e^{\xi_X}),jet_N (e^{\xi_Y}))f-f-w_d\in M_{d'} +M_{jN+1+ord(f)}$. Here we can assume  $d'\ge jN+1+ord(f)\ge d+j$.
 Apply remark \ref{Rem.Thom-Levine.without.Q} to the pair $(M_d,\cA^{(j)})$
 to get: $(\xi_Y+\xi_X)f-w_d\in M_{d+j}$.
  Therefore $T_{\cA^{(j)}}f+M_{d+j}\supseteq M_d$.

  Now apply the same argument to the pair $(M_{d+j},\cA^{(j)})$, and so on.
 One gets $T_{\cA^{(j)}}f+M_{d'}\supseteq M_d$ for any $d'$. Hence the statement.

\item
Apply the argument of case i. to the pair  $(M_d,\cA^{(j)})$ and $N=\lceil\frac{ d -ord(f)}{j}\rceil+1$ to get: $T_{\cA^{(j)}}f+M_{d+j}\supseteq M_d$.
 Iterate this argument for all  pairs  $(M_{d+l},\cA^{(j)})$ and $N=\lceil\frac{ d+l -ord(f)}{j}\rceil+1$,
   with $l=1,\dots,d-1$. We reach  $T_{\cA^{(j)}}f+M_{2d}\supseteq M_d$, and for this we have
    used the $jet_N$-assumptions with $N=\lceil\frac{ d -ord(f)}{j}\rceil+1$, \dots,
     $\lceil\frac{ 2d -1 -ord(f)}{j}\rceil+1$.

And now observe:
\[
M_{2d}=I^d\cdot M_d\sseteq I^d(T_{\cA^{(j)}}f+M_{2d})\sseteq I^d\cdot T_{\cR^{(j)}}f+f^\#(\cb_j) T_{\cL^{(j)}}f+M_{3d}
\stackrel{(f)\sseteq I}{\sseteq}
T_{\cA^{(2j)}}f+M_{3d}.
\]
Iterate this to get $T_{\cA^{(j)}}f+M_{d'}\supseteq M_d$ for all $d'$. Hence the statement.
\epr
\eee

\subsubsection{}
The conclusion of theorem \ref{Thm.Orbits.Aj.vs.TAj} is for filtration closures, ``$\overline{\cA^{(j)}f}$ vs $\overline{T_{\cA^{(j)}} f}$".
 With additional assumptions we get the statement  ``$ {\cA^{(j)}f}$ vs $T_{\cA^{(j)}} f$".

  Let $R_X$ be one of $ {\k[[x]]} ,  {\k\{x\}} ,  {\k\bl x\br}  $, where the local ring $\k$ contains a  field,  \S\ref{Sec.Preliminaries.Rings.Germs}.ii.
  Fix some integers $1\!\le\! j\!<\!d$.
  Suppose either $char(\k)=0$ or
   $char (\k)\! \!>\!\!\lceil\frac{2d-1-ord(f)}{j}\rceil\!+\!1$.
    Take an ideal $(f)\sseteq I\sseteq\cm$.

\bcor\label{Thm.Orbits.Aj.vs.TAj.for.nice.rings}
\bee
\item  Suppose $f$ is $\cK$-finite.  If
   $ \cA^{(j)}f \supseteq \{f\}+I^d\cdot \RpX$  then $T_{\cA^{(j)}} f  \supseteq I^d\cdot \RpX$.
\item
Suppose $I\sseteq\sqrt{\ca_\cR+(\ca_\cA\cap f^\#(y))R_X}.$
 If  $T_{\cA^{(j)}} f  \supseteq I^d\cdot \RpX$ then
      $ \cA^{(j)}f \supseteq \{f\}+I^d\cdot \RpX.$
\eee
\ecor
Geometrically the condition $I\sseteq\sqrt{\ca_\cR+(\ca_\cA\cap f^\#(y))R_X}$ reads: $V(I)\supseteq f^{-1}(f(V(\ca_\cA)))\cap Crit(f)$.

For $\k\{x\}$ with $\k=\R,\C$, $I=\cm$, $j=1,$  this is Theorem 2.5 of \cite{Bruce.du-Plessis.Wall}.
\bpr W.l.o.g. we assume $I\neq0$, thus $V(\ca_\cA)\ssetneq (\k^n,o)$.
\bee
\item  By theorem \ref{Thm.Orbits.Aj.vs.TAj} we have $\overline{T_{\cA^{(j)}}}f\supseteq I^d\cdot \RpX$.
 Then ${T_{\cA^{(j)}}}f\supseteq I^d\cdot \RpX$, by part two of lemma \ref{Thm.Filtration.Closures.TA}.

\item
By    theorem \ref{Thm.Orbits.Aj.vs.TAj} it is enough to show: $ \cA^{(j)}f \supseteq \{f\}+I^N\cdot \RpX$ for some $N\gg1$.
 Part 2 of corollary  \ref{Thm.IFT.A.corollary} gives (for $\cb_Y=\ca_\cA\cap f^\sharp(y)^2$):
      $  \{f\}+I^{N}\cdot  \RpX\sseteq \{f\}+(\ca_\cR+\ca_\cA\cap f^\sharp(y)^2)^4\cdot \RpX\sset \cA f$.

    Moreover, (for $N\gg1$) in the proof of theorem  \ref{Thm.IFT.A.case}  we get the coefficients $c_{**}\in I^j.$
     Therefore the needed coordinate change is an $\cR^{(j)}$-element. Similarly,
     the involved $\cL$-elements are in fact $\cL^{(j)}$ elements. Therefore
   $ \cA^{(j)} f \supseteq \{f\}+I^N\cdot \RpX$.
\epr
\eee

\bex Let $R_X$ be one of $\k[[x]]$, $\k\{x\}$, $\k\bl x\br.$
 Take the filtration $I^\bullet\cdot \RpX$   and the corresponding pairs $(T_{\cA^{(j)}},\cA^{(j)})$ of \S\ref{Sec.3.Filtrations.on.G.TG}.
  Assume $(f)\sseteq I$ and    $\sqrt{\ca_\cK}\supseteq I.$
\bee[\bf i.]

\item Suppose  $T_{\cA^{(0)}}f\supseteq I^{d_\cA}\cdot \RpX$ and  $T_{\cK^{(0)}}f\supseteq I^{d_\cK}\cdot \RpX.$
    Assume  either $char(\k)=0$ or $char(\k)>2(d_\cA+d_\cK)-ord(f).$ Then $\cA^{(1)}f\supseteq \{f\}+I^{d_\cA+d_\cK}\cdot \RpX.$
\bpr
 By part 2 of lemma \ref{Thm.Artin.Rees.for.TA} we get:  $T_{\cA^{(1)}}f\supseteq I^{d_\cA+d_\cK}\cdot \RpX$.
   Thus $\sqrt{\ca_\cA}\supseteq I\supseteq(f).$ Hence $\sqrt{\ca_\cR+(\ca_\cA\cap f^\sharp(y))R_X}=\sqrt{\ca_\cR+(f)}=\sqrt{\ca_\cK}\supseteq I.$
    (See example \ref{Ex.Annihilators.R.K.simple}.iii.)
\\Finally, apply    Part 2 of corollary \ref{Thm.Orbits.Aj.vs.TAj.for.nice.rings}.
 \epr

\item
We compare the conclusion $\cA^{(1)}f\supseteq \{f\}+I^{d_\cA+d_\cK}\cdot \RpX$  to the classical results (for $\k=\R,\C$):
\bei
\item
  \cite[Theorems  0.2 and 2.1]{Gaffney-du Plessis}:  Assume  $f$ is $\cK$-finite.
   If $\cm\cdot T_\cR f+  T_{\cL^{(0)}}f\supseteq I^d\cdot \RpX$ then
  $\cA f\supseteq \{f\}+I^{2d+2}\cdot \RpX.$
\item \cite[\!Theorem 6.2]{Mond-Nuno}:  If $\cm\cdot T_\cR f+  T_{\cL^{(0)}} f\!\supseteq\! \cm^d\cdot \RpX$ then
  $\cA f\!\supseteq\! \{f\}+\cm^{2d}\cdot \RpX.$
  \item
 \cite[\!\!Corollary \!6.3]{Mond-Nuno}: \!If $T_{\cA}f\!\!\supseteq\! \cm^{d_\cA}\!\cdot\! \RpX$ and  $T_{\cK}f\!\!\supseteq\! \cm^{d_\cK}\!\cdot\! \RpX$ then
  \!$\cA \!f\!\!\supseteq\! \{f\}\!+\!\cm^{d_\cA+d_\cK+1}\!\cdot\! \RpX.$
\eei
For most types of maps (if one is far from a stable map) one has: if $I^d\cdot\RpX\sseteq T_\cA f$ then
 $\cm^d\cdot\RpX\sseteq T_{\cA^{(0)}} f\subsetneq T_{\cK^{(0)}} f.$  Thus our bounds extend (and slightly strengthen) the previous results.
\eee
\eex

\subsection{Geometric characterization of $\cA$-determinacy}\label{Sec.A.Geometry}
 Given a map $f\in \Maps,$ the support of $T^1_\cA f$ is the instability locus, $V(\ca_{\cA})\sseteq Crit(f)\sset X$. The ideal $\ca_{\cA }\cap f^\#(y)\sset R_Y$ defines the image
 $f(V(\ca_{\cA}))\sset (\k^p,o).$ The ideal
$R_X(\ca_{\cA }\cap f^\#(y))+\ca_{\cR}$ defines the preimage,
  $f^{-1}(f(V(\ca_{\cA })))\cap Crit(f)$. This preimage contains the instability locus.
 If $f$ is of finite singularity type (i.e. the restriction $f:Crit\to (\k^p,o)$ is finite), then the loci $V(\ca_\cA),$
  $f^{-1}(f(V(\ca_{\cA })))\cap Crit(f)$ have the same dimension.
 Moreover, the two loci often coincide.

Corollary \ref{Thm.IFT.A.corollary} extends the classical determinacy criterion to the generality of $\quots{\k[[x]]}{J},\! \quots{\k\{x\}}{J},\! \quots{\k\bl x\br}{J}$.
 Let $\ca\sset R_X$ be the defining ideal of the locus  $f^{-1}(f(V(\ca_{\cA })))\cap Crit(f).$
 \bcor Assume $\ca_\cA\!\sseteq\!\cm.$
\bee
\item The map $f\!:X\!\to \!(\k^p,o)$ is   $\cA$-determined by its finite jet on the  locus  $f^{-1}(f(V(\ca_{\cA })))\!\cap\! Crit(f).$ Namely,   $\cA f\!\supset\! \{f\}\!+\!\ca^d\!\cdot\! \RpX,$ for $d\!\gg\!1,$ where $\ca\!\sset\! R_X$ is the defining ideal of $f^{-1}(f(V(\ca_{\cA })))\cap Crit(f).$
\item
 If $f$ is $\cA$-finite (i.e. $\sqrt{\ca}=\cm$), then $f$ is $\cA$-finitely determined.

 \item
  Suppose $f$ is $\cK$-finite. Take the filtration $\cm^\bullet\cdot\RpX.$
Suppose    the local ring $\k$ contains a  field, and   $char(\k)=0$ or $char(\k)\gg1$. If $\cA^{(1)} f\supseteq \{f\}+\ca^d\cdot \RpX$ then $T_\cA f\supseteq (\ca_\cR+\ca^d)\cdot \RpX.$
 In particular,  $Supp( T^1_\cA f)\sseteq V(\ca)\sset X$, i.e. $f$ is infinitesimally stable over $X\smin V(\ca)$.
\eee
\ecor

For complex analytic maps, $f:(\C^n,o)\to (\C^p,o),$ with isolated instabilities (i.e. $V(\ca_\cA)\!=\!o\in (\C^n,o)$)   this gives the classical finite determinacy criterion of Mather-Gaffney.

 In the case $\sqrt{\ca_\cA}\ssetneq \cm$, i.e. $f$ has a non-isolated instability, this result seems to be new even for $\C\{x\}$.

\bpr Parts 1 and 2 follow by Corollary \ref{Thm.IFT.A.corollary}.
Part 3 follows by Corollary \ref{Thm.Orbits.Aj.vs.TAj.for.nice.rings}.
 \epr


\vspace{-0.1cm}
\begin{thebibliography}{\!\!\!\!\!99\!\!}
\vspace{-0.1cm}

\bibitem[Abhyankar]{Abhyankar} Sh.-Sh. Abhyankar, {\em  Local analytic geometry.} Singapore: World Scientific. xv, 488 p. (2001).



\bibitem[A.G.L.V.]{AGLV} V. I. Arnold, V.V. Goryunov, O.V. Lyashko, V.A. Vasil'ev, {\em Singularity theory. I.}  Reprint of the original English edition from the series Encyclopaedia of Mathematical Sciences . Springer-Verlag, Berlin, 1998. iv+245 pp


\bibitem[Artin.68]{Artin.68} M. Artin, {\em On the solutions of analytic equations.} Invent. Math. 5 (1968), 277--291.
\bibitem[Artin.69]{Artin} M. Artin, {\em  Algebraic approximation of structures over complete local rings}, Publ. Math. IHES, 36, (1969), 23-58. 


\bibitem[B.K.16]{B.K.motor}
G. Belitski, D. Kerner,  {\em Group actions on filtered modules and finite determinacy. Finding large submodules
in the orbit by linearization}, C. R. Math. Acad. Sci. Soc. R. Can. 38 (2016),
no. 4, 113--153.

\bibitem[B.G.K.22]{BGK.20} A.-F. Boix, G.-M. Greuel, D. Kerner, {\em  Pairs of Lie-type and large orbits of group actions on filtered modules.
 (A characteristic-free approach to finite determinacy.)},  Math. Z. 301 (2022), no. 3, 2415--2463.

\bibitem[Bo.Da.Re.21]{BDR} M. Borovoi, C. Daw, and J. Ren, {\em Conjugation of semisimple subgroups over real number fields of bounded degree},
Proc. Amer. Math. Soc. 149 (2021), no. 12, 4973--4984.

\bibitem[Bou.Gr.Ma.11]{Boubakri.Gre.Mark2011} Y. Boubakri, G.-M. Greuel, Th. Markwig, {\em
Normal forms of hypersurface singularities in positive characteristic.}
Mosc. Math. J. 11, No. 4, 657--683 (2011).

\bibitem[Bou.Gr.Ma.12]{Boubakri.Gre.Mark} Y. Boubakri, G.-M. Greuel, Th. Markwig, {\em Invariants of hypersurface singularities in positive characteristic.}
Rev. Mat. Complut. 25 (2012), no. 1, 61--85.


\bibitem[Bourbaki.CA]{Bourbaki.CA} N. Bourbaki, {\em Alg\`{e}bre commutative},  Fasc. XXVIII, Chap. III, \S4, No. 5 (1961).

\bibitem[Bourbaki.Lie]{Bourbaki.Lie} N. Bourbaki, {\em Elements of Mathematics. Lie groups and Lie algebras. Chapters 1--3},
Hermann, Paris; Addison-Wesley Publishing Co., Reading, Mass. 1975, xxviii+450.


\bibitem[B.dP.W.87]{Bruce.du-Plessis.Wall} J.W. Bruce, A.A. du Plessis, C.T.C. Wall, {\em Determinacy and unipotency.} Invent. Math. 88 (1987), no. 3, 521--554.




\bibitem[Bru.Rob.88]{Bruce-Roberts} J. W. Bruce, R. M. Roberts, {\em Critical points of functions on analytic varieties.}
Topology 27 (1988), no. 1, 57--90.

\bibitem[Bru.Rua.Sai.928]{Bruce.Ruas.Saia} J.W.Bruce, M.A.S.Ruas, M.J.Saia,{\em A note on determinacy.}
 Proc.Amer.Math.Soc. 115 (1992), no. 3, 865--871.


\bibitem[Cut.Sri.97]{Cutkosky-Srinivasan.1997} S.D.Cutkosky, H.Srinivasan, {\em Equivalence and finite determinancy of mappings.}
J. Algebra 188 (1997), no. 1, 16--57.

\bibitem[Damon.84]{Damon84}
J. Damon, {\em The unfolding and determinacy theorems for subgroups of $\cA$ and $\cK$.} Mem. Amer. Math. Soc. 50 (1984), no. 306, x+88 pp.


\bibitem[Damon.88]{Damon88}
J.N. Damon, {\em Topological triviality and versality for subgroups of $\cA$ and $\cK$.} Mem. Amer. Math. Soc. 75 (1988), no. 389, x+106 pp.



\bibitem[Damon.91]{Damon}\!J. Damon,{\em $\cA$-equivalence \!and \!the \!equivalence of sections of images and discriminants.}
Singularity theory and its applications.Pt.I:Geometric aspects of singularities, Proc.Symp.,Warwick 1988-89, Lect.Notes Math.1462,93--121 (1991).

\bibitem[Damon.92]{Damon92} J. Damon, {\em Topological triviality and versality for subgroups of $\cA$ and $\cK$.
II. Sufficient conditions and applications.} Nonlinearity 5 (1992), no. 2, 373--412.


\bibitem[Denef-Lipshitz]{Denef-Lipshitz} J. Denef, L. Lipshitz, {\em Ultraproducts and approximation in local rings. II.}
Math. Ann. 253 (1980), no. 1, 1--28.


\bibitem[Domitrz-Rieger.09]{Domitrz-Rieger}  W. Domitrz,    J. H. Rieger, {\em Volume preserving subgroups of $\cA$ and $\cK$ and singularities in unimodular geometry.} Math. Ann. 345 (2009), no. 4, 783--817.


\bibitem[du Plessis.80]{du Plessis} A. du Plessis, {\em On the determinacy of smooth map-germs.} Invent. Math. 58 (1980), no. 2, 107--160.


\bibitem[Ebeling]{Ebeling} W.Ebeling, {\em Functions of several complex variables and their singularities.}
  Translated from the 2001 German original by Philip G. Spain. Graduate Studies in Mathematics, 83.
   AMS, Providence, RI, 2007. xviii+312 pp.


\bibitem[Eisenbud]{Eisenbud-book} D. Eisenbud, {\em Commutative algebra. With a view toward algebraic geometry.}
 Graduate Texts in Mathematics, 150. Springer-Verlag, New York, 1995.


\bibitem[F.-S. R.]{Ferrer Santos-Rittatore} W. R. Ferrer Santos,  A. Rittatore,
 {\em Actions and invariants of algebraic groups.} Second edition. Monographs and Research Notes in Mathematics. CRC Press, Boca Raton, FL, 2017. xx+459 pp.

\bibitem[Gabrielov.73]{Gabrielov} A.M.Gabri\`{e}lov,{\em Formal relations among analytic functions.}
Izv.Akad.Nauk SSSR Ser.Mat.37 (1973), 1056--1090.

\bibitem[Gaffney.79]{Gaffney}\!\!T.Gaffney,{\em A note on the order of determination of a finitely determined germ.}
Invent.Math.52(1979), no.2,127--130.




\bibitem[Gaffney-du Plessis.82]{Gaffney-du Plessis} T. Gaffney, A. du Plessis, {\em More on the determinacy of smooth map-germs.}
  Invent. Math. 66 (1982), no. 1, 137--163.


\bibitem[Greuel-Kr\"oning.90]{Gre.Kron.90} G.-M. Greuel, H. Kr{\"o}ning, \emph{Simple singularities in positive characteristic}, Math. Z. \textbf{203} (1990),
 339--354.


\bibitem[Gr.Lo.Sh]{Gr.Lo.Sh} G.-M. Greuel, C. Lossen, E. Shustin, {\em Introduction to singularities and deformations.}
Springer Monographs in Mathematics. Springer, Berlin, 2007. xii+471 pp.


\bibitem[Greuel-Nguyen.16]{Greuel-Nguyen.2016} G.-M. Greuel, H. D. Nguyen, {\em Right simple singularities in positive characteristic.}
J. Reine Angew. Math. 712 (2016), 81--106.



\bibitem[Greuel.18]{Greuel-review.sings.char.positive} G. M. Greuel, {\em Singularities in positive characteristic: equisingularity,
 classification, determinacy}. Singularities, algebraic geometry, commutative algebra, and related topics, 37--53, Springer, Cham, 2018


\bibitem[Greuel-Pham.19]{Greuel-Pham.2017} G.-M. Greuel, T.H. Pham, {\em Finite determinacy of matrices and ideals in arbitrary characteristics},
 J. Algebra 530 (2019), 195--214.




\bibitem[Houston-Wik Atique.13]{Houston.Wik.Atique}    K. Houston, R. Wik Atique, {\em    $\cA$-classification of map-germs via $_V\!\cK$-equivalence},
 arXiv:1302.1104.

\bibitem[Jacobson]{Jacobson} N.Jacobson, {\em Lie algebras.} Republication of the 1962 original. Dover Publications,Inc.,New York,1979. ix+331 pp.


\bibitem[Jahnel.00]{Jahnel} J. Jahnel, {\em The Brauer-Severi variety associated with a central simple algebra: a survey.}
Preprint: Linear Algebraic Groups and Related Structures, no. 52 (2000),
https://www.math.uni-bielefeld.de/lag/man/052.pdf.

\bibitem[Kerner.22]{Kerner.Unfoldings} D. Kerner, {\em Unfoldings of maps, the first results on stable maps, and results of
  Mather-Yau/Gaffney-Hauser type in arbitrary characteristic},   arXiv:2209.05071.

\bibitem[Kerner.23]{Kerner.LRAP} D. Kerner, {\em Artin approximation for left-right equivalence of map-germs},
 arXiv:2310.01521.

\bibitem[Kucharz.86]{Kucharz} W. Kucharz, {\em Power series and smooth functions equivalent to a polynomial.} Proc. Amer. Math. Soc. 98 (1986), no. 3, 527--533.


\bibitem[Lafon.65]{Lafon.65} J.-P. Lafon, {\em S\'eries formelles alg\'ebriques,}
C. R. Acad. Sci. Paris 260 (1965), 3238–-3241.

\bibitem[Lafon.67]{Lafon} J.-P. Lafon, {\em Anneaux hens\'{e}liens et th\'{e}or\`{e}me de pr\'{e}paration.}
C.R.Acad. Sci. Paris S\'{e}r. A-B264 (1967),A1161-A1162.


\bibitem[L-P.N-B.O-O.T.]{L-P.N-B.O-O.T.}
B. K. Lima-Pereira,  J.J. Nu\~{n}o-Ballesteros, B, Or\'{e}fice-Okamoto, J.N. Tomazella, {\em
The relative Bruce-Roberts number of a function on a hypersurface.}
Proc. Edinb. Math. Soc. (2) 64 (2021), no. 3, 662--674.

\bibitem[Martinet.76]{Martinet.1976} J. Martinet, {\em D\'{e}ploiements versels des applications diff\'{e}rentiables et classification des applications stables.}
  In Singularit\'{e}s d`applications diff\'{e}rentiables (S\'{e}m., Plans-sur-Bex, 1975), volume 535 of Lecture Notes in Math., pages 1--44.
 Springer, Berlin, 1976.


\bibitem[Martinet]{Martinet.1982}\!\!J.Martinet,{\em  Singularities \!of \!smooth \!functions \!and \!maps}, \!Lecture \!Note \!Series58,Cambridge University Press,1982,256 pp.

\bibitem[Mather.68-69]{Mather1968} J.N. Mather,
{\em Stability of $C^\infty$-mappings. I. The division theorem.} Ann. of Math. (2) 87, 1968, 89--104.
\\{\em Stability of $C^\infty$-mappings. II. Infinitesimal stability implies stability.} Ann. of Math. (2) 89, 1969, 254--291.
\\{\em Stability of $C^\infty$-mappings.III.Finitely determined map-germs.}Publ. Inst. Hautes \'{E}tudes Sci.Publ.Math.No.35, 1968, 279--308.
\\{\em Stability of $C^\infty$-mappings. IV. Classification of stable germs by R-algebras.} Publ. Inst. Hautes \'{E}tudes Sci. Publ. Math. No. 37, 1969, 223--248.


\bibitem[Matsumura]{Matsumura} H. Matsumura, {\em Commutative ring theory,}
Cambridge Studies in Advanced Mathematics, 8. Cambridge etc.: Cambridge University Press. XIII, 320 p.



\bibitem[Mond-Montaldi.91]{Mond-Montaldi} D. Mond, J. Montaldi, {\em Deformations of maps on complete intersections, Damon's $\cK_V$-equivalence and bifurcations.}   `Singularities in geometry and topology', Cambridge: Cambridge University Press. Lond. Math. Soc. Lect. Note Ser. 201, 263-284 (1994).

\bibitem[Mon. Nu\~{n}.-Bal.]{Mond-Nuno} D. Mond, J.-J. Nu\~{n}o-Ballesteros, {\em Singularities of Mappings.
The Local Behaviour of Smooth and Complex Analytic Mappings},
Grundlehren der mathematischen Wissenschaften, 357. Springer, Cham, [2020], 2020. xv+567 pp


\bibitem[Moret-Bailly.21]{Moret-Bailly}     L. Moret-Bailly, {\em     A Henselian preparation theorem},   Israel Journal of Mathematics, 257, 519-531 (2023).



\bibitem[Nguyen.19]{Nguyen} H. D. Nguyen, {\em Right unimodal and bimodal singularities in positive characteristic.} Int. Math. Res. Not. IMRN 2019, no. 6,
 1612--1641.

\bibitem[Pellikaan.88]{Pellikaan88} R. Pellikaan, {\em
Finite determinacy of functions with nonisolated singularities.}
Proc. London Math. Soc. (3) 57 (1988), no. 2, 357--382.

\bibitem[Jo-P\'{e}.Nu-Ba.09]{Perez-Nuno} V.H. Jorge P\'{e}rez, J.J. Nu\~{n}o-Ballesteros,
 {\em  Finite determinacy and Whitney equisingularity of map germs from $\C^n$ to $\C^{2n-1}$.}
Manuscripta Math. 128 (2009), no. 3, 389--410.



\bibitem[Pfister-Popescu.75]{Pfister-Popescu} G.  Pfister, D. Popescu, {\em Die strenge Approximationseigenschaft lokaler Ringe.}
Invent. Math. 30 (1975), no. 2, 145--174.


\bibitem[Pham.16]{Pham.Phd.Thesis}  T. H. Pham, {\em On Finite Determinacy of Hypersurface Singularities and Matrices in Arbitrary Characteristic},
 PhD thesis, 2016.

\bibitem[Popescu.00]{Popescu}D. Popescu,
{\em  Commutative \!Rings \!and \!Algebras. \!Artin \!approximation},   \!Handbook \!of \!Algebra, \!Vol 2, pg. 321-356, 2000.

\bibitem[Rond.18]{Rond}   G. Rond, {\em Artin approximation},  J. Singul. 17 (2018), 108--192.

\bibitem[Ruas.83]{Ruas} M. A. S. Ruas, {\em On the finite $C^l$-determinacy and applications}, PhD thesis, 1983, S\~{a}o-Carlos.


\bibitem[Ruiz]{Ruiz}  J. M. Ruiz, {\em The basic theory of power series.}
Advanced Lectures in Mathematics. Friedr. Vieweg \& Sohn, Braunschweig, 1993. x+134 pp.



\bibitem[Serre]{Serre.Lie}  J.-P. Serre, {\em Lie algebras and Lie groups.} 1964 lectures, given at Harvard University. 2nd ed.
Lecture Notes in Mathematics. 1500. Berlin etc.: Springer-Verlag. vii, 168 p. (1992).

\bibitem[Serre.02]{Serre} J.-P. Serre, {\em Galois cohomology,} Corrected reprint of the 1997 English edition,
Springer Monographs in Mathematics, Springer-Verlag, Berlin, 2002.


\bibitem[Shiota.98]{Shiota.1998} M. Shiota, {\em Relation between equivalence relations of maps and functions.}
  Real analytic and algebraic singularities (Nagoya/Sapporo/Hachioji, 1996), 114--144, Pitman Res. Notes Math. Ser., 381, Longman, Harlow, 1998.

\bibitem[Shiota.10]{Shiota.2010} M. Shiota, {\em Analytic and Nash equivalence relations of Nash maps.} Bull.Lond.Math.Soc. 42 (2010), no. 6, 1055--1064.







\bibitem[Stacks]{Stacks} {\em The Stacks project}, https://stacks.math.columbia.edu/

 
\bibitem[Sansuc.81]{Sansuc} J.-J. Sansuc, {\em Groupe de Brauer et arithm\'etique des groupes alg\'ebriques lin\'eaires sur
un corps de nombres}, J. Reine Angew. Math. 327 (1981), 12--80.




\bibitem[Sun-Wilson.01]{Sun-Wilson} B. Sun, L. C. Wilson, Leslie {\em Determinacy of smooth germs with real isolated line singularities.}
Proc. Amer. Math. Soc. 129 (2001), no. 9, 2789--2797.

\bibitem[Tougeron.66]{Tougeron1966} J.C. Tougeron, {\em Une g\'{e}n\'{e}ralisation du th\'{e}or\`{e}me
 des fonctions implicites.} C. R. Acad. Sci. Paris S\'{e}�r. A--B 262 1966 A487--A489


\bibitem[Wall.81]{Wall-1981} C.T.C. Wall, {\em Finite determinacy of smooth map-germs.} Bull. London
Math. Soc. 13 (1981), no. 6, 481--539.

\bibitem[Wall.95]{Wall-1995} C. T. C. Wall, {\em Classification and stability of singularities of smooth maps.}
  Singularity theory (Trieste, 1991), 920-952, World Sci. Publ., River Edge, NJ, 1995.

\bibitem[Wilson.81]{Wilson} L. C. Wilson, {\em Infinitely determined map germs.} Canadian J. Math. 33 (1981), no. 3, 671--684.

\end{thebibliography}
\end{document}